\documentclass[reqno,english]{amsart}

\makeatletter
\def\@currentlabel{2.1}\label{e:dispaa}
\def\@currentlabel{2.21}\label{e:dispau}
\def\@currentlabel{2.22}\label{e:dispav}
\def\@currentlabel{2.23}\label{e:dispaw}
\def\@currentlabel{2.24}\label{e:dispax}
\def\theequation{\thesection.\@arabic\c@equation}
\makeatother

\newcommand{\ttt}{\tilde }
\newcommand{\LL}{{\tt L}  }
\newcommand{\TT}{{\mathcal T}  }

\newcommand{\QQ}{{\mathcal Q}  }
\newcommand{\nn}{ {\nabla}  }

\newcommand{\by}{ {\bf y } }
\newcommand{\ps}{ {\psi^\perp } }
\newcommand{\py}{ {\tt y } }
\newcommand{\pp}{ {\partial} }
\newcommand{\ww}{{\tt w}  }
\newcommand{\HH}{{{\mathbb H}}  }

\newcommand{\A}{\alpha }
\newcommand{\vp}{\varphi}

\newcommand{\OO}{{\mathcal O}}

\newcommand{\JJ}{{\mathcal J}}

\newcommand{\NN}{ {\mathcal N}}

\newcommand{\R} {\mathbb R}
\newcommand{\cuad}{{\sqcap\kern-.68em\sqcup}}

\newcommand{\BB}{{\tt B}}
\newcommand{\NNN}{{\tt N}}

\newcommand{\KK}{{\mathcal K}}
\newcommand{\dist}{{\rm dist}\, }
\newcommand{\foral}{\quad\mbox{for all}\quad}
\newcommand{\ve}{\varepsilon}

\newcommand{\be}{\begin{equation}}
\newcommand{\ee}{\end{equation}}

\newcommand{\la}{\lambda}

\newcommand{\equ}[1]{(\ref{#1})}

\renewcommand{\theequation}{\thesection.\arabic{equation}}
 
 \newtheorem{lemma}{Lemma}[section]

\newtheorem{teo}{Theorem}
\newtheorem{prop}{Proposition}[section]

\newtheorem{remark}{Remark}[section]
\newcommand{\bremark}{\begin{remark} \em}
\newcommand{\eremark}{\end{remark} }

\title[ the Allen Cahn equation and Minimal Surfaces in $\R^3$]
{Entire Solutions of the Allen-Cahn equation  and Complete Embedded  Minimal Surfaces of Finite Total Curvature in $\R^3$}

\author{Manuel del Pino}
\address{ 
Departamento de
Ingenier\'{\i}a  Matem\'atica and CMM, Universidad de Chile,
Casilla 170 Correo 3, Santiago,
Chile.}
\email{delpino@dim.uchile.cl}

\author{Michal Kowalczyk}
\address{
Departamento de
Ingenier\'{\i}a  Matem\'atica and CMM, Universidad de Chile,
Casilla 170 Correo 3, Santiago,
Chile.}
\email{kowalczy@dim.uchile.cl}

\author{Juncheng Wei}
\address{
Department of Mathematics, Chinese University of Hong Kong, Shatin, Hong Kong
} \email{wei@math.cuhk.edu.hk}

\begin{document}


\begin{abstract}
We consider  minimal
surfaces $M$ which are complete, embedded and have finite total curvature in $\R^3$, and
bounded, entire solutions with finite Morse index of  the Allen-Cahn equation $\Delta u +
f(u) = 0 \  \hbox{in} \ \R^3 $. Here $f=-W'$ with $W$  bistable and
balanced, for instance $W(u) =\frac 14 (1-u^2)^2$. We assume that $M$ has $m\ge 2$ ends, and  additionally that $M$
is non-degenerate, in the sense that its bounded Jacobi fields
are all originated from rigid motions (this is known for instance for a Catenoid and for the Costa-Hoffman-Meeks surface of any genus).
We prove that for any small
$\alpha >0$, the Allen-Cahn equation has  a family of bounded solutions depending on $m-1$ parameters distinct from rigid motions, whose level sets
 are embedded surfaces  lying close to  the blown-up
surface  $M_\alpha := \alpha^{-1} M$,  with ends possibly diverging
logarithmically from $M_\A$. We prove that
these solutions are $L^\infty$-{\em non-degenerate} up to rigid motions, and find that their Morse index
coincides with the index of the minimal surface.   Our construction suggests
parallels of De Giorgi conjecture for general bounded solutions of
finite Morse index.
\end{abstract}

\date{}\maketitle
\tableofcontents

\setcounter{equation}{0}
\section{ \emph{Introduction and main results}}

\subsection{The Allen-Cahn equation and minimal surfaces}
The Allen-Cahn equation in $\R^N$ is the semilinear elliptic problem
\be \Delta u\, + \, f(u)\, =\, 0 \quad \hbox{in } \R^N\, ,
\label{ac}\ee where $f(s) = -W'(s)$ and $W$ is a ``double-well
potential'', bi-stable and balanced, namely \be W(s) > 0\hbox{ if }\
s\ne 1, -1, \quad W(1) =0=W(-1), \quad W''(\pm 1)\, =\, f'(\pm
1)\,=:\, \sigma_{\pm}^2 >0. \label{acn}\ee A typical example of such
a nonlinearity is \be f(u) \,= \, (1-u^2)u\,  \quad \hbox{ for }
W(u)\, =\, \frac 14 (1-u^2)^2, \label{ace}\ee
while we will not make 
use of the special symmetries enjoyed by this example.

Equation \equ{ac} is a prototype for the continuous modeling of phase transition phenomena.
 Let us consider the energy in a subregion region $\Omega$ of $\R^N$
$$
J_\A(v) =   \int_\Omega   \frac \A 2 \,|\nn v|^2  + \frac 1 {4\A}  W(v),
$$
whose Euler-Lagrange equation is a scaled version of $\equ{ac}$,
\be
\A^2 \Delta v\, + \, f(v)\, =\, 0 \quad \hbox{in } \Omega \, .
\label{ac1}\ee
We observe that the constant functions $u=\pm 1$ minimize $J_\A$. They are idealized  as two {\em stable phases} of a material in $\Omega$. It is of
interest to analyze stationary configurations in which the two phases coexist.
Given any subset $\Lambda$ of $\Omega$, any discontinuous function of the form
\be v_*= \chi_\Lambda - \chi_{\Omega\setminus \Lambda} \label{u*}\ee
minimizes the second term in $J_\ve$. The introduction of the gradient term in $J_\A$ makes an $\A$-regularization of $u_*$ a test function for which the energy gets bounded and proportional to the surface area of the {\em interface} $M=\partial \Lambda$, so that in addition to minimizing approximately the second term, stationary configurations should also select asymptotically interfaces $M$ that are stationary for surface area, namely (generalized) minimal surfaces.
This intuition on the Allen-Cahn equation  gave important impulse to the calculus of variations, motivating the development of the
theory of {\em $\Gamma$-convergence} in the 1970's. Modica  \cite{modica1} proved that a family of local minimizers $u_\A$ of $J_\A$ with uniformly
bounded energy must converge in suitable sense to a function of the form \equ{u*} where $\pp \Lambda$ minimizes perimeter. Thus, intuitively,
for each given  $\la \in (-1,1)$,
the level sets $[v_\A =\la]$, collapse as $\A \to 0$ onto the interface $\pp\Lambda$.
 Similar result holds for critical points not necessarily minimizers, see \cite{tonegawa}. For minimizers this convergence is known in very strong sense, see \cite{caffarelli1,caffarelli2}.

\medskip
If, on the other hand, we take such a critical point $u_\A$ and scale it around an interior point $0\in \Omega$, setting
$u_\A (x) = v_\A( \A x)$,
then $u_\A$ satisfies equation \equ{ac} in an expanding domain,
$$
\Delta u_\A + f(u_\A ) = 0 \quad \hbox{in } \A^{-1}\Omega$$
so that letting formally $\A \to 0$ we end up with equation \equ{ac} in entire space. The ``interface'' for $u_\A$ should thus be
around the (asymptotically flat) minimal  surface $M_\A = \A^{-1}M$. Modica's result is based on the intuition that if $M$ happens to be a smooth surface, then
the transition from the equilibria $-1$ to $1$ of $u_\A$ along the normal direction  should take place in the approximate
 form $u_\A (x)\approx w(z)$
where $z$ designates the normal coordinate to $M_\A$. Then $w$ should solve  the ODE problem

\be
w'' + f(w) \, =\, 0\quad \hbox{in }\R, \quad w(-\infty)=-1,\ w(+\infty) =1 \ .
\label{edo}\ee
This solution indeed exists thanks to assumption \equ{acn}. It is strictly increasing and unique up to constant translations. We fix in what follows the unique  $w$ for which
\be
\int_\R t\,w'(t)^2 \, dt \ = \ 0\ .
\label{cdm}\ee
For  example \equ{ace}, we have $w(t) = \tanh\left ( t/{\sqrt{2}}\right )$.
In general $w$ approaches its limits at exponential rates,
$$
w(t) -\, \pm 1\ = \ O(\,e^{-\sigma_{\pm} |t|}\, )\quad  \hbox{as } t\to \pm \infty\ .
$$
Observe then that
$$
J_\A (u_\A) \approx  Area\, (M )\, \int_\R [ \frac 12 {w'}^2 + W(w) ]
$$
which is what makes it plausible that $M$ is critical for area, namely a minimal surface.

\medskip
The above considerations led E. De Giorgi \cite{dg} to formulate in 1978 a celebrated
conjecture on  the Allen-Cahn equation    \equ{ac}, parallel to Bernstein's theorem for minimal surfaces: The level sets
$[u=\la]$ of a bounded entire solution $u$ to \equ{ac},  which is also monotone in one direction, must be hyperplanes, at least for dimension $N\le 8$.
Equivalently, up to a translation and a rotation,
$u=w(x_1)$.
 This conjecture has been proven
in dimensions $N=2$ by Ghoussoub and Gui \cite{gg},
 $N=3$ by Ambrosio and Cabr\'e \cite{cabre},
 and under a mild additional assumption by Savin \cite{savin}.
 A counterexample was recently built for $N\ge 9$ in \cite{dkwdg,dkwdg-n},
  see also \cite{cabreterra,jerison}.  See \cite{FV} for
  a recent survey on the state of the art of this question.

\medskip
The assumption of monotonicity in one direction for the solution $u$ in De Giorgi conjecture implies a form of stability, locally minimizing character for $u$ when compactly supported perturbations are considered in the energy. Indeed, if $Z = \pp_{x_N} u >0,$ then the linearized operator
$
L = \Delta  + f'(u),
$
satisfies maximum principle. This implies stability of $u$, in the
sense that its associated quadratic form, namely the second
variation of the corresponding energy, \be \QQ(\psi, \psi) := \int
|\nn \psi|^2 - f'(u) \psi^2 \label{QQ}\ee satisfies $\QQ(\psi,
\psi)>0$ for all $\psi\ne 0$ smooth and compactly supported.
Stability is a basic ingredient in the proof of the conjecture
dimensions $2,3$ in \cite{cabre,gg}, based on finding a control at
infinity of the growth of the Dirichlet integral. In dimension $N=3$
it turns out that \be \int_{B(0,R)} |\nn u|^2 \, =\, O(R^2)
\label{dirichlet}\ee which intuitively means that the embedded level surfaces
$[u=\la]$ must have a finite number of components outside a large
ball, which are all ``asymptotically flat''. The question whether stability alone
suffices for property \equ{dirichlet} remains open. More generally, it is  believed
that this property is equivalent to {\em finite Morse index} of the solution $u$
(which means essentially that $u$
is stable outside a bounded set). The Morse index $m(u)$  is defined as the maximal dimension of
a vector space  $E$ of compactly supported functions such that
$$
\QQ(\psi,\psi)  <0 \foral \psi \in E\setminus \{0\}   .
$$

\medskip
Rather surprisingly,  basically no examples of finite Morse index entire solutions of the Allen-Cahn equation seem known in dimension $N=3$.
Great progress has been achieved in the last decades, both in the theory of semilinear elliptic PDE like \equ{ac}
and in minimal surface theory  in $\R^3$. While this link traces back to the very origins of the study of \equ{ac} as discussed above, it has only been partially explored in producing new solutions.

\medskip
In this paper we construct a new class of entire solutions to the Allen-Cahn equation in $\R^3$ which have the characteristic \equ{dirichlet}, and also finite Morse index, whose level sets resemble  a large dilation of a given complete, embedded minimal surface $M$, asymptotically flat in the sense that it has {\em finite total curvature}, namely
$$
\int_M |K |\, dV\, < \, +\infty
$$
where $K$ denotes Gauss curvature of the manifold, which is also {\em non-degenerate} in a sense that we will make precise below.

 \medskip

As pointed out by Dancer \cite{dancer},  Morse index is a natural element to attempt
classification of solutions of \equ{ac}. Beyond De Giorgi conjecture, classifying solutions with given Morse index should be a natural  step
towards the understanding of the structure
of the bounded solutions of \equ{ac}. Our main results show that, unlike the stable case, the structure of the set of solutions with finite Morse index is  highly  complex. On the other hand, we believe that our construction contains germs of generality, providing elements to extrapolate
what may be true in general, in analogy with classification of embedded minimal surfaces We elaborate on these issues in \S \ref{22}.

\bigskip

\subsection{Embedded minimal surfaces of finite total curvature}
 The theory of embedded,
minimal surfaces of finite total curvature in $\R^3$, has reached a notable development in the last 25 years. For more than a century, only two examples of such surfaces were known: the plane and the catenoid. The first nontrivial example was found in 1981 by C. Costa,  \cite{Costa1, Costa2}.
The {\em Costa surface} is a genus one minimal surface, complete and properly embedded, which outside a large ball has exactly three components (its {\em ends}), two of which are asymptotically catenoids with the same axis and opposite directions, the third one asymptotic to a plane perpendicular to that axis. The complete proof of embeddedness is due to Hoffman and Meeks \cite{hm1}. In \cite{hm2,hm3} these authors generalized notably Costa's example by exhibiting
a class of three-end, embedded minimal surface, with the same look as Costa's far away, but with an array of tunnels that provides arbitrary genus $k\ge 1$. This is known as the Costa-Hoffman-Meeks surface with genus $k$.


\medskip
{Many other examples of multiple-end embedded minimal surfaces have been found since, see for instance \cite{kap1, traizet} and references therein. In general all these surfaces look like parallel planes, slightly perturbed at their ends by asymptotically logarithmic corrections with a certain number of catenoidal links connecting their adjacent sheets.
In reality this intuitive picture is not a coincidence. }
Using the Eneper-Weierstrass representation,  Osserman  \cite{osserman}  established that any embedded, complete minimal surface with finite total curvature  can be described by a conformal diffeomorphism
of a compact surface (actually of a Riemann surface), with a finite number of its points removed. 
These points correspond to the ends. Moreover, after a convenient rotation, the ends are asymptotically all either catenoids or plane, all of them with parallel axes, see Schoen \cite{schoen 1}. The topology of the surface is thus characterized by the genus of the compact surface and the number of ends, having therefore ``finite topology''.


\medskip
\subsection{Main results}
In what follows  $M$ designates a complete, embedded minimal surface in  $\R^3$ with finite total curvature (to which below we will make a further nondegeneracy assumption).
 As pointed out in  \cite{hk1}, $M$ is orientable and the set
 $\R^3\setminus M$ has exactly two components
$S_+$, $S_-$.
In what follows we fix a continuous choice of unit normal field  $\nu(y)$, which conventionally we  take it to point towards $S_+$.

\medskip
For  $x= (x_1,x_2,x_3)= (x',x_3)\in \R^3$,  we denote $$r = r(x) = |(x_1,x_2)|= \sqrt{x_1^2 + x_2^2} . $$
After a suitable rotation of the coordinate axes,
outside the infinite cylinder  $ r < R_0$ with sufficiently large radius $R_0$, then
$M$ decomposes into a finite number $m$ of unbounded components $M_1,\ldots, M_m$, its {\em ends}. From a result in  \cite{schoen 1}, we know that asymptotically each end of $M_k$  either resembles a plane or a catenoid. 
More precisely,
 $M_k$ can be represented as the graph of a function
$F_k$ of the first two variables,
$$
M_k = \{\, y\in \R^3 \ / \ r(y)> R_0,\   y_3 =  F_k(y')\, \}
$$
where $F_k$ is a smooth function which  can be expanded as
\be
F_k(y') = a_k \log r +  b_k  + b_{ik} \frac {y_i} {r^2} +  O( r^{-3})\quad\hbox{as } r\to + \infty,
\label{ends}\ee
for certain constants  $a_k$, $b_k$, $b_{ik}$,
and this relation can also be differentiated.
Here
\begin{equation}
\label{balfor}
 a_1\leq a_2 \leq \ldots \leq a_m\ ,\qquad  \sum_{k=1}^m a_k \ =\ 0 \, .
\end{equation}
The direction of the normal vector $\nu(y)$ for large $r(y)$ approaches on the ends that of the $x_3$ axis, with alternate signs.
We use the convention that for $r(y)$ large  we have
\begin{equation}
\label{asyends}
\nu (y)\ = \ \frac {(-1)^k } { \sqrt{1 +|\nn F_k(y')|^2}}\, (\nn F_k(y')\, ,\, -1\,)\quad\hbox{if } y\in M_k.
\end{equation}

\medskip
Let us consider the Jacobi operator of $M$
\be
\JJ(h)\, :=\, \Delta_M h + |A|^2 h
\label{JJ}\ee
where $|A|^2= -2K$ is the Euclidean norm of the second fundamental form of $M$. $\JJ$ is the linearization of the mean curvature
operator with respect  to perturbations of $M$ measured along its normal direction.
A smooth function $z(y)$ defined on $M$  is called a {\em Jacobi field} if $\JJ(z) =0$.
 Rigid motions of the surface induce naturally some bounded Jacobi fields:
Associated to respectively translations along coordinates axes and rotation around the $x_3$-axis, are the functions
$$z_{1}(y) = \nu(y)\cdot e_i, \quad y\in M,\quad  i=1,2,3, $$
 \be z_{4}(y) = (-y_2, y_1,0) \cdot \nu(y), \quad y\in M. \label{jfields}\ee


We assume that $M$ is {\em non-degenerate} in the sense that these functions are actually {\em all} the bounded Jacobi fields,
namely
\be
  \{ \ z\in L^\infty (M)\ / \ \JJ(z) = 0\, \} \ = \ {\rm span}\, \{\, z_{1}, z_{2}, z_{3}, z_{4}\, \}\, .
\label{kill}\ee
We denote in what follows by $J$ the dimension $(\le 4)$ of the above vector space.

\medskip
This assumption, expected to be generic for this class of surfaces, is known in some important cases, most notably the catenoid and the Costa-Hoffmann-Meeks
surface which is an example of a three ended $M$ whose genus may be of any order. See  Nayatani \cite{Nay1,Nay2} and Morabito \cite{Mo}. Note that for a catenoid, $z_{04}=0$ so that $J=3$.  Non-degeneracy has been used as a tool to build new minimal surfaces for instance
in  Hauswirth and Pacard \cite{Hau-Pacard}, and
 in P\'erez and Ros \cite{perez-ros}. It is also the basic element, in a compact-manifold version, to build solutions to the small-parameter Allen-Cahn
equation in Pacard and Ritor\'e \cite{pacard}.

\medskip

In this paper we will construct a solution to the Allen Cahn equation whose zero level sets look like a large dilation
of the surface $M$, with ends perturbed logarithmically. 
Let us consider a large dilation of $M$,
$$
M_\A \ := \ \A^{-1} M .
$$
This dilated minimal surface has ends parameterized as
$$
M_{k, \A}\, =\,  \{\, y\in \R^3 \ / \ r(\A y)> R_0,\   y_3 =  \A^{-1}F_k(\A y')\, \}\ .
$$
Let $\beta$ be a vector of given $m$ real numbers with
\be
\beta \ =\ (\beta_1, \ldots, \beta_m ),\qquad \sum_{i=1}^m \beta_i = 0\ .
\label{beta1}\ee
Our first result asserts the existence of a solution $u= u_{\A}$ defined for all sufficiently small $\A>0$ such that given  $\la \in (-1,1)$,
its level set  $[ u_{\A}=\la]$  defines an embedded surface lying at a uniformly bounded  distance in $\A$ from the surface $M_\A$, for points
with  $r(\A y) = O(1)$, while  its $k$-th end, $\ k=1,\ldots, m,$ lies at a uniformly bounded distance from the graph
\be
\ r(\A y)> R_0,\   y_3 \, = \, \, \A^{-1}\, F_k(\A y')  + \beta_k \log |\A y'| \,  .
\label{mkbeta}\ee


\medskip
The parameters $\beta$  must satisfy an additional constraint.
It is  clear that if two ends are parallel, say $a_{k+1} = a_k$,  we need at least that $\beta_{k+1} -\beta_k \ge 0$, for otherwise the ends
would eventually  intersect. Our further condition on these numbers is that these ends in fact diverge at a sufficiently fast rate.  We require
\be
\beta_{k+1}-\beta_k\,> \, 4\, \max\, \{\sigma^{-1}_-, \sigma^{-1}_+\} \quad\hbox{if}\quad  a_{k+1} = a_k \ .
\label{condbeta1}
\ee
Let us consider  the smooth map
\be
 X(y,z) = y +  z \nu(\A y)  ,\quad (y, t)\in M_\A \times \R.
\label{x-1}\ee
$x= X(y,z)$ defines coordinates  inside the image of any region where the map is one-to-one. In particular, let us consider a function $p(y)$ with $$p(y) = (-1)^k \beta_k \log |\A y'| + O(1), \quad k=1,\ldots, m, $$
and $\beta$ satisfying $\beta_{k+1} -\beta_k >\gamma>0$ for all $k$ with $a_k= a_{k+1}$. Then the map $X$ is one-to-one for all small $\A$
in  the region of points $(y,z)$ with $$|z- q(y)| < \frac \delta\A + \gamma \log (1+ |\A y'|) $$
provided that $\delta>0$ is chosen sufficiently small.

\begin{teo}\label{teo1} Let $N=3$ and $M$ be a minimal surface embedded, complete with finite total curvature which is nondegenerate. Then,
given $\beta$ satisfying relations $\equ{beta1}$ and $\equ{condbeta1}$, there exists a bounded solution $u_{\A} $ of equation $\equ{ac}$,
defined for all sufficiently small $\A$, such that
 \be
 u_{\A}(x) = w( z-  q(y) ) + O(\A)\, \foral  x= y +  z \nu(\A y), \quad |z-  q (y)| < \frac\delta\A,
\label {eq} \ee
where the function $q$  satisfies
$$
q(y)\, = \,  (-1 )^k \beta_k\, \log |\A y'| + O(1)  \quad  y\in M_{k,\A} ,\quad k=1,\ldots, m.
$$
In particular, for each given  $\la \in (-1,1)$, the level set $[u_\A=\la]$ is an embedded surface that
 decomposes for all sufficiently small $\A$ into $m$ disjoint components (ends) outside a bounded set. 
The $k$-th end lies at $O(1)$ distance from the graph
$$  y_3 \, = \, \, \A^{-1}\,F_k(\A y)   + \beta_k \log |\A y'|. $$

\end{teo}

\bigskip The solution predicted by this theorem depends, for fixed $\A$,
on $m$ parameters. Taking into account the constraint $\sum_{j=1}^m \beta_j =0$ this gives $m-1$ independent parameters corresponding to
logarithmic twisting of the ends of the level sets.  Let us observe that consistently, the combination $\beta \in \hbox{Span}\,\{(a_1,\ldots, a_m)\}$
can be set in correspondence with moving $\A$ itself, namely with a dilation parameter of the surface. We are thus left with $m-2$ parameters
for the solution in addition to $\A$.
Thus, besides the trivial rigid motions of the solution,  translation along the coordinates axes, and rotation about the $x_3$ axis, this family of
solutions depends exactly on $m-1$ ``independent'' parameters. Part of the conclusion of our second result is that the bounded kernel of the linearization of equation \equ{ac} about one of these solutions is made up exactly of the generators of the rigid motions, so that in some sense
the solutions found are $L^\infty$-isolated, and the set of bounded solutions nearby is actually $m-1+J$-dimensional. A result parallel to this one,     in which the {\em moduli space} of the minimal surface $M$ is described by a similar number of parameters, is found in \cite{perez-ros}.

\bigskip
 Next we discuss the connection of the Morse index of the solutions of Theorem \ref{teo1} and the {\em index} of the minimal surface $M$, $i(M)$, which has a similar definition
relative to the quadratic form for the Jacobi operator:
The number $i(M)$ is  the largest dimension for a vector spaced $E$ of compactly supported smooth functions in $M$ with
$$
\int_M |\nn k|^2\, dV \, -\, \int_M |A|^2 k^2\, dV \  < 0 \foral k \in E\setminus \{0\} .
$$
We point out that for complete, embedded surfaces, finite index is equivalent to finite total curvature, see \cite{gulliver} and also \S 7 of \cite{hk1} and  references therein.   Thus, for our surface $M$,
$i(M)$ is indeed finite. Moreover, in the Costa-Hoffmann-Meeks surface it is known that  $i(M) = 2l-1$ where $l$ is the genus of $M$.  See \cite{Nay1}, \cite{Nay2} and \cite{Mo}.

\medskip
Our second result is that the Morse index and non-degeneracy of $M$ are transmitted  into the linearization of equation \equ{ac}.

\begin{teo}\label{teo2}
 Let $u_{\A} $ the solution of problem $\equ{ac}$ given by  Theorem \ref{teo1}. Then for all sufficiently small $\A$ we have
 $$ m(u_{\A})\  =\  i(M). $$
 Besides, the solution is non-degenerate, in the sense that any bounded solution of
 $$
 \Delta \phi + f'(u_{\A}) \phi =0 \quad\hbox{ in }\R^3
 $$
 must be a linear combination of the functions $Z_{i},\ i=1,2,3,4$ defined as
 $$
 Z_{i} = \partial_i u_{\A}, \quad i=1,2,3,\quad Z_{4} = -x_2\partial_1 u_{\A}+ x_1 \pp_2 u_{\A}.
 $$

\end{teo}

\bigskip
We will devote the rest of this paper to the proofs of Theorems 1 and 2.

\setcounter{equation}{0}
\section{ \emph{ The Laplacian  near  $M_\A$}}

\medskip
\subsection{ The Laplace-Betrami Operator of $M_\A$}

Let $D$ be the set
$$
D\ =\ \{ \py\in \R^2 \ /\  |\py|> R_0\}.  $$
We can parameterize the end $M_k$ of $M$ as
\be
\py \in D \longmapsto  y:= Y_k(\py ) \ = \  \py_i e_i + F_k(\py) e_3
 \, .
\label{Yk}\ee
and $F_k$ is the function in \equ{ends}. In other words, for $y =(y',y_3)\in M_k$ the coodinate $\py$ is just defined as $\py = y'$.
We want to represent $\Delta_M$--the Laplace-Beltrami operator of $M$--with respect to these coordinates.
For the coefficients of the metric $g_{ij}$ on $M_k$ we have
$$
\pp_{\py_i}Y_k =   e_i +  O\left (  r^{-1} \right ) e_3
$$
so that
\be
g_{ij}(\py)  = \left < \pp_{i}Y_k ,\pp_{j}Y_k \right > = \delta_{ij} +   O\left (  r^{-2}\right ) ,
\label{gij}
\ee
where $r= |\py|$.
The above relations ``can be differentiated'' in the sense that differentiation makes the terms  $O( r^{-j} )$ gain corresponding negative powers of $r$.
Then we find the representation

\be
\Delta_M = \frac 1{ \sqrt{\det g_{ij} }}
 \partial_i(  \sqrt{\det g_{ij} }\,  g^{ij} \partial_j) = \Delta_{\py} +   O( r^{-2})\partial_{ij}  +  O( r^{-3})\ \partial_i   \quad\hbox{on } M_k \, .
\label{DeltaM}\ee

\medskip
The normal vector to $M$ at $y\in M_k$  $k=1,\ldots, m$, corresponds to
$$
\nu (y ) \ =\ (-1)^k\frac 1{\sqrt{ 1+ |\nn F_k(\py )|^2}} \, \left ( \, \partial_i F_k(\py) e_i - e_3 \, \right)
\ , \quad y = Y_k(\py)\in M_k\
$$
so that
\be
\nu (y) \ =\ (-1)^k e_3 +   \A_k  r^{-2} \, \py_i e_i   + O(r^{-2}) \ , \quad y = Y_k(\py)\in M_k\ .
\label{normal-1}\ee
Let us observe for later reference that since
$
\partial_i \nu =   O( r^{-2}),
$
then the principal curvatures of $M$, $k_1,k_2$ satisfy
$
k_l = O( r^{-2}) $.
In particular, we have that
\be
|A(y)|^2 = k_1^2 + k_2^2 = O( r^{-4}) .
\label{A2}\ee
To describe the entire manifold  $M$ we consider a finite number $N\ge m+1$ of local parametrizations


\be
 \py \in {\mathcal U}_k \subset \R^2\, \longmapsto\, y = Y_k(\py ),\quad  Y_k\in C^\infty (\bar{\mathcal U}_k )  , \quad k=1,\ldots, N.
\label{chart}\ee
For $k=1,\ldots, m$  we choose them to be those in \equ{Yk}, with ${\mathcal U}_k = D$,  so that $Y_k({\mathcal U}_k) = M_k$, and
 $\bar {\mathcal U}_k$ is bounded for $k= m+1,\ldots, N$. We require then that
  $$
  M= \bigcup_{k=1}^N Y_k({\mathcal U}_k ).
  $$
 We remark that the Weierstrass representation of $M$
 implies that we can actually take $N=m+1$, namely only one extra parametrization is needed to describe the bounded complement of the ends in $M$.
We will not use this fact.  In general, we represent for $y\in Y_k(  {\mathcal U}_k)$,
\be
\Delta_M  =  a_{ij}^0(y)\partial_{ij} + b_i^0(y) \partial_i,\quad  y =Y_k(\py),\quad \py \in {\mathcal U}_k,
\label{lapM}\ee
where $a_{ij}^0$ is a uniformly elliptic matrix and the index $k$ is not made explicit in the coefficients. For $k=1,\ldots, m$ we have
\be
a_{ij}^0(y) = \delta_{ij}  + O(r^{-2}),\quad  b_i^0 = O(r^{-3}), \quad\hbox{as } r(y) = | \py|\to \infty.
\label{lapM1}\ee

\medskip
The parametrizations set up above induce naturally a description of the expanded manifold
$M_\A = \A^{-1} M $ as follows. Let us consider
the functions
\be
Y_{k\A}\, :\, {\mathcal U}_{k\A} := \A^{-1} {\mathcal U}_k  \to M_\A,\quad \py\mapsto  Y_{k\A} (\py ):=  \A^{-1} Y_k(\A \py ),\quad  k=1,\ldots, N.
\label{chart2}\ee
Obviously we have
$$
  M_\A = \bigcup_{k=1}^N Y_{k\A}  ({\mathcal U}_{k\A} ).
  $$
The computations above lead to the following representation for the operator $\Delta_{M_\A}$:
\be
\Delta_{M_\A}   =  a_{ij}^0(\A y)\partial_{ij} + b_i^0(\A y) \partial_i,\quad  y =Y_{k\A} (\py),\quad \py \in {\mathcal U}_{k\A},
\label{lapMA}\ee
where $a_{ij}^0$, $b_i^0$ are the functions in \equ{lapM}, so that for $k=1,\ldots, m$ we have
\be
a_{ij}^0 = \delta_{ij}  + O(r_\A ^{-2}),\quad  b_i^0 = O(r_\A ^{-3}), \quad\hbox{as } r_\A (y) := |\A \py|\to \infty.
\label{lapMA1}\ee

\bigskip
\subsection{The Euclidean Laplacian near $M_\A$}\label{lma}

We will describe in coordinates relative to $M_\A$ the Euclidean Laplacian $\Delta_x$, $x\in\R^3$, in a setting needed for the proof of our main results.
Let us consider a smooth function $h: M \to \R$, and the smooth map $X_h$ defined as
\be
X_h \, :\, M_\A \times \R \to \R^3,\quad (y,t)\longmapsto X_h(y,t):= y + (t+ h(\A y)\,)\,\nu(\A y)
\label{Xh}\ee
where $\nu$ is the unit normal vector to $M$. Let us consider an open subset $\OO$ of  $M_\A \times \R  $ and assume
that the map $X_h |_\OO $ is one to one, and that it defines a diffeomorphism onto its image ${\mathcal N}= X_h(\OO)$. Then
$$ x= X_h(y,t), \quad (y,t)\in \OO, $$
defines smooth coordinates to describe the open set ${\mathcal N}$ in $\R^3$. Moreover,
the maps
$$
x\ =\  X_h(\, Y_{k\A} (\py)\,,t) , \quad (\py,t) \in   ({\mathcal U}_{k\A} \times \R) \, \cap \OO,\quad k=1,\ldots, N,
$$
define local  coordinates $(\py,t)$  to describe the region
$ \mathcal{N}$.
We shall assume in addition that for certain small number $\delta>0$, we have
\be
\OO \subset \{ (y,t)\ /\ |t+ h(\A y)| < \frac \delta\A \,\log (2+ r_\A(y)\, )\, \}.
\label{condO}\ee

 We have the validity of the following expression for the Euclidean Laplacian operator in $\mathcal N$.

\begin{lemma}\label{lapfinal} For $x= X_h(y,t)$, $(y,t)\in \OO$ with  $y= Y_{k\A}(\py)$, $\py \in {\mathcal U}_{k\A}$,  we have
the validity of the identity
$$
\Delta_x   =  \, \partial_{tt}  + \Delta_{M_\A}  - \A^2[ (t+h)|A|^2 + \Delta_M h]  \partial_t\ -
2 \A \, a_{ij}^0\, \partial_jh \partial_{it}  \ +
$$
$$
 \A (t+h) \, [ a_{ij}^1\partial_{ij}   - 2\A\, a_{ij}^1\, \partial_ih \partial_{jt}
  \, +\, \A \,b_i^1\,(\partial_i   - \A  \partial_ih \partial_t)\,  \, ] \ +
 $$

\be
  \A^3 (t+h)^2b^1_3  \partial_t\, +\, \A^2 [\,a_{ij}^0 + \A(t+h) a_{ij}^1 \,]\partial_ih\partial_jh \, \partial_{tt}\ .
\label{laplacianfinal}\ee
Here, in agreement with  $\equ{lapMA}$,
$
\Delta_{M_\A}  =  a_{ij}^0(\A y)\partial_{ij} + b_i^0(\A y) \partial_i.
$\\
The functions $a_{ij}^1$, $b_i^1$, $b_3^1$ in the above expressions appear evaluated at the pair \\ $(\A y , \A(t+h(\A y) )$,
while the functions $h$, $\pp_ih$, $\Delta_M h$, $|A|^2$, $a_{ij}^0$, $b_i^0$ are evaluated at $\A y$
In addition, for $k=1,\ldots, m$, $l=0,1$,
$$
a_{ij}^l = \delta_{ij}\delta_{0l}  + O(r_\A^{-2}),\quad  b_i^l = O(r_\A^{-3}), \quad  b_3^1 = O(r_\A^{-6}) \ ,
$$
as $ r_\A (y) = |\A \py|\to \infty$, uniformly in their second variables.
The notation $\partial_j h$ refers to $\pp_j [ h \circ Y_k] $.

\end{lemma}
We postpone proof of this fact for the appendix. The proof actually yields that the coefficients $a_{ij}^1$ and $b_i^1$ can be further expanded as follows:

$$
a_{ij}^1  = a_{ij}^1 (\A y,0)  + \A (t+h)\,a_{ij}^{(2)} (\A y, \A (t+h)) =: a_{ij}^{1,0} + \A (t+h)a_{ij}^2,
$$
with $a_{ij}^{(2)} = O( r_\A^{-3})$, and similarly
$$
 b_{j}^1 = b_j^1 (\A y,0)  + \A (t+h)\,b_{j}^{(2)} (\A y, \A (t+h))=: b_j^{1,0} + \A (t+h)b_{j}^2,
 $$
 with $b_{j}^{(2)} = O( r_\A^{-4})$.
As an example of the previous formula, let us compute the Laplacian of a function that separates variables
$t$ and $y$, that will be useful in \S \ref{aproz} and  \S 11.

\begin{lemma}\label{formula}
Let
$
v(x)\, = \, k(y)\, \psi(t)\, .
$
Then the following holds.

$$
\Delta_x v \  =  \,  k \psi''  \, +
 \psi \, \Delta_{M_\A}k  \,-\,  \A^2[ (t+h)|A|^2 + \Delta_M h]\,k\, \psi'  \, - \, 2 \A\, a_{ij}^0\, \partial_jh \partial_{i}k  \, \psi'  \ +
$$

$$
 \A (t+h) \, [ a_{ij}^{1,0}\partial_{ij}k \psi   -  2\A\, a_{ij}^{1,0}\, \partial_jh \partial_{i}k  \, \psi' \,  +
 \A (b_i^{1,0}\partial_ik\,\psi    - \A b_i^{1,0} \partial_ih\,k\, \psi')\, \, ] \ +
 $$

 $$
 \A^2 (t+h)^2 \, [ a_{ij}^2\partial_{ij}k \psi   - 2\A\, a_{ij}^2\, \partial_jh \partial_{i}k \, \psi' \,  +
 \A (b_i^2\partial_ik\,\psi    - \A b_i^2 \partial_ih\,k\, \psi')\, \, ] \ +
$$

\be
  \A^3 (t+h)^2b^1_3 \, k\, \psi' \, +\, \A^2 [\,a_{ij}^0 + \A(t+h) a_{ij}^1 \,]\partial_ih\partial_jh \,k\, \psi'' \ .
\label{r88}\ee

\end{lemma}


\emph{}

\setcounter{equation}{0}
\section{ \emph{ Approximation of the solution and preliminary discussion}}\label{aproz}

\subsection{Approximation of order zero and its projection}
Let us consider a function $h$ and sets $\OO$ and $\NN$ as in \S \ref{lma}. Let $x= X_h(y,t)$ be the coordinates introduced in \equ{Xh}.
At this point we shall make a more precise  assumption about the function $h.$ We need the following preliminary result whose proof we postpone
for \S \ref{logjac}.

\medskip
We consider a fixed  $m$-tuple of real numbers ${\bf \beta} = (\beta_1,\ldots , \beta_m)$ such that
\be
\sum_{i=1}^m \beta_j = 0.
\label{beta}\ee

\begin{lemma}\label{lemin6}
Given any real numbers $\beta_1,\ldots, \beta_m$ satisfying $\equ{beta}$,
there exists a smooth function $h_0(y)$ defined on $M$ such that
$$
\JJ (h_0)= \Delta_M h_0 + |A|^2 h_0  =0 \quad \hbox{in } M,
$$
$$
h_0(y) = (-1)^j \beta_j \log r +   \theta  \hbox{ as } r\to \infty \quad \hbox{in } M_j \ \foral y\in M_j \, ,
$$
where $\theta$ satisfies
\be
\|\theta\|_\infty + \|r^2D\theta\|_\infty < +\infty \ .
\label{nb}\ee
\end{lemma}

\bigskip

We fix  a function $h_0$ as in the above lemma and consider a function  $h$ in the form
$$ h=  h_0 +h_1 .$$
We allow  $h_1$ to be a parameter which we will adjust.
For now we will assume that for a certain constant $\KK$ we have
\be
\| h_1\|_{L^\infty (M)} + \|(1+ r^2) Dh_1\|_{L^\infty (M)}  \ \le \ \KK \A \, .
\label{ass}\ee
We want to find a solution to $$S(u):= \Delta_x u + f(u)=0.$$
We consider in the region ${\mathcal N}$  the approximation
$$u_0(x):= w(t)  = w(z- h_0(\A y)- h_1(\A y) )  $$
where $z$ designates the normal coordinate to $M_\A$. Thus, whenever $\beta_j\ne 0$, the  level sets $[u_0=\la]$ for a fixed $\la\in (-1,1)$ 
departs logarithmically from
the end $\A^{-1}M_j$ being still asymptotically catenoidal, more precisely
it is described as the graph
$$
y_3 =   ( \A^{-1} a_j   +  \beta_j)\, \log r \, +\,  O(1) \hbox{ as } r\to \infty   .
$$
Note that, just as in the minimal surface case, the coefficients of the ends are balanced in the sense that they add up to zero.

\medskip
It is  clear that if two ends are parallel, say $a_{j+1} = a_j$,  we need at least that $\beta_{j+1} -\beta_j \ge 0$, for otherwise the ends
of this zero level set would eventually  intersect.
We recall that our further condition on these numbers is that these ends in fact diverge at a sufficiently fast rate:
\be
\beta_{j+1}-\beta_j \,> \, 4\, \max\, \{\sigma^{-1}_-, \sigma^{-1}_+\} \quad\hbox{if}\quad  a_{j+1} = a_j \ .
\label{condbeta}
\ee
We will explain later the role of this condition.
Let us evaluate the error of approximation $S(u_0)$. Using Lemma \ref{formula}  and the fact that
$w''+ f(w)=0$, we find

$$
S(u_0) := \Delta_x u_0 +f(u_0) \  =  \,
$$

$$
 - \A^2[ |A|^2h_1  + \Delta_M h_1 ]\,w'  \ +
$$

$$
  - \A^2 |A|^2\, tw' \ + 2 \
  \A^2  a_{ij}^0 \,\partial_ih_0 \partial_jh_0\, w''\  +
$$

$$
\A^2 \, \,a_{ij}^0
\, ( 2\partial_ih_0 \partial_j h_1  + \partial_ih_1 \partial_j h_1\,) \, w'' \, +\,
$$

$$
 2 \A^3(t+ h_0 +h_1) a_{ij}^1 \, \partial_i(h_0 +h_1) \partial_j(h_0+ h_1) \, w''\, +\,
 $$

 \be
\A^3 (t+ h_0 +h_1)  b_i^1 \partial_i(h_0 +h_1)\,w' +  \A^3 (t+h_0+ h_1)^3b^1_3 w'
\label{error1-1}
\ee
where the formula above has been broken into ``sizes'', keeping in mind that $h_0$ is fixed while $h_1=O(\A)$.
Since we want that $u_0$ be as close as possible to be a solution of \equ{ac}, then  we would like to choose $h_1$ in such a way that the quantity \equ{error1-1} be as small as possible.
Examining the above expression, it does not look like we can do that in absolute terms. However part of the error could be made smaller by adjusting $h_1$. Let us consider the
``$L^2$-projection" onto $w'(t)$ of the error for each fixed $y$, given by

$$
\Pi(y) := \int_{- \infty}^{\infty}   S(u_0)(y,t)\, w'(t)\, dt\
$$
where for now, and for simplicity we assume the coordinates are defined for all $t$, the difference with the integration is taken in all the actual
domain for $t$ produces only exponentially small terms in $\A^{-1}$.
Then we find
$$
\Pi(y)=  \A^2 (\Delta_M h_1 + h_1|A|^2)\int_{- \infty}^{\infty}  {w'}^2 dt  + \A^3\partial_{i} (h_0+h_1)\int_{- \infty}^{\infty} b_i^1(t+h_0+h_1) {w'}^2 dt\
+
$$
\be
\A^3  \partial_i (h_0 + h_1)\partial_j(h_0 + h_1) \int_{- \infty}^{\infty} (t+h_0+h)a_{ij}^1 w'' w' dt  +  \A^3 \int_{- \infty}^{\infty}
(t+h_0+ h_1)^3b^1_3 {w'}^2 dt\
\label{projection}\ee
 where we have used $\int_{-\infty }^\infty  t{w'}^2\, dt = \int_{-\infty}^\infty w''w'\, dt = 0$ to get rid in particular of the terms
  of order $\A^2$.

\medskip
Making all these ``projections'' equal to zero amounts to a nonlinear differential equation for $h$ of the form
\be
{\mathcal J}(h_1) = \Delta_M h_1 + h_1|A( y)|^2  = G_0(h_1) \quad  y \in M
\label{probjac}\ee
where $G_0$ is easily checked to be a contraction mapping of small constant in $h_1$,
  in the ball radius $O(\A)$ with the $C^1$ norm defined by the expression in
the left hand side of inequality \equ{ass}. This is where the nondegeneracy assumption on the Jacobi operator ${\mathcal J}$ enters, since we
would like to invert it,  in  such a way to set up equation \equ{probjac} as a fixed point problem for a contraction mapping of a ball of the form
\equ{ass}.

\subsection{Improvement of approximation} \label{u1}
The previous considerations are  not sufficient since even after adjusting optimally $h$, the error in absolute value does not necessarily decrease.
As we observed, the  ``large'' term in the error, $$-\A^2|A|^2tw' + \A^2a_{ij}^0 \pp_i h_0\pp_j h_0\, w''$$  did not contribute to the projection.
In order to eliminate, or reduce the size of this remaining part $O(\A^2)$ of the error, we  improve the approximation through the following argument.
Let us
consider the differential equation
$$
\psi_0''(t) + f'(w(t))\psi_0 (t) =  tw'(t),
$$
which has a unique bounded solution with $\psi_0(0) =0$, 
given explicitly by the formula 
$$
\psi_0(t) = w'(t)  \int_0^t w'(t)^{-2} \int_{-\infty}^s sw'(s)^2ds \, .
$$
Observe that this function is well defined and it is bounded since $\int_{-\infty}^\infty sw'(s)^2ds  = 0$ and $w'(t)\sim e^{ - \sigma_\pm |t|}$ as $t\to \pm \infty$, with
$\sigma_\pm >0$.  Note also that $\psi_1(t) = \frac 12 tw'(t)$ solves
$$
\psi_1''(t) + f'(w(t))\psi_1 (t) =  w''(t) \, .
$$
We consider as a second approximation
\be
u_1 = u_0  + \phi_1 ,\quad \phi_1(y,t) :=  \A^2 |A(\A y)|^2 \psi_0 (t)  - \A^2 a_{ij}^0 \pp_i h_0\pp_j h_0 (\A y)\, \psi_1 (t) \ .
\label{u1-1}\ee

Let us observe that
$$
S(u_0 + \phi) = S(u_0) + \Delta_x \phi + f'(u_0)\phi +  N_0(\phi), \quad N_0(\phi) = f(u_0+\phi) -f(u_0) -f'(u_0)\phi \, .
$$
We have that $$\partial_{tt}\phi_1 + f'(u_0)\phi_1 =  \A^2 |A(\A y)|^2 tw'  - \A^2 a_{ij}^0 \pp_i h_0\pp_j h_0 (\A y)\, w'' \, .$$
Hence  we get that the largest remaining term in the error is canceled.
Indeed, we have
$$
S(u_1) =   S(u_0)- (2 \A^2 a_{ij}^0 \pp_i h_0\pp_j h_0\, w'' -\A^2 |A(\A y)|^2 tw') +  [ \Delta_x -\partial_{tt}] \phi_1 + N_0(\phi_1) .
$$
Since $\phi_1$ has size of order $\A^2$, a smooth dependence in $\A y$ and it is of size $O( r_\A^{-2} e^{-\sigma|t|})$
using  Lemma \ref{formula}, we readily check  that the ``error created''
$$
[ \Delta_x -\partial_{tt}] \phi_1 + N_0(\phi_1) : =   - \A^4 \, (\,|A|^2 t\psi_0' - a_{ij}^0 \pp_i h_0\pp_j h_0\,t \psi_1'\, )\, \Delta h_1\ +  R_0
$$
satisfies
$$
|R_0(y,t)| \le C \A^3(1+ r_\A(y))^{-4} e^{-\sigma|t|} .
$$
Hence  we have eliminated the $h_1$-independent term $O(\A^2)$ that did not contribute to the projection $\Pi(y)$, and replaced it by one smaller and with faster
decay. Let us be slightly more explicit for later reference. We have

$$
S(u_1) := \Delta u_1 +f(u_1) \  =  \,
$$

$$
 - \A^2[ |A|^2h_1  + \Delta_M h_1 ]\,w'  +
\A^2  \,a_{ij}^0
\, ( \partial_ih_0 \partial_j h_1 + \partial_ih_1 \partial_j h_0 + \partial_ih_1 \partial_j h_1\,) \, w''
$$

\be
 - \A^4 \, (\,|A|^2 t\psi_0' - a_{ij}^0 \pp_i h_0\pp_j h_0\,t \psi_1'\, )\, \Delta_M h_1\, + 2
\A^3(t+ h) a_{ij}^1 \, \partial_i h \partial_jh \, w'' +
R_1
\label{error1}
\ee
where
$$R_1= R_1(y,t,h_1 (\A y) , \nn_M h_1 (\A y) )$$ with
$$
 | D_{\imath} R_1(y, t, \imath, \jmath)| + | D_{\jmath} R_1(y, t, \imath, \jmath)| +| R_1(y,t, \imath, \jmath)| \le C \A^3(1+ r_\A(y))^{-4} e^{-\sigma|t|}
$$
and the constant $C$ above possibly depends on the number $\KK$ of condition \equ{ass}.


\medskip
The above arguments are in reality the way we will actually solve the problem: two separate, but coupled steps are involved:
(1) Eliminate the parts of the error that do not contribute to the projection $\Pi$ and (2) Adjust $h_1$  so that
the projection $\Pi$ becomes identically zero.

\medskip
\subsection{The condition of diverging ends}\label{diverg}
Let us explain the reason to introduce condition \equ{condbeta} in the parameters $\beta_j$. To fix ideas, let us assume that we have two consecutive planar ends of $M$, $M_j$ and $M_{j+1}$,
namely with $a_j=a_{j+1}$ and with $d= b_{j+1} - b_j >0$. Assuming that the normal in $M_j$ points upwards, the coordinate
 $t$ reads approximately as
$$ t=  x_3 - \A^{-1}b_j  -h \quad\hbox{ near } M_{j\A} ,\quad  t=  \A^{-1}b_{j+1}- x_3-h  \quad\hbox{ near } M_{j+1\A} .
 $$
 If we let $h_0\equiv 0$ both on $M_{j\A}$ and $M_{j+1\A}$ which are separated at distance $d/\A$, then a good  approximation in the entire region
 between  $M_{j\A}$ and $M_{j+1\A}$ that matches the parts of $w(t)$ coming both from $M_j$ and $M_{j+1}$
 should read  near  $M_j$ approximately as
$$ w(t) +  w(\A^{-1}d -t ) - 1 .$$
When computing the error of approximation, we observe that the following additional term arises near $M_{j\A}$:
$$
E:= f( \,w(t) +  w(\A^{-1}d -t ) - 1\,) \,- \,f( w(t)) - f(\,w(\A^{-1}d -t )\,) \ \sim
$$
$$
\sim \ [f'(w(t))-f'(1)\, ]\,  (\,w(\A^{-1}d -t ) - 1\,) \, .
$$
Now in the computation of the projection of the error this would give rise to
$$
\int_{-\infty}^\infty  [\, f'(w(t))-f'(1)\, ]\, (\,w(\A^{-1}d -t )\, -\, 1\,) \, w'(t)\,dt\  \sim
c_* e^{-\sigma_+ \, \frac d\A  } .
$$
where $c_*\ne 0$ is a constant. Thus equation \equ{probjac} for $h_1$ gets modified with a term which even though very tiny, it
has no decay as $|y|\to +\infty$ on $M_j$, unlike the others involved in the operator $G_0$ in \equ{probjac}. That terms eventually dominates and the equation for $h_1$
for very large $r$ would read in $M_j$ as $$\Delta_M h_1 \sim
 e^{- \frac \sigma \A  }\ne 0, $$
which is inconsistent with the assumption that $h$ is bounded.  Worse yet, its solution would be quadratic thus eventually
intersecting  another end.  This nuisance is fixed with the introduction of $h_0$ satisfying condition \equ{condbeta}. In that case the term $E$
created above
will now read near $M_{j\A}$ as
$$
E \sim  C e^{-\sigma_+ \frac d\A } \, e^{- (\beta_{j+1}- \beta_j)\log r_\A }\, e^{-\sigma|t|}\,  = \,   O( e^{-\frac \sigma \A} r_\A^{- 4} e^{-\sigma|t|})$$ which is qualitatively of the same type of the other terms  involved in the computation of the error.

\subsection{The global first approximation}
The approximation $u_1(x)$ in \equ{u1} will be sufficient for our purposes, however it is so far defined only in a region
of the type ${\mathcal N}$ which we have not made precise yet.
Since we are assuming that $M_\A$ is connected, the fact that $M_\A$ is properly embedded implies
that $\R^3 \setminus  M_\A$ consists of precisely two components $S_-$ and $S_+$. Let us use the convention that $\nu$ points in the direction of $S_+$. Let us consider the function
 $\HH$ defined in $ \R^3 \setminus  M_\A$ as
\be
\HH(x)\, := \,
\left\{
\begin{matrix}
\ 1 &\quad\hbox{if } x\in S_+ \\
-1 &\quad\hbox{if } x\in S_- \\
\end{matrix}
\right .  \ .
\label{HH}\ee
Then our approximation $u_1(x)$ approaches $\HH(x)$ at an exponential rate $O( e^{-\sigma_{\pm} |t|})$ as $|t|$ increases.
 The global approximation we will use  consists simply of interpolating $u_1$ with $\HH$ sufficiently well-inside $\R^3\setminus M_\A$  through a cut-off in $|t|$.  In order
to avoid the problem described in \S \ref{diverg} and  having the coordinates $(y,t)$ well-defined, we consider this cut-off to be supported in a region $y$-dependent that expands logarithmically in $r_\A$.  Thus we will actually consider a region $\NN_\delta$ expanding at the ends, thus
becoming wider as $r_\A\to \infty$ than the set ${\mathcal N}_\delta^\A$ previously considered, where the coordinates
are still well-defined.

\medskip
We consider the open set $\OO$ in $M_\A \times \R$ defined as
\be
\OO\, =\, \{\, (y,t)\in M_\A\times \R,\quad |t + h_1(\A y)| <  \frac \delta\A  \, +  \,  4\, \max\, \{\sigma^{-1}_-, \sigma^{-1}_+\}  \log (1 + r_\A(y))=: \rho_\A(y) \, \}
\label{region}\ee
where $\delta$ is small positive number. We consider the the region $\NN=:\NN_\delta$ of points $x$ of the form
$$
x=X_h(y,t) =  y + (t+ h_0(\A y) + h_1(\A y))\, \nu(\A y),\quad (y,t)\in \OO,
$$
namely $\NN_\delta = X_h (\OO)$.
The coordinates $(y,t)$ are well-defined in $\NN_\delta$ for any sufficiently small $\delta$: indeed the map $X_h$ is one to one in $\OO$
  thanks to assumption \equ{condbeta} and the
fact that $h_1 = O(\A)$.  Moreover, Lemma \ref{lapfinal} applies in $\NN_\delta$.

\bigskip
Let $\eta(s)$ be a smooth cut-off function with $\eta(s) =1$ for $s<1$ and $=0$ for $s>2$.
and define
\be
\eta_\delta (x) \, := \,
\left\{
\begin{matrix}
\ \eta( \, |t + h_1(\A y)| - \rho_\A(y)-3)  &  \quad\hbox{if } x\in \NN_\delta \, , \\
 0 &\quad\hbox{if } x \not \in  \NN_\delta
\end{matrix}
\right.
\label{etadelta}\ee
where $\rho_\A$ is defined in \equ{region}.
 Then we let our global approximation $\ww(x)$ be simply defined as
\be
\ww  \, := \, \eta_\delta u_1  + (1-\eta_\delta )\HH
\label{global}\ee
where  $\HH$ is given by \equ{HH} and $u_1(x)$ is just understood to be $\HH(x)$  outside $\NN_\delta$.

\medskip
Since $\HH$ is an exact solution in $\R^3\setminus M_\delta$, the global error of approximation is simply computed as

\be
S(\ww) \  =\ \Delta \ww + f(\ww)\ =   \ \eta_\delta  S(u_1)\ \  + E
\label{error}\ee
where
$$
E=  2\nabla \eta_\delta\nabla u_1  + \Delta \eta_\delta (u_1 -\HH) \,  +\,   f(\eta_\delta u_1  + (1-\eta_\delta )\HH)\,  )\, -\, \eta_\delta f(u_1) \, .
$$

\medskip
The new error terms created are of exponentially small size $O( e^{-\frac \sigma \A})$ but have in addition decay with $r_\A$. In fact we have
$$
|E|\le Ce^{-\frac \delta \A}\, r_\A^{-4} .
$$
Let us observe that  $ |t + h_1(\A y)|  = |z-h_0(\A y)|$ where $z$ is the normal coordinate to $M_\A$, hence $\eta_\delta$ does not depend
on $h_1$, in particular  the term
$\Delta \eta_\delta$ does involves second derivatives of $h_1$ on which  we have not made assumptions yet.


\setcounter{equation}{0}
\section{ \emph{The proof of Theorem \ref{teo1}}}

The proof of Theorem \ref{teo1}  involves  various ingredients whose  detailed proofs are fairly technical. In order to keep the presentation
as clear as possible,  in this section we carry out the proof, skimming it from several (important) steps, which we state as lemmas or propositions, with complete proofs  postponed for the subsequent sections.

\medskip
We look for a solution $u$ of the Allen Cahn equation \equ{ac} in the form
\be
u = \ww + \vp
\ee
where $\ww$ is the global approximation defined in \equ{global} and $\vp$ is in some suitable sense small.
Thus we need to solve the following problem
\be
\Delta\vp + f'(\ww)\vp = -S(\ww)  -N(\vp)
\label{nl}\ee
where $$N(\vp) = f(\ww +\vp )- f(\ww) - f'(\ww)\vp .$$

\medskip
Next we introduce various  norms that we will use to set up a suitable functional analytic scheme for solving problem \equ{nl}.
For a function $g(x)$ defined in $\R^3$, $1<p\le +\infty$,  $\mu>0$, and $\A>0$ we write
$$
\| g\|_{p,\mu,*} := \sup_{x\in \R^3} (1+ r(\A x))^\mu \|g \|_{L^p(B(x,1))} ,\quad r(x',x_3) = |x'|\ .
$$

On the other hand, given numbers $\mu\ge 0$, $0< \sigma< \min\{\sigma_+, \sigma_-\}$, $p>3$,  and
functions $g(y,t)$ and $\phi(y,t)$ defined in $M_\A \times \R$
we consider the norms
\be
\|g\|_{p,\mu,\sigma}\, :=\, \sup_{(y,t)\in M_\A\times \R}  r_\A(y)^\mu\, e^{\sigma|t| }  \left ( \int_{B( (y,t),1)} |f|^p\, dV_\A \right )^{\frac 1p} .
\label{normp}\ee
Consistently we set
\be
\|g\|_{\infty,\mu,\sigma}\, :=\,  \sup_{(y,t)\in M_\A\times \R} r_\A(y)^\mu\, e^{\sigma|t| } \, \|f\|_{L^\infty ( B( (y,t),1) )}
\label{norminfty}\ee
 and  let
 \be
\|\phi\|_{2,p,\mu,\sigma} := \|D^2\phi\|_{p,\mu,\sigma}+ \|D\phi\|_{\infty,\mu,\sigma}+
 \|\phi\|_{\infty,\mu,\sigma} \, .
\label{norm2p}\ee
We consider also for a function $g(y)$ defined in $M$ the  $L^p$-weighted norm

\begin{equation}
\|f\|_{p,\beta}\, :=\,  \left(  \int_{M }
|f(y) |^p\,(1+ |y|^\beta\,)^{p } \, d V(y )\, \right )^{1/p} \ = \  \| \,(1+ |y|^{\beta})\, f \, \|_{L^p(M )}
\label{*p}\end{equation}
where $p>1$ and $ \beta>0$.

\medskip
We assume in what follows,  that for a certain constant $\KK>0$  and $p>3$ we have that the parameter function $h_1(y)$ satisfies
\be
\|h_1\|_*:= \| h_1\|_{L^\infty (M)} + \|(1+ r^2)D h_1\|_{L^\infty (M)} +  \| D^2 h_1\|_{p, 4-\frac 4p}
 \ \le \ \KK \A\  .
\label{ass1}\ee

\medskip
 Next we reduce problem \equ{nl} to solving one qualitatively similar (equation \equ{e4} below) for a function $\phi(y,t)$ defined in the whole space  $M_\A \times \R$.


\subsection{Step 1: the gluing reduction}\label{step0}  We will follow the following procedure.
Let us consider
again $\eta(s)$, a smooth cut-off function with $\eta(s) =1$ for $s<1$ and $=0$ for $s>2$,
and define
\be
\zeta_n (x) \, := \,
\left\{
\begin{matrix}
\ \eta( \, |t+h_1(\A y) | - \frac \delta\A +n )  &  \quad\hbox{if } x\in \NN_\delta  \\
 0 &\quad\hbox{if } x \not \in  \NN_\delta
\end{matrix}
\right .  \ .
\label{zetan}\ee
We look for a solution $\vp(x)$ of problem \equ{nl} of the following form
\be
\vp(x) = \zeta_2(x)\phi(y,t) + \psi(x)
\label{formphi}
\ee
where  $\phi$ is defined in entire $M_\A\times \R$, $\psi(x)$ is defined in $\R^3$ and
$\zeta_2(x)\phi(y,t)$ is understood as zero outside $\NN_\delta$.

We compute, using that $\zeta_2\cdot\zeta_1 = \zeta_1$,
$$
S(\ww +\vp)\, =\, \Delta \vp +f'(\ww )\vp  + N(\vp) + S(\ww)  =
$$

$$
   \zeta_2\,\left [\, \Delta\phi + f'(u_1)\phi \, + \, \zeta_1 (f'(u_1) + H(t)) \psi\,  +
\zeta_1 N(\psi + \phi) +   S(u_1)  \, \right ] \  +
$$

$$\Delta\psi -
[\, (1-\zeta_1) f'(u_1)  + \zeta_1 H(t)\, ] \psi\  +
$$

\be
(1-\zeta_2)S(\ww) + (1-\zeta_1)N(\psi + \zeta_2\phi)  + 2\nn\zeta_1\nn\phi + \phi \Delta \zeta_1
\label{dd}\ee

where $H(t)$ is any smooth, strictly negative function satisfying

$$
H (t) \, = \,
\left\{
\begin{matrix}
\  f'(+1)  &  \quad\hbox{if }  t> 1\, , \   \\
 f'(-1) &\quad\hbox{if }  t< -1
\end{matrix}
\right .  \ .
$$
Thus, we will have constructed a solution $\vp= \zeta_2\phi +\psi$ to problem \equ{nl} if we require that
the pair $(\phi,\psi)$ satisfies the following coupled system

\be
 \Delta\phi + f'(u_1)\phi \, + \, \zeta_1 (f'(u_1) - H(t)) \psi\,  +
\zeta_1 N(\psi + \phi) +   S(u_1) = 0 \ \hbox{for } |t|< \frac \delta\A + 3
\label{e1}\ee

$$
 \Delta\psi +
[\, (1-\zeta_1) f'(u_1)  + \zeta_1 H(t)\, ] \psi\  +
$$

\be
(1-\zeta_2)S(\ww) + (1-\zeta_1)N(\psi + \zeta_2\phi)  + 2\nn\zeta_1\nn\phi + \phi \Delta \zeta_1 \,=\, 0 \quad\hbox{ in } \R^3\, .
 \label{e2}\ee

\medskip
In order to find a solution to this system we will first extend equation \equ{e1} to entire $M_\A\times \R$
in the following manner. Let us set

\be
\BB (\phi) =   \zeta_4 [\Delta_x - \partial_{tt} - \Delta_{y,M_\A}\,]\, \phi
\label{BB}\ee

where $\Delta_x$ is expressed in $(y,t)$ coordinates using expression \equ{laplacianfinal} and  $\BB(\phi)$ is understood to be zero
for $|t+h_1| > \frac \delta\A + 5$. The other terms in equation \equ{e1} 
are simply extended as zero beyond the support of $\zeta_1$.
Thus we consider the extension of equation \equ{e1} given by

$$
\partial_{tt}\phi\, + \, \Delta_{y,M_\A}\phi\, + \BB(\phi) \, +\, f'(w(t))\phi = - \ttt S(u_1)
$$

\be
 - \left \{ [f'(u_1)- f'(w)]\phi \, + \, \zeta_1 (f'(u_1) - H(t)) \psi\,  +
\zeta_1 N(\psi + \phi)\right \} \hbox{ in }\in M_\A\times \R,
\label{e3}\ee
where we set, with reference to expression \equ{error1},

\medskip
$$
 \ttt S(u_1) =
 - \A^2[ |A|^2h_1  + \Delta_M h_1 ]\,w'  +
\A^2  \,a_{ij}^0
\, ( 2\partial_ih_0 \partial_j h_1 + \partial_ih_1 \partial_j h_1\,) \, w''
$$
\be
 - \A^4 \, (\,|A|^2 t\psi_0' - a_{ij}^0 \pp_i h_0\pp_j h_0\,t \psi_1'\, )\, \Delta h_1\, +
\zeta_4 \,[\, \A^3(t+ h) a_{ij}^1 \, \partial_i h \partial_jh \, w'' +
R_1(y,t)\,]
\label{error2}
\ee

\medskip
\noindent
and, we recall
$$R_1= R_1(y,t,h_1 (\A y) , \nn_M h_1 (\A y) )$$ with
\be
 | D_{\imath} R_1(y, t, \imath, \jmath)| + | D_{\jmath} R_1(y, t, \imath, \jmath)| +| R_1(y,t, \imath, \jmath)| \le C \A^3(1+ r_\A(y))^{-4} e^{-\sigma|t|} .
\label{R11}\ee

\medskip
In summary $\ttt S(u_1)$ coincides with $S(u_1)$ if  $\zeta_4=1 $ while outside the support of $\zeta_4$,
their parts that are not defined for all $t$
are cut-off.

\medskip
To solve the resulting system \equ{e2}-\equ{e3}, we find first solve equation \equ{e2} in $\psi$ for a given $\phi$ a
small function in absolute value. Noticing that the potential
$[\, (1-\zeta_1) f'(u_1)  + \zeta_1 H(t)\, ]$ is uniformly negative, so that the linear operator is qualitatively like $\Delta -1$ and using
contraction mapping principle, a solution $\psi= \Psi(\phi)$  is found according to the following lemma, whose detailed proof we carry out in
\S \ref{lpsi1}.

\begin{lemma}\label{lemapsi1}
For all sufficiently small $\A$ the following holds.
Given $\phi $ with \\ $\|\phi\|_{2,p,\mu,\sigma}\le 1$,
there exists a unique solution $\psi= \Psi (\phi)$ of problem $\equ{e2}$ such that
\be \|\psi\|_X := \|D^2\psi\|_{p,\mu,*} + \|\psi\|_{p,\mu,*} \, \le\, Ce^{-\frac {\sigma \delta}{\A}} . \label{x-2}\ee
Besides, $\Psi$ satisfies the Lipschitz condition

\be
\|\Psi(\phi_1) - \Psi(\phi_2) \|_X \,\le\, C\,e^{ -\frac {\sigma \delta}{\A}} \|\phi_1 -\phi_2\|_{2,p,\mu,\sigma} \ .
\label{lipspsi}\ee

 \end{lemma}

Thus we  replace  $\psi= \Psi(\phi)$  in  the first equation \equ{e1} so that by setting
\be
\NNN(\phi) := \BB(\phi) + [f'(u_1)- f'(w)]\phi \, + \, \zeta_1 (f'(u_1) - H(t)) \Psi(\phi)\,  +
\zeta_1 N( \Psi(\phi) + \phi),
\label{NNN}\ee
our problem is reduced to finding a solution $\phi$  to the following nonlinear, nonlocal problem in $M_\A\times \R$.

\be
\partial_{tt}\phi\, + \, \Delta_{y,M_\A}\phi\, + f'(w)\phi\, = \, -\tilde S(u_1) - \NNN(\phi)
\quad \hbox{ in } M_\A\times \R.
\label{e4}\ee

\medskip
Thus, we concentrate in the remaining of the proof in solving equation \equ{e4}.
As we hinted in \S \ref{u1}, we will find a solution of problem \equ{e4} by considering two steps:
(1) ``Improving the approximation'', roughly solving for
$\phi$ that eliminates the part of the error that does not contribute to the ``projections''\  $\int [\ttt S(U_1 ) + \NNN(\phi) ]w'(t) dt$,\  which amounts to a nonlinear problem in $\phi$,
 and (2) Adjust $h_1$ in such a way
that the resulting projection is actually zero.
Let us set up the scheme for step (1) in a precise form.

\subsection{ Step 2:  Eliminating  terms not contributing to projections}\label{step1}
Let us consider the problem of finding a function $\phi(y,t)$ such that for a certain function
$c(y)$ defined in $ M_\A$, we have

\begin{align}
\begin{aligned}
\partial_{tt}\phi\, + \, \Delta_{y,M_\A}\phi\,& = \, - \ttt S(u_1) - \NNN(\phi)\,  +\, c (y) w'(t)
\quad \hbox{ in } M_\A\times \R
,\\
\int_\R \phi(y,t)\,w'(t)\,dt&=0,  \foral  y\in  M_\alpha \, .
\end{aligned}
\label{nonlinear0}
\end{align}
Solving this problem for $\phi$ amounts to ``eliminating the part of the error that does not contribute to the projection'' in problem \equ{e4}. To justify
this phrase let us consider the associated linear problem in $M_\A\times \R$

\begin{align}
\begin{aligned}
 \partial_{tt}\phi + \Delta_{y,M_\A} \phi   +f'(w(t))\phi&=  g(y,t)  +  c (y) w'(t)  , \foral (y,t) \in  M_\A\times \R
,\\
\int_{-\infty} ^ {\infty} \phi(y,t)\,w'(t)\,dt&=0,  \foral  y\in  M_\A \, .
\end{aligned}
\label{p1}
\end{align}

\medskip
Assuming that the corresponding operations can be carried out, let us multiply the equation by $w'(t)$ and integrate in $t$ for fixed $y$. We find
that
$$
\Delta_{y,M_\A} \int_\R \phi(y,t)\,w'\,dt + \int_\R \phi(y,t)\, [ w''' + f'(w)w' ]\, dt = \int_\R g\,w' + c(y) \int_\R {w'}^2\, .
$$
The left hand side of the above identity is zero and then we find that
\be
c(y) =  - \frac {\int_\R g(y,t ) w'dt }{\int_\R {w'}^2dt } \, ,
\label{cy}\ee
hence a $\phi$ solving problem \equ{p1}.
$\phi$ {\em precisely}  solves or {\em eliminates} the part of $g$ which does not contribute to the projections
in the equation $\Delta\phi   +f'(w)\phi =  g$, namely the same equation with $g$ replaced by $\ttt g$ given by
\be
 \tilde g(y,t) =  g(y,t) -   \frac {\int_\R f(y, \cdot ) w' }{\int_\R {w'}^2 } \, w'(t)\, .
\label{io}\ee
 The term
$c(y)$ in problem \equ{nonlinear0}  has a similar role, except that we cannot find it so explicitly.

\medskip

In order to solve problem \equ{nonlinear0} we need to devise a  theory to solve problem \equ{p1}
where we consider a class of right hand sides $g$ with a qualitative behavior similar to that of the error
$S(u_1)$. As we have seen in  \equ{error2}, typical elements in this error are  of the type
$O((1+ r_\A(y))^{-\mu} e^{-\sigma|t|} )$, so this is the type of functions $g(y,t)$ that we want to consider. This is actually the motivation to introduce the norms \equ{normp}, \equ{norminfty} and \equ{norm2p}.
We will prove that problem \equ{p1} has a unique solution $\phi$ which respects the size of $g$ in norm \equ{normp}
up to its second derivatives, namely in the norm \equ{norm2p}. The following fact holds.

\begin{prop}\label{prop1}
Given $p>3$, $\mu \ge 0$ and $0< \sigma < \min\{\sigma_-, \sigma_+\}$,
there exists a constant $C>0$ such that for all sufficiently small $\A >0$ the following holds.
Given $f$ with $\|g\|_{p,\mu,\sigma}< +\infty$, then Problem $\equ{p1}$ with $c(y)$ given by $\equ{cy}$, has a unique  solution $\phi$ with
$\|\phi \|_{\infty,\mu,\sigma} < +\infty$.
This  solution
satisfies in addition that
\be
\|\phi \|_{2,p,\mu,\sigma}\, \le \, C \|g\|_{p,\mu,\sigma}\ .
\label{cota}
\ee
\end{prop}

We will prove this result in \S \ref{linear0} . After Proposition
\ref{prop1}, solving Problem \equ{nonlinear0} for a small $\phi$ is
easy using the small Lipschitz character of the terms involved in
the operator $\NNN(\phi)$ in \equ{NNN} and contraction mapping
principle.  The error term $\ttt S(u_1)$
satisfies
\be \|\ttt S(u_1) + \A^2 \Delta h_1  w' \|_{p,4,\sigma}  \le C\A^3 .
\label{errorsize}\ee
 Using this, and the fact that $\NNN(\phi)$ defines a contraction mapping in a ball center zero and radius
$O(\A^3)$ in  $\|\ \|_{2,p,4,\sigma}$, we conclude the existence of a unique small solution $\phi$ to problem \equ{nonlinear0}
whose size is $O(\A^3)$  for this norm.
This solution $\phi$ turns out to define an operator in $h_1$ $\phi =\Phi(h_1)$
which is Lipschitz in the norms $\|\ \|_*$ appearing in condition \equ{ass1}. In precise terms, we have the validity of the following result, whose detailed
proof we postpone for \S \ref{pr2}.

\begin{prop}\label{prop2} Assume $p>3$, $0\le \mu\le 3$, $0< \sigma < \min\{ \sigma_+,\sigma_-\}$.
There exists a $K>0$ such that problem $\equ{nonlinear1}$ has a unique solution $\phi =\Phi(h_1)$ such that
$$
 \|\phi\|_{2,p,\mu,\sigma}\, \le K \A^3\, .
 $$
Besides, $\Phi$ has a Lipschitz dependence on $h_1$ satisfying $\equ{ass1}$ in the sense that
 \be
 \|\Phi (h_1)- \Phi(h_2)  \|_{2,p,\mu,\sigma} \le C \A^2 \|h_1-h_2\|_{*} .
 \label{rw}\ee
\end{prop}

\subsection{ Step 3: Adjusting $h_1$ to make the projection zero}\label{step2}
In order to conclude the proof of the theorem, we have to carry out the second step, namely adjusting $h_1$, within a region of the form
\equ{ass1} for suitable $\KK$ in such a way that the ``projections'' are identically zero, namely making zero the function $c(y)$ found for the solution $\phi=\Phi(h_1)$ of  problem \equ{nonlinear0}.
Using expression \equ{cy} for $c(y)$ we find that

\be
 c(y)\int_\R {w'}^2  = \int_\R \ttt S(u_1) w'\, dt  + \int_{\R}\NNN(\Phi(h_1)\,)\, w'\, dt \, .
\label{cyy}\ee

Now, setting  $c_*:= \int_\R  {w'}^2 dt$ and using same computation employed to derive formula \equ{projection}, we find from expression \equ{error2} that

$$
\int_\R   \ttt S(u_1)(y,t)\, w'(t)\, dt\ = \,
   - c_*\, \A^2 (\Delta_M h_1 + h_1|A|^2)\, +\,  c_* \A^2 G_1(h_1)
  $$
  where
  $$
   c_* G_1(h_1) \, =\,
 - \A^2 \,\Delta h_1\, (\,|A|^2 \int_\R   t\psi_0'w'\, dt  - a_{ij}^0 \pp_i h_0\pp_j h_0\,\int_\R t \psi_1'w'\, dt \, ) \ +
 $$
 \be
\, \A  \, \partial_i (h_0+ h_1) \partial_j (h_0+ h_1) \,  \int_\R \zeta_4 (t+ h)a_{ij}^1w''w'\, dt\,  +
\A^{-2} \int_\R \zeta_4\, R_1(y,t, h_1, \nn_M h_1\, )\, w'\, dt
\label{G1}
\ee
and we recall that $R_1$ is of size $O(\A^3)$ in the sense \equ{R11}. Thus, setting
\be
c_* G_2(h_1)\, :=\, \A^{-2} \int_{\R}\NNN(\Phi(h_1)\,)\, w'\, dt, \quad G(h_1):= G_1(h_1) + G_2(h_1),
\label{G}\ee
we find that the equation  $c(y) =0$ is equivalent to the problem
\be
{\mathcal J}(h_1) = \Delta_M h_1 + |A|^2 h_1 =  G(h_1) \quad\hbox{ in } M.
\label{probjac1}\ee

Therefore, we will have proven Theorem \ref{teo1} if we find a function $h_1$ defined on $M$ satisfying constraint \equ{ass1}
for a suitable $\KK$ that solves equation \equ{probjac1}. Again, this is not so direct  since the operator $\JJ$ has a nontrivial bounded kernel.
Rather than solving directly \equ{probjac1}, we consider first a projected version of this problem, namely that of finding $h_1$ such that
for certain scalars $c_1,\ldots, c_J$ we have
$$
\JJ(h_1) = G(h_1)   + \sum_{i=1}^J  \frac {c_i}{1+ r^4}\, \hat{z}_i\quad\hbox{ in } M,
$$
\be
\int_M  \frac {\hat{z}_ih} {1+ r^4}\, dV\, =\, 0, \quad i=1,\ldots J.
\label{probjac2}\ee
Here $\hat{z}_1, ..., \hat{z}_J$ is a basis of the vector space of bounded  Jacobi fields.

\medskip
In order to solve problem \equ{probjac2} we need a corresponding linear invertibility theory.   This leads us to consider
the linear problem
\be
\JJ(h) = f  + \sum_{i=1}^J  \frac {c_i}{1+ r^4}\, \hat{z}_i\quad\hbox{ in } M,
\nonumber \ee
\be
\int_M  \frac {\hat{z}_ih} {1+ r^4}\, dV\, =\, 0, \quad i=1,\ldots J.
\label{projac2}\ee
Here $\hat{z}_1, ..., \hat{z}_J$ are bounded, linearly independent Jacobi fields, and $J$ is the dimension of the vector space of bounded Jacobi fields. 

\medskip
We will prove in \S \ref{proi} the following result.

\begin{prop}\label{proj} Given $p>2$ and
$f$ with $\|f\|_{{p,4-\frac 4p}} <+\infty$, there exists a unique bounded solution $h$ of problem
$\equ{projac2}$.  Moreover, there exists a positive number $C= C(p,M)$ such that
\begin{equation}
\|h\|_*:= \|h\|_\infty \, + \, \|\,(1+|y|^2)\, D h\|_\infty \, + \, \| \, D^2h\, \|_{p,4-\frac 4p}\, \le\, C\|f\|_{{p,4-\frac 4p}} \, .
\label{est}\end{equation}
\end{prop}

 Using the fact that $G$  is a small operator of size $O(\A)$ uniformly on functions $h_1$ satisfying \equ{ass1}, Proposition \ref{proj} and
 contraction mapping principle yield the following result, whose detailed proof we carry out in \S \ref{nonljac}.

\begin{prop}\label{proii} Given $p>3$, there exists a number $\KK>0$ such that for all sufficiently small $\A>0$ there is
 a unique solution $h_1$ of problem
$\equ{probjac2}$ that satisfies constraint $\equ{ass1}$.
\end{prop}

\subsection{Step 3: Conclusion}
At the last step we prove that the constants $c_i$ found in equation
\equ{probjac2} are in reality all zero, without the need of
adjusting any further parameters but rather as a consequence of the
natural invariances of the of the full equation. The key point is to realize what equation has been solved so far.

First we observe the following.  For each $h_1$ satysfying \equ{ass1}, the pair $(\phi,\psi)$ with $\phi =\Phi(h_1)$, $\psi = \Psi(\phi)$,
 solves the system
\be
 \Delta\phi + f'(u_1)\phi \, + \, \zeta_1 (f'(u_1) - H(t)) \psi\,  +
\zeta_1 N(\psi + \phi) +   S(u_1) =  c(y)w'(t) \ \hbox{for } |t|< \frac \delta\A + 3
\nonumber \ee

$$
 \Delta\psi +
[\, (1-\zeta_1) f'(u_1)  + \zeta_1 H(t)\, ] \psi\  +
$$

\be
(1-\zeta_2)S(\ww) + (1-\zeta_1)N(\psi + \zeta_2\phi)  + 2\nn\zeta_1\nn\phi + \phi \Delta \zeta_1 \,=\, 0 \quad\hbox{ in } \R^3\, .
 \nonumber \ee
Thus setting
\be
\vp(x) = \zeta_2(x)\phi(y,t) + \psi(x), \quad u = \ww + \vp\, ,
\nonumber
\ee
we find from formula \equ{dd} that

$$
\Delta u + f(u) = S(\ww +\vp)\, =\,  \zeta_2 c(y)\, w'(t)\ .
$$
On the other hand choosing $h_1$ as that given in Proposition \ref{proii} which solves problem \equ{probjac2}, amounts precisely to
making
$$
c(y) =  c_*\A^2  \sum_{i=1}^J c_i \frac {\hat z_i ( \A y)} {1+ r_\A(y)^4}
$$
for certain scalars $c_i$. In summary, we have found $h_1$ satisfying constraint \equ{ass1} such that
\be u=  \ww +  \zeta_2(x)\Phi(h_1) + \Psi( \Phi(h_1)\, ) \label{uu}\ee
solves the equation
 \be \Delta u + f(u) = \sum_{j=1}^J
\frac{\ttt c_i}{1+ r_\A^4}\, \hat{z}_i(\A y) w'(t) \label{equ}\ee
where $\ttt c_i = c_*\A^2  c_i$.  Testing
equation \equ{equ} against the generators of the rigid motions $\partial_i u$
$i=1,2,3$,  $-x_2\partial_1u + x_1\partial_2u$, and using the
balancing formula for the minimal surface and the zero average of
the numbers $\beta_j$ in the definition of $h_0$, we find a system
of equations that leads us to $c_i=0$ for all $i$, thus conclude
the proof. We will carry out the details in \S \ref{conclusion}.

\medskip
In sections \S \ref{linear0}-\ref{conclusion} we will complete the proofs of the intermediate steps of the program designed in this section.


\setcounter{equation}{0}
\section{   \emph{The linearized operator}}\label{linear0}
In this section we will prove Proposition \ref{prop1}. 
At the core of the proof of the stated a priori estimates is the fact that the one-variable solution
$w$ of \equ{ac} is {\em nondegenerate} in  $L^\infty(\R^3)$ in the sense that
 the linearized operator
$$
L(\phi) = \Delta_y\phi + \partial_{tt}\phi + f'(w(t))\phi ,\quad (y,t)\in \R^3=\R^2\times\R,
$$
is such that the following property holds.
\begin{lemma}\label{lemma l1}
Let $\phi$ be a bounded, smooth solution of the problem
\begin{equation}
L(\phi) = 0 \quad\hbox{in } \R^2\times \R.
\label{l3}\end{equation}
Then $\phi(y,t) = Cw'(t) $ for some $C\in \R$.
\end{lemma}

\proof
We begin by  reviewing some known facts about the one-dimensional operator
 $
L_0(\psi)= \psi'' + f'(w)\psi.
$
Assuming that $\psi(t)$ and its derivative decay sufficiently fast as $|t|\to +\infty$ and
defining  $\psi(t) = w'(t)\rho(t) $,
we get that
$$
\int_\R [|\psi'|^2 - f'(w)\psi^2]\,dt\, =\, \int_\R L_0(\psi) \psi \, dt \, =\,  \int_\R  {w'}^2|\rho '|^2\, dt ,
$$
therefore this quadratic form is positive unless $\psi$ is a
constant multiple of $w'$. Using this and a standard compactness
argument we get that there is a constant $\gamma >0$ such that
whenever $\int_\R \psi w' = 0$ with $\psi\in H^1(\R)$ we have that
\be \int_\R (\, |\psi'|^2 - f'(w)\psi^2\,)\, dt\, \ge \, \gamma
\int_\R (\, |\psi'|^2 + |\psi|^2\,)\, dt . \label{q}\ee Now, let
$\phi$ be a bounded solution of equation \equ{l3}. We claim that
$\phi$ has exponential decay in $t$, uniform in $y$. Let us consider
a small number $\sigma>0$ so that for a certain  $t_0 >0$  and all
$|t| > t_0$ we have that
$$
f'(w) <  - 2\sigma^2 .
$$
Let us consider for $\ve>0$  the function
$$
g_\ve (t,y) =   e^{-\sigma(|t|-t_0)}  + \ve  \sum_{i=1}^2 \cosh( \sigma y_i)
$$
Then for $|t|>t_0$ we get that
$$
L(g_\delta ) < 0 \quad \hbox{if } |t| >t_0 .
$$
As a conclusion, using maximum principle, we get
$$
|\phi| \le \|\phi\|_\infty \, g_\ve \quad \hbox{if } |t| >t_0 ,
$$
and letting $\ve \to 0$ we then get
$$
|\phi(y,t)| \ \le\ C\|\phi\|_\infty   e^{-\sigma|t|}  \quad \hbox{if } |t| >t_0 \ .
$$
Let us  observe the following fact: the function
$$
\ttt \phi(y,t) =
\phi(y,t) - \left ( \int_\R w' (\zeta)\, \phi(y,\zeta)\, d \zeta \right )\, \frac{ w'(t)} {\int_\R {w'}^2}
$$
also satisfies $L(\ttt \phi ) = 0$ and, in addition,
\be
\int_\R   w' (t)\, \tilde \phi(y,t)\, d t  = 0\foral y\in \R^2 .
\label{orti}\ee
In view of the above discussion, it turns out that the function
$$
\vp (y) := \int_\R \ttt \phi^2(y,t)\, dt
$$
is well defined. In fact so are its first and second derivatives by elliptic regularity of $\phi$, and differentiation under the
integral sign is thus justified. Now, let us observe that
$$
\Delta_y \vp (y) = 2 \int_\R \Delta_y \ttt \phi \cdot \ttt \phi \, dt  + 2\int_\R |\nabla_y \ttt \phi |^2
$$
and hence
\begin{align}
\begin{aligned}
0 &= \int_\R  (L(\ttt \phi) \cdot \ttt \phi )\\
& =
\frac{1}{2}\Delta_y \vp
- \int_\R  |\nabla_y \ttt \phi |^2
 \, dz  -
\int_\R (\, |\ttt \phi _t|^2 - f'(w)\ttt \phi ^2\,)\, dt\, .
\end{aligned}
\label{ew}
\end{align}

Let us observe that because of relations \equ{orti} and \equ{q}, we have that
$$
\int_\R (\, |\ttt\phi _t|^2 - f'(w)\ttt \phi ^2\,)\, dt\, \ge \gamma \vp .
$$
It follows then that

$$\frac{1}{2}\Delta_y \vp   -\gamma \vp   \ge 0. $$
Since $\vp$ is bounded, from maximum principle we find that
$\vp$ must be identically equal to zero.
But this means
\be
\phi(y,t) = \left ( \int_\R w' (\zeta)\, \phi(y,\zeta)\, d \zeta \right )\, \frac{ w'(t)} {\int_\R {w'}^2}. \label{uuu}\ee
Then the bounded function
$$
g(y) = \int_\R w_\zeta (\zeta)\, \phi(y,\zeta)\, d \zeta
$$
satisfies the equation
\be
\Delta_y g = 0, \quad \mbox{in}\ \R^2.
\label{f}\ee
Liouville's theorem  implies that $g\equiv$ constant and relation \equ{uuu}
yields
 $\phi(y,t) = Cw'(t)$ for some $C$.
This concludes the proof.\qed

\bigskip

\subsection{ A priori estimates}
We shall consider problem \equ{p1} in a slightly more general form, also in a domain finite in $y$-direction.
For a large number $R>0$ let us set
$$
M_\A^R := \{y\in M_\A \ /\ r(\A y) < R\}
$$
and consider the variation of Problem \equ{p1} given by

\begin{align}
\begin{aligned}
 \partial_{tt}\phi + \Delta_{y,M_\A} \phi   +f'(w(t))\phi &=  g(y,t)  +  c (y) w'(t)   \quad\hbox{in }   M_\A^R\times \R
,\\
 \phi\,&=0,  \quad\hbox{on }  \partial M_\A^R\times \R,\\
\int_{-\infty} ^ {\infty} \phi(y,t)\,w'(t)\,dt&=0  \foral  y\in  M_\A^R,
\end{aligned}
\label{p11-1}
\end{align}
where we allow $R= +\infty$ and
$$
c(y) \int_\R {w'}^2 dt =  - \int_\R  g(y,t) \,w'\, dt \ .
$$

\medskip
We begin by proving  a priori estimates.

\begin{lemma}\label{apriori}

 Let us assume that  $0<\sigma< \min\{\sigma_-,\sigma_+\}$ and $\mu\ge 0$.  Then there exists a constant $C>0$ such that for all small $\A$ and all large $R$, and every solution $\phi$ to Problem $\equ{p11}$ with
 $\|\phi\|_{\infty, \mu, \sigma} <+\infty$
and right hand side $g$ satisfying
 $\|g\|_{p,\mu,\sigma} < +\infty $ we have
\be
\| D^2 \phi\|_{p,\mu,\sigma } + \|D\phi\|_{\infty,\mu,\sigma} + \|\phi\|_{\infty,\mu,\sigma } \le  C\|g\|_{p,\mu,,\sigma}.
\label{1.1}\ee
\end{lemma}

\proof For the purpose of the a priori estimate, it clearly suffices to consider 
the case $c(y)\equiv 0$.
By local elliptic estimates, it is enough to show that
\be
\|\phi\|_{\infty,\mu,\sigma } \le   C\|g\|_{p,\mu,\sigma}.
\label{2}\ee
Let us assume by contradiction that \equ{2} does not hold. Then we have sequences $\A=\A_n\to 0$, $R=R_n\to \infty$,
$g_n$ 
with  $\|g_n\|_{p,\mu,\sigma}\to 0$, $\phi_n$ with $\|\phi_n\|_{\infty,\mu,\sigma} =1$  such that

\begin{align}
\begin{aligned}
 \partial_{tt}\phi_n + \Delta_{y,M_\A} \phi_n   +f'(w(t))\phi_n  &=  g_n     \quad\hbox{in }   M_\A^R\times \R
,\\
 \phi_n\,&=0  \quad\hbox{on }  \partial M_\A^R\times \R,\\
\int_{-\infty} ^ {\infty} \phi_n(y,t)\,w'(t)\,dt&=0  \foral  y\in  M_\A^R \, .
\end{aligned}
\label{p12}
\end{align}

Then we can find points $(y_n,t_n)\in  M_\A^R\times \R$ such that

$$
 e^{-\sigma|t_n|}(1+ r(\A_n y_n))^\mu \,  |\phi_n(y_n,t_n)| \ge \frac 12.
 $$

We will consider different possibilities.  
We may assume that either   $ r_\A(y_n) = O(1)$ or  $r_\A(y_n)\to +\infty$.

\subsubsection{Case  $r(\A_n y_n )$  bounded.}
We have  $\A_n y_n$ lies within a
bounded subregion of $M$, so we may assume that
$$ \A_n y_n \to \tilde y_0\in M .$$
Assume that
$\tilde y_0  \in Y_k( {\mathcal U}_k)$ for one of the local parametrization of $M$.  
We consider $\ttt \py _n, \ttt \py_0\in {\mathcal U}_k$ with $Y_k( \ttt \py_n )= \A_n y_n,  $ $Y_k( \ttt \py_0 )= \ttt y_0$.

\medskip
On $\A_n^{-1} Y_k( {\mathcal U}_k)$, $M_\A$ is parameterized by
$Y_{k,\A_n} ( \py ) = \A_n^{-1} Y_k( \A_n \py) $, $\py \in \A_n^{-1}{\mathcal U}_k$. Let us consider
the local change of variable,
$$
\py =  \A^{-1} \ttt \py_n + \by.
$$

\subsubsection{ Subcase $t_n$ bounded} 
 Let us assume first that
$ |t_n| \le C .$
Then, setting $$\ttt \phi_n(\by,t ):= \ttt \phi_n (\A^{-1} \ttt \py_n + \by, t),$$
the local equation becomes
$$
a_{ij}^0 (\ttt \py_n + \A_n \by )\pp_{ij}\ttt \phi_n
 + \A_nb_j^0 (\ttt \py_n + \A_n \by )\pp_j\ttt \phi_n   +  \partial_{tt}\ttt \phi_n  + f'(w (t) )\ttt \phi_n
=  \ttt g_n (\by ,t)
$$
where   $\ttt g_n (\by ,t): =  g_n(  \ttt \py_n + \A \by , t)$.
We observe that this expression is valid for $\py$ well-inside the domain $\A^{-1}{\mathcal U}_k$ which is expanding to entire $\R^2$.
Since $\ttt\phi_n$ is bounded, and $\ttt g_n \to 0$ in $L^p_{loc}(\R^2)$, we obtain local uniform $W^{2,p}$-bound.  Hence we may assume, passing to a subsequence, that
$\ttt \phi_n $ converges uniformly in compact subsets of $\R^3$ to a function $\ttt \phi (\by, t)$ that satisfies

$$
a_{ij}^0 (\ttt \py )\pp_{ij}\ttt \phi
   +  \partial_{tt}\ttt \phi  + f'(w (t) )\ttt \phi
= 0 \, .
$$
Thus $\ttt \phi$ is non-zero and bounded. After a rotation and stretching of coordinates, the constant coefficient
operator $a_{ij}^0 (\ttt \py )\pp_{ij}$ becomes $\Delta_{\by}$. Hence Lemma \ref{lemma l1}
implies that, necessarily,
$\ttt \phi (\by,t) = Cw'(t)$. On the other hand, we have
$$
0= \int_\R \ttt \phi_n (\by,t) \, w'(t)\, dt \longrightarrow \int_\R \ttt \phi (\by,t) \, w'(t)\, dt\  \quad\hbox{as } {n\to \infty}.
$$
Hence, necessarily $\ttt\phi \equiv 0$. But we have $ (1+ r(\A_n y_n))^\mu\,|\ttt \phi_n ( 0, t_n)| \ge \frac 12 $, and since $t_n$ and
$r(\A_ny_n)$ were bounded,
the local uniform convergence implies $\ttt \phi \ne 0$. We have reached a contradiction.

\medskip
\subsubsection{ Subcase $t_n$ unbounded}
If $y_n$ is in the same range as above, but, say, $t_n\to +\infty$, the situation is similar. The variation is that we define
now
$$
 \ttt \phi_n (\by, t) =    e^{\sigma(t_n + t) } \phi_n ( \A_n^{-1} \py_n + \by, t_n+  t), \quad  \ttt g_n (\by, t) =    e^{\sigma(t_n + t) } g_n ( \A_n^{-1}\py_n +  \by, t_n+  t).
 $$
 Then $ \ttt \phi_n $ is uniformly bounded, and $\ttt g_n \to 0$ in $L^p_{loc} (\R^3)$.
 Now $\ttt \phi_n$ satisfies 
 $$
{a_{ij}^{0}} ( \py_n + \A_n \by)\, \pp_{ij} \ttt \phi_n \, + \,\pp_{tt} \ttt \phi_n \,+ \, \A_n b_j ( \py_n + \A_n \by)\, \pp_{j} \ttt \phi_n \
 $$
 $$
 - 2\sigma \,\pp_t \ttt \phi_n  \, + \,(f'(w(t+t_n ) +\sigma^2 )\, \ttt \phi_n \,=\, \ttt g_n .
 $$
 We fall into the limiting situation
\be
a_{ij}^{*} \, \pp_{ij} \ttt \phi \, + \,\pp_{tt} \ttt \phi  \,
 - \,2\sigma \,\pp_t \ttt \phi  \, - \,( \sigma_+^2 -\sigma^2 )\, \ttt \phi \,=\,  0  \quad \hbox{in }\R^3
\label{positive}\ee
where $a_{ij}^{*}$ is a positive definite, constant  matrix and $\ttt \phi\ne 0$.
But since, by hypothesis   $\sigma_+^2 -\sigma^2 >0$, maximum principle implies that $\ttt \phi \equiv 0$. We obtain a contradiction.

\medskip
\subsubsection{Case  $r(\A_n y_n ) \to +\infty$.}
In this case we may assume that the sequence $\A_n y_n$ diverges along one of the ends, say $M_k$. Considering now the
parametrization associated to the end, $ y = \psi_k(\py )$,  given by \equ{Yk}, which inherits that for $M_{\A_n, k}$, $ y = \A_n^{-1} \psi_k(\A_n \py )$. Thus in this case $ a_{ij}^0 (\ttt \py_n + \A_n  \by ) \to \delta_{ij}$,  uniformly in compact subsets
 of $\R^2$.

 \subsubsection{ Subcase $t_n$  bounded}
Let us assume  first that the sequence $t_n$ is bounded and set
$$
\ttt  \phi_n ( \by, t) =    (1+ r ( \ttt\py_n + \A_n \by ))^\mu  \,  \phi_n (  \A_n^{-1}\ttt \py_n + \by  , t_n + t ) . $$
 Then
 $$
\pp_j ( r_{\A_n}^{-\mu} \ttt \phi_n ) =   -\mu \A \, r^{-\mu-1} \pp_j r\ttt\phi + r^{-\mu} \pp_j \ttt \phi
$$
$$
\pp_{ij} (r_{\A_n}^{-\mu} \ttt \phi_n) =  \mu(\mu+1)\A^2 r^{-\mu-2} \pp_i r\pp_j r \phi
 -\mu\A^2 r^{-\mu-1}\pp_{ij} r  \ttt\phi  - \mu\A r^{-\mu-1} \pp_j r\pp_i\ttt\phi
$$
$$
+ r^{-\mu}\pp_{ij} \ttt \phi - \mu \A  r^{-\mu-1}\pp_i r  \pp_j\ttt \phi  \ .
$$
Now $\pp_i r = O(1)$, $\pp_{ij} r = O( r^{-1})$, hence we have
$$
\pp_j ( r_{\A_n}^{-\mu} \ttt \phi_n ) =    r^{-\mu}\,\left [\,
\pp_j \ttt \phi    + O(\A r_\A^{-1}) \ttt \phi \, \right]\, ,
$$
$$
\pp_{ij} (r_{\A_n}^{-\mu} \ttt \phi_n) = r_{\A} ^{-\mu}\, \left [ \, \pp_{ij} \ttt \phi + O(\A r_\A^{-1}) \pp_i \ttt \phi +
O(\A^{2} r_\A^{-2}) \ttt \phi\, \right ] ,
$$
and the equation satisfied by $\ttt \phi_n$ has therefore the form
 $$
 \Delta_\by \ttt \phi_n \, + \,\pp_{tt} \ttt \phi_n \, + \, o(1) \pp_{ij} \ttt \phi_n + \, o(1)\, \pp_{j} \ttt \phi_n \ + \ o(1)\, \ttt \phi_n  + f'(w(t))\ttt \phi_n= \ttt g_n .
 $$
 where $\ttt \phi_n$ is bounded, $\ttt g_n \to 0$ in $L^p_{loc}(\R^3)$. From elliptic estimates, we also get uniform bounds for
  $\|\pp_j\ttt \phi_n\|_\infty$ and $\|\pp_{ij}\ttt\phi_n\|_{p,0,0}$.
  In the limit we obtain a $\ttt \phi \ne 0$ bounded, solution of
   \be
  \Delta_{\by} \ttt \phi + \pp_{tt} \ttt\phi + f'(w(t))\ttt \phi = 0, \quad  \int_\R \ttt \phi (\by,t) \, w'(t)\, dt\ = \ 0\,  ,
  \label{li1}\ee
a situation which is discarded in the same way as before if $\ttt \phi$ is defined in $\R^3$.
There is however, one more possibility which is that $ r(\A_n y_n) - R_n = O(1)$. In such a case we would see in the limit equation \equ{li1}  satisfied
 in a half-space, which after a rotation in the $\by$-plane can be assumed to be
 $$ H = \{ (\by,t) \in \R^2\times \R \ / \  \py_2 < 0\, \},\quad \hbox{with } \phi(\py_1,0,t) =0 \foral (\py_1,t)\in \R^2. $$
 By Schwarz's reflection, the odd extension of $\ttt \phi$, which achieves for $y_2 >0$,\\ $\ttt\phi (y_1,y_2,t)= -\ttt\phi (y_1,-y_2,t)$,
 satisfies the same equation, and thus we fall into one of the previous cases, again finding a contradiction.

\subsubsection{Subcase  $t_n$ unbounded}
Let us assume now  $|t_n|\to +\infty$. If $t_n\to +\infty$ we define
 $$
\ttt  \phi_n ( \by, t) =    (1+ r ( \ttt\py_n + \A_n \by ))^\mu \, e^{t_n + t} \,  \phi_n (  \A_n^{-1}\ttt \py_n + \by  , t_n + t ) . $$
In this case we end up in the limit with a $\ttt \phi \ne 0$ bounded and satisfying the equation
$$
\Delta_{\by} \ttt \phi \, + \,\pp_{tt} \ttt \phi  \,
 - \,2\sigma \,\pp_t \ttt \phi  \, - \,( \sigma_+^2 -\sigma^2 )\, \ttt \phi \,=\,  0
$$
either in entire space or in a Half-space under zero boundary condition. This implies again $\ttt \phi =0$, and a contradiction has been reached
that finishes the proof of the a priori estimates. \qed

\bigskip
\subsection{ Existence: conclusion of proof of Proposition \ref{prop1}}
Let us prove now existence.
We assume first that $g$ has compact support in $M_\A \times \R$.
\begin{align}
\begin{aligned}
 \partial_{tt}\phi + \Delta_{y,M_\A} \phi   +f'(w(t))\phi &=  g(y,t)  +  c (y) w'(t)   \quad\hbox{in }   M_\A^R\times \R
,\\
 \phi\,&=0,  \quad\hbox{on }  \partial M_\A^R\times \R,\\
\int_{-\infty} ^ {\infty} \phi(y,t)\,w'(t)\,dt&=0  \foral  y\in  M_\A^R,
\end{aligned}
\label{p11}
\end{align}
where we {allow $R= +\infty$} and
$$
c(y) \int_\R {w'}^2 dt =  - \int_\R  g(y,t) \,w'\, dt \ .
$$
Problem \equ{p11} has a weak formulation which is the following.
Let
$$ H = \{ \phi \in H^1_0( M_\A^R \times \R) \ /\  \int_\R \phi(y,t)\,w'(t)\,dt =0  \foral  y\in  M_\A^R\, \} \, . $$
$H$ is a closed subspace of $H_0^1( M_\A^R \times \R) $, hence a Hilbert space when endowed with  its natural norm,
$$
\|\phi\|_H^2 = \int_{M_\A^R} \int_\R (\, |\pp_t\phi|^2 + |\nn_{M_\A} \phi |^2 - f'(w(t)\, \phi^2 \, )\, dV_\A\, dt\ .
$$
$\phi$
is then a weak solution of Problem \equ{p11} if
$\phi\in H$ and satisfies
$$
a(\phi, \psi) := \int_{M_\A^R \times \R} \left (\, \nn_{M_\A} \phi \cdot \nn_{M_\A} \psi\, -\, f'(w(t))\, \phi\, \psi\,\right )\, dV_\A\, dt\, =
 $$
 $$- \int_{M_\A^R \times \R}  g\, \psi\, dV_\A\, dt \foral \psi\in H. $$
It is standard to check that a weak solution of problem \equ{p11}  is also classical provided that $g$ is regular enough.
Let us observe that because of the orthogonality  condition defining $H$ we have that
$$
\gamma \int_{M_\A^R \times \R} \psi^2 \, dV_\A\, dt\ \le \ a(\psi, \psi) \foral \psi \in H.
$$
Hence the bilinear form $a$ is coercive in $H$, and existence of a unique weak solution follows from Riesz's theorem. If $g$ is regular and compactly supported,
$\psi$ is also regular. Local elliptic regularity implies in particular that $\phi$ is bounded.
Since for some $t_0>0$, the equation satisfied by $\phi$ is
\be
\Delta \phi + f'(w(t))\, \phi = c(y) w'(t), \quad |t|> t_0 , \quad y\in M_\A^R,
\label{www}\ee
and $c(y)$ is bounded, then enlarging $t_0$ if necessary, we see that for $\sigma < \min\{ \sigma_+,\sigma_-\}$,
the function $v(y,t):= Ce^{-\sigma|t|} + \ve e^{\sigma|t|}$ is a positive
supersolution of equation \equ{www}, for a large enough choice of $C$ and arbitrary $\ve>0$. Hence $|\phi| \le  Ce^{-\sigma|t|}$, from maximum principle. Since $M_\A^R$ is bounded, we conclude that
$\|\phi\|_{p,\mu,\sigma} < +\infty  $.  From Lemma \ref{apriori} we obtain that
if $R$ is large enough then
\be
 \|D^2 \phi\|_{p,\mu,\sigma} + \|D \phi\|_{\infty ,\mu,\sigma} + \| \phi\|_{\infty ,\mu,\sigma} \le C \|g\|_{p,\mu,\sigma}
\label{coota}\ee
Now let us consider Problem \equ{p11} for $R= + \infty$, allowed above, and for $\|g\|_{p,\mu,\sigma} <+\infty$. Then solving the equation for
finite $R$ and suitable compactly supported $g_R$, we generate a sequence of approximations $\phi_R$ which is uniformly controlled in $R$
by the above estimate. If $g_R$ is chosen so that $g_R\to g$ in $L^p_{loc} ( M_\A \times \R)$ and $\|g_R\|_{p,\mu,\sigma} \le C\|g\|_{p,\mu,\sigma}$, We obtain that $\phi_R$ is locally uniformly bounded, and by extracting a subsequence, it converges uniformly locally over compacts to
a solution $\phi$ to the full problem which respects the estimate \equ{cota}. This concludes the proof of existence, and hence that of the proposition. \qed

\medskip
\setcounter{equation}{0}
\section{ \emph{The Jacobi operator}}\label{jacobi}

We consider this section the
problem of finding a function $h$ such that for certain constants $c_1,\ldots, c_J$,
\begin{equation}\label {proj1}
\JJ(h) = \Delta_M h  + |A|^2h = f  + \sum_{j=1}^J  \frac {c_i}{1+ r^4}\, \hat{z}_i\quad\hbox{ in } M,
\ee
\be\label {proj2}
\int_M  \frac {\hat{z}_ih} {1+ r^4}\, =\, 0, \quad i=1,\ldots , J
\end{equation}
and prove the result of Proposition \ref{proj}. We will also deduce the existence
of Jacobi fields of logarithmic growth as in Lemma \ref{lemin6}.
We recall the definition of the norms $\|\ \|_{p,\beta}$ in \equ{*p}.

\medskip
Outside of a ball of sufficiently large radius $R_0$, it is natural to parameterize each end of $M$, $y_3= F_k(y_1,y_2)$ using the Euclidean coordinates $  \py =(y_1,y_2) \in \R^2$.
The requirement in $f$  on each end amounts to
 $\tilde f \in  L^p( B(0, 1/R_0))$ where
\begin{equation}\label{tf}
\ttt f (\by ) := |\by|^{-4}f ( |\by|^{-2}\by ) \, .
\end{equation}
Indeed, observe that
$$
\|\ttt f \|_{L^p(B(0, 1/R_0))}^p = \int_{B(0, 1/R_0)} |\by|^{-4p} | \, f(  |\by|^{-2} \by ) \,  |^p\, d \by\,=\, \int_{\R^2\setminus B(0,R_0)}
|\py|^{4(p-1) } | f (\py) |^p\,  d\py\ .
$$
In order to prove the proposition we need some a priori estimates.

\begin{lemma} \label{lem8.1}
Let $p>2$. For each  $R_0>0 $ sufficiently large there exists a constant $C>0$ such that if  $$\|f\|_{{p,4-\frac 4p}} + \|h\|_{L^\infty(M)} <+\infty$$
and  $h$ solves
$$
\Delta_M h  + |A|^2h = f, \quad    y\in M,\quad |y |> R_0 \, ,
$$
then 
 $$
 \| h\|_{L^\infty ( |y | > 2R_0 )}  + \|\,| y |^2 Dh\|_{L^\infty (| y | > 2R_0)}    + \| \,|y|^{4- \frac 4p} D^2h\|_{L^p (|y | > 2R_0)}
\ \le
$$
$$
C \, [ \, \|f\|_{{p,4-\frac 4p}} \, + \, \|h\|_{L^\infty( R_0 <|y|< 3R_0)}\  ] \, .
$$
\end{lemma}

\proof
Along each end $M_k$ of $M$,
$\Delta_M$ can be expanded in the coordinate $\tt y$ as
$$
\Delta_M = \Delta +   O( |\py |^{-2}) D^2 + O(|\py |^{-3}) D .
$$
A solution of  $h$  of equation \equ{proj1} satisfies
$$
\Delta_M h  + |A|^2h = f, \quad |\py |> R_0
$$
for a sufficiently large $R_0$.  Let us consider a Kelvin's transform

$$h(\py ) = \ttt h( \py /|\py |^2) .$$ Then
we get
$$
\Delta h (\py ) =  |\py |^{-4} (\Delta \ttt h) (\py /|\py |^2) \, .
$$
Besides
$$
O( |\py |^{-2}) D^2 h(\py )  + O(|\py |^{-3}) Dh(\py ) =  O( |\py |^{-6 }) D^2\ttt h (\py /|\py |^2) +   O(|\py |^{-5})D\ttt h (\py /|\py |^2) \, .
$$
Hence
$$
 (\Delta_M  h) (\by /|\by |^2) =  |\by |^4 \left [ \, \Delta\ttt h (\by )   +  O( |\by |^{2}) D^2\ttt h (\by ) +   O(|\by |)D\ttt h (\by )\right ] \, .
$$
Then $\ttt h$ satisfies the equation
$$
\, \Delta\ttt h   +  O( |\by |^{2}) D^2\ttt h  +   O(|\by |)D\ttt h  + O(1) h = \ttt f (\by )  ,\quad  0< |\by | < \frac 1{R_0}
$$
where $\ttt f$ is given by \equ{tf}.
The operator above satisfies maximum principle in  $B(0, \frac 1{R_0})$ if $R_0$ is fixed large enough.
This, the fact that
 $\ttt h$ is bounded, and $L^p$-elliptic regularity for $p>2$ in two dimensional space imply
that
$$
\|\ttt h\|_{L^\infty ( B(0,1/2R_0))}  + \| D \ttt h\|_{L^\infty ( B(0,1/2R_0))}  + \| D^2 \ttt h\|_{L^p ( B(0,1/2R_0))} \ \le\
$$

$$ C[ \|\ttt f\|_{L^p((B(0, 1/{R_0}))} +  \|\ttt h\|_{L^\infty (
1/3R_0 < |\py | <  1/{R_0})} ] \, \le \,
$$

$$
 C \, [\, \|f\|_{{p,4-\frac 4p}} +
\|h\|_{L^\infty(B( R_0< |\py|< 3R_0)) }\, ]\, .
$$

Let us observe that
$$
\|\ttt h\|_{L^\infty ( B(0,1/2R_0))} = \|h\|_{L^\infty ( |\py|> 2 R_0 )},\quad
$$
$$
\|D\ttt h\|_{L^\infty ( B(0,1/2R_0))} = \|\,|\py |^2\, Dh\|_{L^\infty ( |\py|> 2 R_0 )}.
$$
Since
$$
|D^2h(\py) | \le  C(\,  |\py|^{-4}  \, |D^2\ttt h (\,  |\py|^{-2} \py)| + |\py|^{-3} |D\ttt  h (\, |\py|^{-2} \py)| \, )
$$
then
$$
|\py|^{4 -\frac 4p} |D^2h(\py)| \ \le \     C(\,  |\py|^{-4/p }  |D^2\ttt h (\, |\py|^{-2} \py)| +      |\py|^{- \frac 4p -1 }  |D\ttt  h ( |\py|^{-2} \py)| \, ).
$$
Hence
$$
\int_ {|\py|> 2 R_0 } |\py|^{4p-4} |D^2h|^p  d\py \  \le
$$
$$
C(\, \int_{B(0,1/2R_0)}  |D^2\ttt h(\by) |^p \, d\by \, +  \,    \|D\ttt  h\|^p_{L^\infty ( B(0,1/2R_0))}\int_{|\py|> 2 R_0 } |\py|^{-4- p }d \py\, ).
$$

It follows that
 $$
\| h\|_{L^\infty ( |\py | > 2R_0 )}  + \|\,|\py |^2 Dh\|_{L^\infty (| \py | > 2R_0)} + \|\,|\py |^{4-\frac 4p} D^2h\|_{L^p (| \py | > 2R_0)} \ \le
 $$

 $$C \,[ \, \|f\|_{{p,4-\frac 4p}} + \|h\|_{L^\infty(B( R_0< |\py|< 3R_0))}\  ] \, .
$$
Since this estimate holds at each end, the result of the lemma follows, after possibly changing slightly the value $R_0$. \qed

\bigskip

\begin{lemma}\label{lemin2}
Under the conditions of Lemma \ref{lem8.1},
assume that $h$ is a bounded solution of Problem $\equ{proj1}$-$\equ{proj2}$. Then the a priori estimate \equ{est}
holds.
\end{lemma}

\proof

Let us observe that this a priori estimate in Lemma \ref{lem8.1} implies in particular that
the Jacobi fields $\hat{z}_i$ satisfy$$\nn \hat{z}_i (y)= O(|y|^{-2})\quad\hbox{ as }|y|\to +\infty .$$

\medskip
Using $\hat{z}_i$ as a test function in a ball $B(0,\rho)$ in $M$ we obtain
$$
\int_{\partial B(0,\rho)} ( h {\partial_\nu} \hat{z}_i -  \hat{z}_i {\partial_\nu} \hat{z}_i)\,   + \int_{|y|<\rho} (\Delta_M \hat{z}_i + |A|^2 \hat{z}_i )\, h  =
$$
$$
\int_{|y|<\rho}  f\hat{z}_i \ +    \sum_{j=1}^J  c_j \int_M \frac{\hat{z}_i \hat{z}_j}{1+ r^4}  .
$$
Since the boundary integral in the above
identity is of size $O(\rho^{-1})$  we get
\be
\int_{M}  f\hat{z}_i \ +    \sum_{j=1}^J  c_j \int_M \frac{\hat{z}_i \hat{z}_j}{1+ r^4}= 0
\label{cj}\ee
so that  in particular
\begin{equation}\label{estc}
|c_j| \le C \|f\|_{{p,4-\frac 4p}}\foral j=1,\ldots, J.
\end{equation}

In order to prove the desired estimate, we assume by contradiction that there are sequences $h_n, f_n$ with $\|h_n\|_\infty = 1$ and
$\|f_n\|_{{p,4-\frac 4p}} \to 0 $, such that
$$
\Delta_M h_n  + |A|^2h_n = f_n + \sum_{j=1}^J \frac{c_i^n \hat{z}_i}{1+ r^4}
$$
$$
\int_M \frac{h_n \hat{z}_i}{1+ r^4} = 0 \foral i=1,\ldots, J.
$$
Thus according estimate \equ{estc}, we have that $c_i^n \to 0$. From Lemma \ref{lem8.1} we find
 $$
\| h_n\|_{L^\infty ( |y| > 2R_0 )}   \le C [ o(1) + \|h_n\|_{L^\infty(B(0,3R_0))} ] \, .
$$
The latter inequality implies that
$$
\|h_n\|_{L^\infty(B(0,3R_0))} \ge \gamma >0 .
$$
Local elliptic estimates imply a $C^1$ bound for $h_n$ on bounded sets. This implies the presence of a subsequence $h_n$ which we denote the same way
such that $h_n\to h$ uniformly on compact subsets of $M$, where $h$ satisfies
$$
\Delta_M h  + |A|^2h = 0 \, .
$$
$h$ is bounded hence, by the nondegeneracy assumption, it is a linear combination of the functions $\hat{z}_i$. Besides $h\ne 0$
and satisfies
$$
\int_M \frac{h\hat{z}_i}{1+ r^4}  = 0 \foral i=1,\ldots, J \, .
$$
The latter relations imply $h = 0$, hence a contradiction that proves the validity of the a priori estimate. \qed

\bigskip
\subsection{Proof of Proposition \ref{proj}} \label{proi}
 Thanks to Lemma \ref{lemin2} it only remains to prove existence of a bounded solution
to problem \equ{proj1}-\equ{proj2}.  Let $f$ be as in the statement of the proposition.
Let us consider the Hilbert space $H$ of functions  $h\in H_{loc }^1(M)$ with
$$
\|h\|_H^2:= \int_M |\nn h|^2 + \frac 1 {1+ r^4} |h|^2 < + \infty \, ,
$$
 $$
\int_M \frac{1}{1+r^4} h \hat{z}_i = 0 \foral i=1,\ldots, J \, .
$$
Problem \equ{proj1}-\equ{proj2} can be formulated in weak form as that of finding $h\in H$
with
 $$
\int_M \nn h \nn \psi  -|A|^2h\,\psi =  - \int_{M} f\psi \foral \psi \in H \, .
$$
In fact, a weak solution $h\in H$ of this problem must be bounded thanks to elliptic regularity, with the use of
Kelvin's transform in each end for the control at infinity.
Using that $|A|^2 \le   C r^{-4},$ Riesz representation theorem
and the fact that
$H$ is compactly embedded in
$
L^2( (1+ r^4)^{-1} dV )
$ (which follows for instance by inversion at each end), we see that
this weak problem can be written as an equation of the form
$$
h  - T( h) = \ttt f
$$
where $T$ is a compact operator in $H$ and $\ttt f \in H$ depends linearly on $f$.
When $f=0$, the a priori estimates found yield that necessarily $h=0$. Existence of a  solution
 then follows from
Fredholm's alternative. The proof is complete. \qed

\bigskip

\bigskip
\subsection{ Jacobi fields of logarithmic growth. The proof of Lemma \ref{lemin6}}\label{logjac}
 We will use the theory developed above to construct  Jacobi fields with logarithmic growth as $r\to +\infty$, whose existence we stated and
 use to set up the initial approximation in Lemma \ref{lemin6}.
One of these Jacobi fields  is the generator of dilations
 of the surface, $z_0(y) = y\cdot \nu(y)$.
We will prove next that there are  another $m-2$ linearly independent logarithmically growing  Jacobi fields.

\medskip
Let us consider an $m$-tuple of numbers $\beta_1,\ldots, \beta_m$  with  $\sum_j \beta_j =0$, and any smooth function $p(y)$ in $M$ such that on each  end $M_j$ we have that for sufficiently large $r= r(y)$,
$$
p(y) = (-1)^j\beta_j \log r(y) ,\quad y\in M_j
$$
for certain numbers $\beta_1,\ldots, \beta_m$  that we will choose later.
To prove the result of Lemma \ref{lemin6} we need to find a solution $h_0$ of the equation
$
\JJ(h_0) = 0$ of the form
$h_0 = p +h$ where $h$ is bounded.
This amounts to solving
\be
\JJ(h) = -\JJ(p) \, .
\label{jaco}\ee
Let us consider the cylinder $ C_R = \{ x\in \R^3\ /\ r(x)< R\}$ for a large $R$. Then
$$
\int_{M\cap C_R}  \JJ(p)\, z_{3} dV = \int_{M\cap C_R}  \JJ(z_{3}) z_{3} dV + \int_{\partial C_R \cap M} (z_{3} \partial_n p - p\partial_n z_{3})\, d\sigma(y) \, .
$$
Hence
$$
\int_{M\cap C_R}  \JJ(p)\, z_{3} dV =  \sum_{j=1}^m \int_{\partial C_R \cap M_j} (z_{3} \partial_n p - p\partial_n z_{3})\, d\sigma(y) \, .
$$
Thus using the graph coordinates on each end, we find
$$
 \int_{M\cap C_R}  \JJ(p)\, z_{3} dV\  =  $$
 $$\sum_{j=1}^m (-1)^j\,\left [\,  \frac {\beta_j}R \,\int_{|\py|= R} \nu_3  d\sigma(\py)\, - \,
\beta_j \log R  \int_{|\py|= R} \partial_r \nu_{3}  d\sigma(\py)\ \right] \,  + \, O(R^{-1}).
$$
We have that, on each end $M_j$,
$$
\nu_3(\py ) =  \frac {(-1)^j}{\sqrt{ 1+ |\nn F_k(\py )|^2}} = (-1)^j + O( r^{-2}) , \quad \partial_r\nu_3(\py) = O(r^{-3}).
$$
Hence we get
$$
 \int_{M\cap C_R} \JJ(p)\, z_{3} dV\ = \  2\pi \sum_{j=1}^m   \beta_j  + O( R^{-1}) \, .
 $$
It is easy to see, using the graph coordinates that $\JJ(p) = O(r^{-4})$ and it is hence integrable. We pass to the limit $R\to +\infty$
and get
\be
 \int_{M} \JJ(p)\, z_{3} dV\ = \  2\pi \sum_{j=1}^m  \beta_j = 0\  .
\label{bj}\ee

We make a similar integration for the remaining bounded Jacobi fields.
For $z_{i} = \nu_i(y)$ $i=1,2$ we find
$$
 \int_{M\cap C_R}  \JJ(p)\, z_{2} dV =  \sum_{j=1}^m (-1)^j\,\left [\,  \frac {\beta_j}R \,\int_{|\py|= R} \nu_2  d\sigma(\py)\, - \,
\beta_j \log R  \int_{|\py|= R} \partial_r \nu_2  d\sigma(\py)\ \right] \,  + \, O(R^{-1}).
$$
Now, on $M_j$,
$$
\nu_2(\py ) =  \frac {(-1)^j}{\sqrt{ 1+ |\nn F_k(\py )|^2}}  = (-1)^j a_j\frac {x_i} {r^2}  + O( r^{-3}) , \quad \partial_r\nu_2(\py) = O(r^{-2}).
$$
Hence
$$
 \int_{M} \JJ(p)\, z_{i} dV\ = \ 0\  i=1,2.
$$
Finally, for
$z_{4}(y) = (-y_2, y_1,0)\cdot \nu(y)$ we find on $M_j$, 

$$
(-1)^j z_{4}(\py) = -\by_2\partial_2F_j   + \by_1\partial_1 F_j \ =\ b_{j1} \frac { \py_2}{ r^2}  -  b_{j2} \frac {\py_1}{ r^2} + O(r^{-2}) , \quad \partial_r z_{4} = O(r^{-2})
$$
and hence again
$$
 \int_{M} \JJ(p)\, z_{4} dV\ = \ 0\  .
$$
From the solvability  theory
developed,
we can then find a bounded solution to
the problem
$$
\JJ (h)  = -\JJ(p) + \sum_{j=1}^J q c_j  \hat{z}_{j} \, .
$$
Since $\int_M \JJ(p) z_{i}  dV\ =\ 0 $ and hence $ \int_M \JJ(p) \hat{z}_i dV=0$, relations \equ{cj} imply that $c_i=0$ for all $i$.

We have thus  found a bounded solution to equation \equ{jaco} and the proof is concluded. \qed

\medskip

\begin{remark}\label{jacobilog}{\em
Observe that, in particular, the explicit Jacobi field $z_0(y)= y\cdot \nu(y)$ satisfies that
$$
z(y) = (-1)^j a_j \log r + O(1) \foral y\in M_j
$$
and we have indeed $\sum_j a_j =0$.
Besides this one, we thus have the presence of another $m-2$ linearly independent
Jacobi fields with $|z(y)| \sim \log r $ as $r\to +\infty$, where $m$ is the number of ends.

\medskip
These are in reality {\em all } Jacobi fields with exact logarithmic growth. In fact if $\JJ(z) =0$
and
\be |z(y)| \le C\log r \, ,
\label{log}\ee
 then  the argument in the proof of Lemma \ref{lem8.1} shows that
 the  Kelvin's inversion $\ttt z(\py)$ as in the proof of Lemma  \ref{lemin2} satisfies near the origin
$
\Delta \ttt z =  \ttt f$ where $\ttt f$ belongs to any $L^p$ near the origin,
so it must equal a multiple of $\log |\py|$ plus a regular function.
It follows that on $M_j$ there is a number $\beta_j$ with
$$
z(\py ) =  (-1)^j\beta_j \log|\py| + h
$$
where $h$ is smooth and bounded.  The computations above force $\sum_j \beta_j =0$. It follows from Lemma  \ref{lemin6}
that then $z$ must be equal to one of the elements there predicted plus a bounded Jacobi field.
We conclude in particular that the dimension of
the space of Jacobi fields satisfying \equ{log} must be at most $m-1 + J$, thus recovering a fact stated in Lemma 5.2 of  \cite{perez-ros}.}
\end{remark}

\medskip
\setcounter{equation}{0}
\section{  \emph{Reducing the gluing system and solving the projected problem}}\label{nonli}
In this section we prove Lemma \ref{lemapsi1}, which reduces the gluing system \equ{e2}-\equ{e3} to solving the nonlocal equation \equ{e4}
and prove Proposition \ref{prop2} on solving the nonlinear projected problem \equ{nonlinear0}, in which the basic element is linear theory stated in Proposition \ref{prop1}.
In what follows we refer to  notation and objects introduced in  \S \ref{step0}, \S \ref{step1}.

\subsection{Reducing the gluing system}
Let us consider  equation \equ{e2} in the gluing system \equ{e2}-\equ{e3},
\be
 \Delta\psi -W_\A(x)\psi
+
(1-\zeta_2)S(\ww) + (1-\zeta_1)N(\psi + \zeta_2\phi)  + 2\nn\zeta_1\nn\phi + \phi \Delta \zeta_1 \,=\, 0 \quad\hbox{ in } \R^3\,
 \label{e23}\ee
where
$$W_\A(x) := [\, (1-\zeta_1) f'(u_1)  + \zeta_1 H(t)\, ] \ .$$

\subsubsection{Solving the linear outer problem}
 We consider first the linear problem
\be
\Delta \psi - W_\A(x)\psi + g(x)\ =\ 0 \quad \hbox{in } \,  \R^3
\label{pv}\ee
We observe that globally we have  $ 0< a < W_\A(x) < b$ for certain constants $a$ and $b$. In fact we can take
$ a= \min\{ \sigma_-^2,\sigma_+^2\}-\tau$ for arbitrarily small $\tau >0$.

\medskip
We consider for the purpose the norms for $1<p\le +\infty$,
$$
\| g\|_{p,\mu} := \sup_{x\in \R^3} (1+ r(\A x))^\mu \|g \|_{L^p(B(x,1))} ,\quad r(x',x_3) = |x'|\ .
$$

\begin{lemma}\label{lemapsi}
Given $p>3$, $\mu \ge 0$, there  is a $C>0$  such that for all sufficiently small $\A$ and any $g$ with $\| g\|_{p,\mu}< +\infty$  there exists a unique $\psi$ solution to Problem $\equ{pv}$ with $\| \psi \|_{\infty ,\mu} < +\infty$. This solution
satisfies in addition,
\be
\| D^2\psi \|_{p,\mu}  +  \| \psi \|_{\infty ,\mu}  \, \le \, C\|g \|_{p,\mu}.
\label{ene}\ee

\end{lemma}

\proof

We claim that the a priori estimate
\be\label{apr} \| \psi \|_{\infty ,\mu}  \, \le \, C\|g \|_{p,\mu}\ee
holds for solutions $\psi$ with $\| \psi \|_{\infty ,\mu} < +\infty$ to problem \equ{pv} with $\| g\|_{p,\mu}< +\infty$ provided that $\A$ is small enough. This and local elliptic estimates in turn implies the validity of \equ{ene}.
To see this, let us assume the opposite, namely the existence $\A_n\to 0$, and  solutions $\psi_n$ to equation \equ{pv} with
$\| \psi_n \|_{\infty ,\mu} =1$,  $\| g_n\|_{p,\mu}\to 0$.
Let us consider a point $x_n$ with
$$
(1+ r(\A_n x_n))^\mu \psi_n (x_n) \ge \frac 12
$$
and define
$$
\ttt\psi_n (x) = (1+ r(\A_n (x_n+ x))^\mu \psi_n (x_n+ x), \quad \ttt g_n (x) = (1+ r(\A_n (x_n+ x))^\mu g_n (x_n+ x),
$$
$$
\ttt W_n (x) =  W_{\A_n} (x_n+ x).
$$
Then, similarly to what was done in the previous section, we  check that the equation satisfied by $\ttt \psi_n$
has the form
$$
\Delta \ttt\psi_n  - \ttt W_n(x) \ttt\psi_n + o(1)\nn  \ttt\psi_n + o(1)\ttt\psi_n = \ttt g_n .
$$
$\ttt \psi_n$ is uniformly bounded. Then elliptic estimates imply $L^\infty$-bounds for the gradient and the existence of a subsequence
uniformly convergent over compact subsets of $\R^3$ to a bounded solution $\ttt \psi\ne 0$ to an  equation of the form
$$
\Delta \ttt\psi  -  W_*(x) \ttt\psi =0 \quad \hbox{in }\R^3
$$
where $0< a\le W_*(x)\le b$.  But maximum principle makes this situation impossible, hence estimate \equ{apr} holds.

\medskip
Now, for existence, let us consider $g$ with $\| g\|_{p,\mu}< +\infty$ and a collection of approximations $g_n$ to $g$ with $\| g_n\|_{\infty,\mu}< +\infty$,
$g_n\to g $ in $L^p_{loc} (\R^3)$ and $\|g_n\|_{p,\mu}\le C\|g\|_{p,\mu}$.  The problem
$$
\Delta \psi_n - W_n(x)\psi_n = g_n \quad\hbox{in }\R^3
$$
can be solved since
this equation has a positive supersolution of the form\\
$
C_n (1 + r(\A x)\, )^{-\mu},
$
provided that $\A$ is sufficiently small, but independently of $n$. Let us call $\psi_n$ the solution thus found, which satisfies
$\|\psi_n\|_{\infty, \mu}< +\infty$. The a priori estimate shows that
$$
\| D^2\psi_n \|_{p,\mu}  +  \| \psi_n \|_{\infty ,\mu}  \, \le \, C\|g \|_{p,\mu}.
$$
and passing to the local uniform limit up to a subsequence, we get a solution $\psi$
to problem \equ{pv}, with $\|\psi\|_{\infty, \mu}< +\infty$. The proof is complete. \qed

\bigskip
\subsubsection{The proof of Lemma \ref{lemapsi1}}\label{lpsi1}
%


Let us call $\psi:=  \Upsilon (g)$ the solution of Problem \equ{pv} predicted by Lemma \ref{lemapsi}. Let us write
Problem \equ{e23} as fixed point problem in the space $X$ of $W^{2,p}_{loc}$-functions $\psi$ with
$\|\psi\|_X  < +\infty$,
\be
\psi = \Upsilon( g_1 + K(\psi)\,)
\label{fp}\ee
where
$$
g_1=  (1-\zeta_2)S(\ww) + \, 2\nn\zeta_1\nn\phi + \phi \Delta \zeta_1\,,
 \quad K(\psi) = (1-\zeta_1)N(\psi + \zeta_2\phi)\ .
$$
Let us consider a function $\phi$ defined in $M_\A\times \R$ such that
$
\|\phi\|_{2,p,\mu, \sigma} \le 1.
$
Then,
$$
 |\, 2\nn\zeta_1\nn\phi + \phi \Delta \zeta_1\, | \, \le \, C e^{- \sigma \frac {\delta}\A} \, (1+ r(\A x))^{-\mu} \|\phi\|_{2,p,\mu, \sigma}.
$$
We also have that
$ \|S(\ww)\|_{p,\mu, \sigma} \le C\A^3,$ hence
$$
|(1-\zeta_2)S(\ww)| \le C e^{- \sigma \frac {\delta}\A} \, (1+ r(\A x))^{-\mu}
$$
and
$$\|g_1\|_{p,\mu} \le  Ce^{- \sigma \frac {\delta}\A}. $$
Let consider the set $$\Lambda \ =\ \{ \psi\in X \ /\ \|\psi\|_{X} \le   A e^{- \sigma \frac {\delta}{\A}}\} ,$$
for a large number $A>0$.
Since
$$
|\, K(\psi_1) - K(\psi_2)\, |\, \le\, C (1-\zeta_1)\, \sup_{t\in (0,1)} \, |t\psi_1 + (1-t) \psi_2 + \zeta_2\phi| \, |\psi_1 -\psi_2|\, ,
$$
we find that
$$
\|\, K(\psi_1) - K(\psi_2)\, \|_{\infty, \mu}\, \le\,  C\, e^{-\sigma \frac {\delta}{\A}}\|\,\psi_1 - \psi_2\, \|_{\infty, \mu}
$$
while
$\|K(0)\|_{\infty, \mu }\, \le \, C\,e^{- \sigma \frac {\delta}{\A}}.$ It follows that the right hand side of equation \equ{fp} defines a contraction mapping of $\Lambda$, and hence a unique solution $\psi = \Psi( \phi) \in \Lambda$ exists, provided that the number $A$ in the definition of $\Lambda$ is taken sufficiently large and $\|\phi\|_{2,p,\mu,\sigma} \le 1$. In addition,
it is direct to check the Lipschitz dependence of $\Psi$  \equ{lipspsi} on $\|\phi\|_{2,p,\mu,\sigma} \le 1$. \qed

\bigskip
Thus, we  replace  replace $\psi= \Psi(\phi)$ into the equation \equ{e3}  of the gluing system \equ{e2}-\equ{e3}
 and  get the (nonlocal) problem,
\be
\partial_{tt}\phi\, + \, \Delta_{y,M_\A}\phi\, = \, -\ttt S(u_1) - \NNN(\phi)
\quad \hbox{ in } M_\A\times \R
\label{e41}\ee
where
\be
\NNN(\phi) :=  \underbrace{\BB(\phi) + [f'(u_1)- f'(w)]\phi }_{\NNN_1(\phi)} \, + \, \underbrace{\zeta_1 (f'(u_1) - H(t)) \Psi(\phi) }_{\NNN_2(\phi)}\,  +
 \underbrace{\zeta_1 N( \Psi(\phi) + \phi) }_{\NNN_3(\phi)},
\label{NNN1}\ee
which is what we concentrate in solving next.

\bigskip
\subsection{Proof of Proposition \ref{prop2}}\label{pr2}

\medskip
We recall from \S \ref{step1} that Proposition \ref{prop2} refers to
solving the projected problem 

\begin{align}
\begin{aligned}
\partial_{tt}\phi\, + \, \Delta_{y,M_\A}\phi\,& = \, -\ttt S(u_1) - \NNN(\phi)\,  +\, c (y) w'(t)
\quad \hbox{ in } M_\A\times \R,
\\
\int_\R \phi(y,t)\,w'(t)\,dt&=0,  \foral  y\in  M_\alpha,
\end{aligned}
\label{nonlinear1}
\end{align}
and then adjust $h_1$ so that $c(y)\equiv 0$.
Let
$\phi = T(g)$ be the linear operator providing the solution in Proposition \ref{prop1}.
Then  Problem \equ{nonlinear1} can be reformulated as the fixed point problem

\be
\phi \ = \ T(  -\ttt S(u_1) - \NNN(\phi)\, ) =:  \TT (\phi), \quad \| \phi\|_{2,p,\mu,\sigma} \le 1
\label{fp11}\ee
which is equivalent to
\be
\phi \ = \ T(  -\ttt S(u_1) + \A^2\Delta h_1 \, w'  - \NNN(\phi)\, ), \quad \| \phi\|_{2,p,\mu,\sigma} \le 1,
\label{fp1}\ee
 since the term added has the form $\rho(y)w'$ which thus adds up to $c(y)w'$.  The reason to absorb this term is that  because of
 assumption \equ{ass1},
$\|\A^2\Delta h_1 \, w'\|_{p,4,\sigma} = O(\A^{3-\frac 2p})$ while the remainder has a priori size slightly smaller, $O(\A^{3})$.


\subsubsection{Lipschitz character of $\NNN$}
We will solve Problem \equ{fp1} using contraction mapping principle, so that we need to give account of a suitable Lipschitz property for the
operator $\TT$. 
We claim the following.

\medskip{\bf Claim.}  \label{lipschitz}{\em
We have that for a certain constant $C>0$ possibly depending  on $\KK$ in $\equ{ass1}$ but independent of $\A>0$, such that for any
$\phi_1$, $\phi_2$  with
$$\|\phi_l\|_{2,p,\mu,\sigma} \le  K\A^3, $$

 \be
 \|\NNN(\phi_1) -\NNN(\phi_2)\|_{p,\mu+1 ,\sigma} \,\le\,  C \, \A  \, \|\phi_1 -\phi_2\|_{2,p,\mu,\sigma}
  \label{N}\ee

where the operator $\NNN$ is defined in $\equ{NNN1}$.
}

\medskip
We study the Lipschitz character of the operator $\NNN$  through analyzing each of its components.
 Let us start with $N_1$. This is a second order linear operator with  coefficients of order $\A$ plus a decay of order at least $O(r_\A^{-1})$.
We recall that $\BB = \zeta_2 B$ where in coordinates
$$
B\,= \,  (f'(u_1)-f'(w)) - \A^2[ (t+h_1)|A|^2 + \Delta_M h_1]  \partial_t\ -  2\A\, a_{ij}^0 \partial_jh \partial_{it}
\ +
$$

$$
 \A (t+h) \, [ a_{ij}^1\partial_{ij}   - \A\, a_{ij}^1(\, \partial_jh \partial_{it} + \partial_ih \partial_{jt}) + \A (b_i^1\partial_i   - \A b_i^1 \partial_ih \partial_t)\, ) \, ] \ +
 $$

\be
  \A^3 (t+h)^2b^1_3  \partial_t\, +\, \A^2 [\,a_{ij}^0 + \A(t+h) a_{ij}^1 )\,]\partial_ih\partial_jh \, \partial_{tt}
\label{B1}\ee
where, we recall,
$$
a_{ij}^1 = O(r_\A^{-2}),\quad a_{ij}^1 = O(r_\A^{-2}),\quad  b_i^1 = O(r_\A^{-3}), \quad  b_i^3 = O(r_\A^{-6}),
$$
$$
f'(u_1)-f'(w) = O(\A^2 r_\A^{-2} e^{-\sigma|t|}) \quad \partial_jh = O(r_\A^{-1}), \quad |A|^2 =O(r_\A^{-4}) \, .
$$

We claim that
 \be
 \| \NNN_1(\phi)\|_{p,\mu+1,\sigma} \, \le\, C\,\A\, \| \phi\|_{2,p,\mu,\sigma} .
 \label{N1}\ee
 The only term of $N_1(\phi)$ that requires a bit more  attention
 is
 $
 \A^2 (\Delta h_1)(\A y)  \pp_t \phi\, .
 $
 We have
 $$
 \int_{B ( (y,t), 1) }  | \A^2 (\Delta h_1)(\A z)  \pp_t \phi |^p\, dV_\A(z)\, d\tau \ \le \
 $$
 $$
  C\,\A^{2p} \| \pp_t \phi \|_{L^\infty ( B ( (y,t), 1) } \, (1+ r_\A(y)\, )^{- 4p+4 }
 \int_{B ( (y,t), 1) }  | (1+ r_\A(z)\,)^{4- \frac 4p }(\Delta h_1)(\A z)|^p\,| dV_\A(z)\ \le
 $$
 $$
  C\,\A^{2p-2} \|\Delta h_1\|^p_{L^p(M)} e^{-p\sigma|t|} (1+ r_\A(y))^{-p \mu -4p +4 } \|\nn \phi\|_{\infty,\mu,\sigma},
$$
and hence in particular for $p\ge 3$,
$$
\|  \A^2 (\Delta h_1)(\A y)  \pp_t \phi \|_{p,\mu+ 2,\sigma} \le \, C\, \A^{2-\frac 2p} \| h_1\|_{*}\, \| \phi \|_{2,p,\mu,\sigma}
 \le C\,\A^{3- \frac 2p} \| \phi\|_{2,p,\mu,\sigma}.
$$
 Let us consider now  functions $\phi_l$ with
 $$
 \| \phi_l \|_{2,p,\mu,\sigma} \le 1, \quad l=1,2.
 $$
 Now, according to Lemma \ref{lemapsi1}, we get that
 \be
\|\NNN_2(\phi_1) -\NNN_2(\phi_2)\|_{p,\mu,\sigma} \, \le \, C \, e^{-\sigma \frac \delta{\A}} \| \phi_1 - \phi_2\|_{p,\mu,\sigma}\, .
  \label{N2}\ee
Finally, we also have that
$$
|\NNN_3(\phi_1) -\NNN_3(\phi_2)\,| \ \le
$$
$$
C \zeta_1  \sup_{t\in (0,1)} \, |t(\Psi(\phi_1) + \phi_1) + (1-t) (\Psi(\phi_2) + \phi_2)| \,
[\, |\phi_1 -\phi_2| + |\Psi (\phi_1) - \Psi(\phi_2) | ]\, ,
 $$
 hence
 \be
 \|\NNN_3(\phi_1) -\NNN_3(\phi_2)\|_{p,2\mu,\sigma} \,\le\,  C \, (\,\|\phi_1\|_{\infty,\mu,\sigma} + \|\phi_2\|_{\infty,\mu,\sigma}\, + e^{-\sigma \frac \delta\A } ) \, \|\phi_1 -\phi_2\|_{\infty,\mu,\sigma}.
  \label{N3}\ee
From \equ{N1}, \equ{N2} and \equ{N3},  inequality \equ{N} follows. The proof of the claim is concluded.

\subsubsection{Conclusion of the proof of Proposition \ref{prop2}}
 The first observation is that  choosing $\mu \le 3$, we get
\be \| \ttt S(u_1) + \A^2\Delta h_1 w' \|_{p, \mu ,\sigma}\ \le \ C\A^3 . \label{N0}\ee
Let us assume now that $\phi_1, \phi_2\in B_\A$
 where
 $$
 B_\A = \{ \phi \ /\ \|\phi\|_{2,p,\mu,\sigma}\, \le K \A^3 \}
 $$
 where $K$ is a constant to be chosen.
Then we observe that for small $\A$
$$ \| \NNN (\phi) \|_{p,\mu+1,\sigma} \le C \A^4 , \foral \phi\in B_\A, $$
where $C$ is independent of $K$. Then, from relations \equ{N0}-\equ{N3} we see that if $K$ is fixed large enough independent of $\A$, then
the right hand side of equation \equ{fp} defines an operator that applies $B_\A$ into itself, which is also a contraction mapping of $B_\A$ endowed with the
norm $\|\ \|_{p,\mu\sigma}$, provided that $\mu\le 3.$ We conclude, from contraction mapping principle, the existence
of $\phi$ as required.

\medskip
The Lipschitz dependence \equ{rw} is a consequence of series of lengthy but straightforward considerations of the
Lipschitz character in $h_1$ of the operator in the right hand side of equation \equ{fp} for the norm $\|\ \|_*$ defined in \equ{est}. 
Let us recall expression \equ{B1} for the operator $B$,
and consider as an example,  two terms that depend
linearly on $h_1$:
$$
A(h_1,\phi):=   \A\,  a_{ij}^0\, \partial_jh_1 \partial_{it}\phi\, .
$$
Then
$$
|A(h_1,\phi)|\, \le\, C\A |\partial_jh_1|\, | \partial_{it}\phi\, .
$$
 Hence
$$
\|A(h_1,\phi)\|_{p,\mu +2 ,\sigma} \, \le\, C\A \| (1+r_\A^2)\,\partial_jh_1 \|_\infty\, \| \partial_{it}\phi\,\|_{p,\mu,\sigma}\, \le\,
C\A^4 \|h_1\|_*\, \|\phi\|_{2,p,\mu,\sigma}  .
$$
Similarly, for $A(\phi,h_1) = \A^2\Delta_M h_1\,  \partial_t\phi$ we have
$$
|\, A(\phi,h_1)\,|\  \le\ C\A^2 |\Delta_M  h_1 (\A y)|\,(1+r_\A)^{-\mu} e^{-\sigma|t|} \|\phi\|_{2,p,\mu,\sigma} \, .
$$
Hence
$$
\|\A^2\Delta_M h_1\,  \partial_t\phi\,\|_{p,\mu+ 2,\sigma} \ \le\ C\A^{5-\frac 2p} \|h_1\|_* \, \|\phi\|_{2,p,\mu,\sigma}.
$$
We should take into account that some terms involve nonlinear, however mild dependence, in $h_1$. We recall for instance that
$a_{ij}^1 =  a_{ij}^1(\A y, \A(t+ h_0 +h_1))$. Examining the rest of the terms involved we find that  the whole operator $\NNN$ produces a dependence on $h_1$ which is Lipschitz with small constant, and gaining decay in $r_\A$,
\be
\|\NNN(h_1,\phi) -  \NNN(h_2,\phi)\|_{p,\mu+1,\sigma} \le C\A^2 \|h_1 - h_2\|_*\,\|\phi\|_{2, p,\mu,\sigma} .
\label{ln}\ee
Now, in the error term $$ {\mathcal R} = -\ttt S(u_1) + \A^2\Delta h_1 w', $$ we have that
\be\| {\mathcal R} (h_1) - {\mathcal R} (h_2)\|_{p,3,\sigma}  \le   C\,\A^2\, \|  h_1- h_2 \|_*\, .
\label{ed}\ee
To see this, again we go term by term in expansion \equ{error2}. For instance the linear term
$\A^2 \, \,a_{ij}^0\partial_ih_0 \partial_j h_1\, w''. $
We have
$$
 |\A^2 \, \,a_{ij}^0
\, \partial_ih_0 \partial_j h_1 |\, \le\, C\,\A^2\, (1+  r_\A)^{-3}\, e^{-\sigma|t|}\, \|  h_1 \|_*
$$
so that
$$
\|\A^2 \, \,a_{ij}^0
\,\partial_ih_0 \, \partial_j h_1\|_{p,3,\sigma}  \, \le\, C\,\A^2\, \|  h_1\|_* ,
$$
the remaining terms are checked similarly.

Combining estimates \equ{ln}, \equ{ed} and the fixed point characterization \equ{fp} we obtain the desired Lipschitz dependence \equ{rw} of $\Phi$.
This concludes the proof. \qed

\setcounter{equation}{0}
\section{ \emph{The reduced problem: proof of Proposition \ref{proii}}}\label{nonljac}
In this section we prove Proposition \ref{proii} based on the linear theory provided by Proposition \ref{proj} 
Thus, we want to solve the problem
\be
{\mathcal J}(h_1) = \Delta_M h_1 + h_1|A|^2  = G(h_1) +  \sum_{i=1}^J  \frac {c_i}{1+ r^4}\, \hat{z}_i \quad\hbox{ in } M \, ,
\label{probjac2-2}\ee
$$
\int_M \frac {h_1 \hat{z}_i}{1+ r^4}\,dV\, =\, 0 \foral i=1,\cdots, J \, ,
$$
where the linearly independent Jacobi fields $\hat{z}_i$ will be chosen in (\ref{zdef1}) and (\ref{normalized}) of \S 8,
and $G =G_1+ G_2$ was defined in \equ{G1}, \equ{G}.  We will use contraction mapping principle to determine the existence of
a unique  solution $h_1$ for which constraint \equ{ass1}, namely
\be
\|h_1\|_*:= \| h_1\|_{L^\infty (M)} + \|(1+ r^2)D h_1\|_{L^\infty (M)} +  \| D^2 h_1\|_{p, 4-\frac 4p}
 \ \le \ \KK \A\  ,
\label{ass2}\ee
is satisfied after fixing $\KK$ sufficiently large.

\medskip
We need to analyze the size of the operator $G$, for which the crucial step is the following estimate.

\begin{lemma} \label{lemin4}
Let $\psi(y,t)$ be a function defined in $M_\A \times \R$  such that
$$
\|\psi\|_{p,\mu,\sigma}\, := \,  \sup_{ (y,t)\in M_\A \times \R }  e^{\sigma|t|}(1+  r_\A^\mu\, )\,   \| \psi \|_{L^p( B((y,t) , 1)} < +\infty
$$
for $\sigma, \mu \ge 0$.  
The function defined in $M$ as
$$ q(y) :=      \int_\R  \psi(y/\A ,t)\, w'(t)\, dt $$
satisfies
\begin{equation}
\label{qest}
\| q\|_{p, a} \,  \le \,   C\, \|\psi\|_{p,\mu,\sigma}
\end{equation}
provided that $$
\mu  > \frac 2p + a \, .
$$
In particular, for any $\tau>0$,
\be
\| q\|_{p, 2-\frac 2p - \tau } \,  \le \,   C\,   \|\psi\|_{p,2,\sigma}
\label{w1}\ee
and
\be
\| q\|_{p, 4- \frac 4p} \,  \le \,   C\,  \|\psi\|_{p,  4,\sigma} \, .
\label{w2}\ee

\end{lemma}

\proof
We have that for $|y|> R_0$
$$
\int_{|y|> R_0}   |y|^{ap} \left |\,  \int_\R   \psi(y /\A ,t)\, w'(t)\, dt\,     \right |^p  \, dV \ \le\   C
\, \int_\R  \, w'(t)\, dt\,  \int_{|y|> R_0} |y|^{ap} \, |\psi(y /\A ,t)|^p    \, dV \, .
$$
Now
$$
\int_{|y|> R_0} |y|^{ap } \, |\psi(y /\A ,t)|^p    \, dV\ =\  \A^{ap+2}
\int_{|y|> R_0/\A } |y|^{ap } \, |\psi(y ,t)|^p    \, dV_\A \
$$
and
$$
\int_{|y|> R_0/\A } |y|^{ap } \, |\psi(y ,t)|^p    \, dV_\A \ \le \ C \sum_{i \ge [ R_0/\A] }  i^{ap}
\int_{  i < |y| < i+1 }  \, |\psi(y ,t)|^p    \, dV_\A \ .
$$
Now, $i < |y| < i+1 $ is contained in $O(i)$ balls with radius one centered at points of the annulus, hence
$$
{\int_{  i < |y| < i+1 }  \, |\psi(y ,t)|^p    \, d V_\A\  \le \  C e^{-\sigma p |t|} i^{1 -\mu p } \, \|\psi\|_{p,\mu}^p}
$$
$$
{ \le \  C e^{-\sigma p |t|}  \|\psi\|_{p,\mu}^p \int_{ i  < |y| < i+1}  (1+r_\alpha)^{-\mu p}  d V_\alpha}
$$
$$
{ \le \  C e^{-\sigma p |t|}  \|\psi\|_{p,\mu}^p \int_{ i  < |y| < i+1}  | \alpha y|^{-\mu p}  d V_\alpha }
$$
$$
{ \le \  C e^{-\sigma p |t|}  \|\psi\|_{p,\mu}^p  \alpha^{-\mu p} i^{1-\mu p}  } \, .
$$

Then we find
$$
\| \, |y|^ a\, q\|^p_{L^p (|y|> R_0) }  \ \le \ C\,  \A^{ap-\mu p+2}
\|\psi\|_{p,\mu}^p \sum_{i \ge [ R_0/\A] }  i^{ap -\mu p +1 } \, .
$$

The sum converges if
$
\mu  > \frac 2p + a
$
and in this case
$$
\| \, |y|^ a\,q\|^p_{L^p (|y|> R_0) }  \ \le \  C\,  \A^{ap-\mu p+2 }   \A^{-ap +\mu p -2}   \|\psi\|^p_{p,\mu}\ =\  C\,      \|\psi\|^p_{p,\mu}
$$
so that
$$
\| \, |y|^ a\,q\|_{L^p (|y|> R_0) }  \ \le \  C\,      \|\psi\|_{p,\mu} .
$$
Now, for the inner part $|y |< R_0$ in $M$, the weights play no role. We have
$$
\int_{| y| < R_0}  \, |\psi(\py /\A ,t)|^p    \, d V\ =\  \A^{2}
\int_{|y|< R_0/\A } \, |\psi( y ,t)|^p    \, dV_\A  \le
$$
$$
{ \ C\A^{2}\,  \sum_{i \le  R_0/\A  }
\int_{  i < |y| < i+1 }  \, |\psi(y ,t)|^p    \, d V_\A \
 \le C\A^{2}\,\|\psi\|_{p,\mu}^p  e^{-\sigma p |t|}\,
 \sum_{i \le  R_0/\A  }  i \,}
$$
$$
{ \le   C \|\psi\|_{p,\mu}^p  e^{-\sigma p |t|}\, } \, .
$$

Hence if $\mu >\frac{2}{p}+a$ we  finally get
$$
\| q\|_{p,a} \, \le \,  C\,   \|\psi\|_{p,\mu} \
$$
and the proof of (\ref{qest}) is concluded. Letting $(\mu, a)= (2, 2-\frac{2}{p}-\tau), (\mu, a)= (4, 4-\frac{4}{p})$ respectively in (\ref{qest}), we obtain (\ref{w1}) and (\ref{w2}).  \qed

\bigskip

Let us apply this result to $\psi (y,t) = \NNN (\Phi (h_1)\,) $ to estimate the size of the operator $G_2$ in \equ{G}.
For $\phi =\Phi(h_1)$ we have that
$$
G_2(h_1)(y)  :=  c_*^{-1}\A^{-2} \int_\R \NNN (\phi)( y/\A ,t) \, w'\, dt
$$
satisfies
$$
\| G_2(h_1)\|_{p, 4- \frac 4p} \le  C\A^{-2}  \| \NNN (\phi)\|_{p, 4,\sigma} \ \le\ C\, \A^{2}.
$$
On the other hand, we have that, similarly, for $\phi_l = \Phi(h_l)$, $l=1,2$,
$$
\| G_2(h_1) -G_2(h_2) \|_{p, 4- \frac 4p} \,\le \, C\A^{-2}  \| \NNN (\phi_1, h_1) - \NNN(\phi_2,h_2) \|_{p, 4,\sigma} .
$$
Now,
$$
\| \NNN (\phi_1, h_1) - \NNN(\phi_1,h_2) \|_{p, 4,\sigma} \le C\A^2 \|h_1 -h_2\|_*  \|\phi_1\|_{2,p,3,\sigma}, \le C\A^5 \|h_1 -h_2\|_*,
$$
according to inequality \equ{ln}, and
$$
\| \NNN (\phi_1, h_1) - \NNN(\phi_2,h_1) \|_{p, 4,\sigma} \le C\A^2 \|\phi_1-\phi_2\|_{p, 3,\sigma} \le C\A^4 \|h_1 -h_2\|_*\ .
$$
We conclude then that
$$
\| G_2(h_1) -G_2(h_2) \|_{p, 4- \frac 4p} \,\le \, C\, \A^{ 2} \|h_1 -h_2\|_*\ .
$$
In addition, we also have that
$$
\|G_2(0)\|_{p, 4- \frac 4p} \le C\A^2.
$$
for some $C>0$ possibly dependent of $\KK$.
On the other hand, it is similarly checked that the remaining small operator $G_1(h_1)$ in \equ{G1}  satisfies
$$
\| G_1(h_1) -G_1(h_2) \|_{p, 4- \frac 4p} \,\le \, C_1\, \A  \|h_1 -h_2\|_*\  .
$$
A simple but crucial observation we make is that
   $$
   c_* G_1( 0) \, =\,
\, \A  \, \partial_i h_0 \partial_j h_0 \,  \int_\R \zeta_4 (t+ h_0)a_{ij}^1w''w'\, dt\,  +
\A^{-2} \int_\R \zeta_4\, R_1(y,t, 0, 0\, )\, w'\, dt
$$
so that for a constant $C_2$ independent of $\KK$ in \equ{ass2} we have
$$
\quad \|G_1(0)\|_{p, 4- \frac 4p} \le C_2\A \, .
$$
In all we have that the operator $G(h_1)$ has an $O(\A)$  Lipschitz constant, and in addition satisfies
$$
\| G (0) \|_{p, 4- \frac 4p} \le 2C_2\A .
$$
Let $h=T(g)$ be the linear operator defined by  Proposition \ref{proj}. Then we consider the problem  \equ{probjac2-2} written as
the fixed point problem
\be
h_1\, = \, T(\, G(h_1)\, ), \quad \|h\|_*\le   \KK\A .
\label{fp3}\ee
We have
$$
 \| T(\, G(h_1)\, ) \|_* \le \|T\| \,\| G (0) \|_{p, 4- \frac 4p}  + C\A \|h_1\|_*\, .
$$
Hence fixing $\KK > 2C_2\|T\|$, we find that for all $\A$ sufficiently  small, the operator $T\, G$ is a contraction mapping of the ball  $\|h\|_*\le   \KK\A$ into itself. We thus have the existence of a unique solution of the fixed problem \equ{fp3},
namely a unique solution $h_1$ to problem \equ{probjac2-2} satisfying \equ{ass2} and the proof of Proposition \ref{proii}
 is concluded. \qed

\setcounter{equation}{0}
\section{ \emph{Conclusion of the proof of Theorem \ref{teo1}}} \label{conclusion}

We denote  in what follows
$$ r (x)= \sqrt{x_1^2 + x_2^2}, \quad \hat r = \frac 1r (x_1,x_2,0),\quad \hat \theta = \frac 1r (-x_2,x_1,0) \, .$$

We consider the four Jacobi fields associated to rigid motions,  $z_{1},\ldots , z_{4}$ introduced in \equ{jfields}.
 Let $J$ be the number of bounded, linearly independent Jacobi fields of $\JJ$. By our assumption and the asymptotic expansion of the ends (\ref{asyends}), $ 3\leq J\leq 4$. (Note that when $M$ is a catenoid, $ z_{4}=0$ and $ J=3$.) Let us choose
\begin{equation}
\label{zdef1}
\hat{z}_{j}= \sum_{l=1}^4 d_{jl} z_{0l}, j=1, ..., J
\end{equation}
be normalized such that
\begin{equation}
\label{normalized}
\int_M q(y) \hat{z}_i \hat{z}_j=0, \ \mbox{for} \ i \not = j, \int_M q (y) \hat{z}_i^2=1, i, j=1,\cdots, J \, .
\end{equation}
In what follows we fix the function  $q$ as
\be
q(y)\  :=\ \frac 1{1 + r(y)^4}  \ .
\label{py}\ee

\medskip
So far we have built, for certain constants $\ttt c_i$  a solution $u$ of  equation \equ{equ}, namely
$$
\Delta u_{} +  f(u_{}) =   \sum_{j=1}^J \ttt c_i \hat{z}_i(\A y) w'(t) q(\A y) \zeta_2
$$
where  $u_{}$, defined in \equ{uu}  satisfies the following properties
\begin{equation}\label{gato}
u_{} (x) =  w( t) + \phi (y,t)
\end{equation}
near the manifold, meaning this $x = y + (t+ h(\A y)\, )\,\nu(\A y) $ with
$$ y\in M_\A ,\quad |t| \le  \frac  \delta \A  + \gamma \log (2 + r(\A y) ) . $$
The function $\phi$ satisfies in this region the estimate
\begin{equation}\label{grad}
|\phi| + |\nabla \phi|  \le  C \A^2 \frac 1 {1+  r^2 (\A y)}  e^{-\sigma|t|} \, .
\end{equation}

Moreover, we have the  validity of the global estimate
 $$
 |\nn u_{} (x) | \le   \frac C {1+  r^3 (\A x)}  e^{-\sigma\frac \delta\A  } \ .
$$
We introduce the functions
$$
Z_{i}(x) = \partial_{x_i} u_{} (x),\ i=1,2,3, \quad Z_{4}(x) = - \alpha x_2 \partial_{x_2} u_{}  + \alpha x_1 \partial_{x_2} u_{} \, .
$$
From the expansion \equ{gato} we see that

$$
\nn u_{}(x) = w'(t )\,  \nabla t \, + \,  \nabla \phi  .
$$
Now,
$
 t = z- h(\A y)
$
where $z$ designates normal coordinate to $M_\A$. Since $\nabla z = \nu = \nu (\A y) $ we then get
$$
\nn t = \nu(\A y ) -  \A \nn h(\A y).
$$
Let us recall that $h$ satisfies  $h= (-1)^k\beta_k \log r  + O(1) $
along the $k$-th end, and
$$
\nn h =  (-1)^k \frac{ \beta_k} r \hat r + O(r^{-2}) \, .
 $$
From estimate \equ{grad} we
we find that
\begin{equation}
\nn u_{}(x) = w'(t) (\nu -\A (-1)^k \frac{ \beta_k} {r_\A}  \hat r   )  +  O( \alpha   r_\A^{-2} e^{-\sigma |t|} )  .
\label{for}\end{equation}
From here we get that near the manifold, 

\begin{equation}
Z_{i} (x) = w'(t)\, (z_{i}(\A y) -\A (-1)^k \frac{\beta_k}{r_{\A}} \hat{r} e_i)  + O( \alpha r_\A^{-2} e^{-\sigma |t|} ), \quad i=1,2, 3,
\label{estvi}\end{equation}

\begin{equation}
 Z_{4} (x) =   w'(t )\, z_{04} (\A y) +   O( \alpha r_\A^{-1} e^{-\sigma |t|} ).
\label{estv4}\end{equation}
Using the characterization \equ{equ} of the solution $u$ and barriers (in exactly the same way as
in Lemma \ref{lemin100} below which estimates eigenfunctions of the linearized operator), we find the following estimate
for $r_\A (x) > R_0$:
\be
|\nn u_{} (x)| \le  C\sum_{k=1}^m e^{-\sigma| x_3 - \A^{-1} ( F_k(\A x') + \beta_j\A  \log|\A x'| \, )\, | } \ .
\label{end11}\ee


We claim that
\begin{equation}
 \int_{\R^3}  ( \Delta u_{} +  f(u_{})) Z_{i}(x) \, dx = 0 \foral i=1,\ldots, 4
\label{ee1-1}\end{equation}
so that
\begin{equation}
 \sum_{j=1}^J\ttt  c_j  \int_{\R^3 } q(\A x)\hat{z}_j(\A y)  w'(t)   Z_{i} (x) \, \zeta_2 \, dx \, =\, 0  \foral i =1,\ldots , 4.
\label{sis}\end{equation}
Let us accept this fact for the moment. Let us observe that from estimates \equ{estvi} and \equ{estv4},
$$
\A^2 \int_{\R^3 }q(\A x) \hat{z}_j(\A y) w'(t)  \sum_{l=1}^4 d_{il} Z_{l}(x) \, \zeta_2 \, dx   =  \int_{-\infty}^\infty  w'(t)^2 dt \int_M  q\, \hat{z}_j \hat{z}_i dV + o(1)
$$
with $o(1)$ is small with $\alpha.$  Since the functions $\hat{z}_i$ are linearly independent on any open set because they solve an homogeneous
elliptic PDE, we conclude that the matrix with the above coefficients is invertible. Hence from \equ{sis} and (\ref{normalized}), all $\ttt c_i$'s are necessarily zero.
We have thus found a solution to the Allen Cahn equation \equ{ac} with the properties required in Theorem \ref{teo1}.

\medskip
It remains to prove identities (\ref{ee1-1}). The idea is to use the invariance of $ \Delta + f(u)$ under
rigid translations and rotations. This type of Pohozaev identity argument has been used in a number of places, see for instance \cite{gui-cross}.

\medskip
In order to prove that the identity  \equ{ee1-1}  holds for $i=3$, we consider a large number $R>> \frac 1\A$ and the infinite cylinder
$$C_{R} = \{ x \ /\ x_1^2 + x_2^2 < R^2 \}. $$
Since in $C_R$ the quantities involved in the integration approach zero at exponential rate as  $|x_3|\to +\infty$
uniformly in $(x_1,x_2)$, we have that

$$
\int_{C_R } (\Delta  u_{}  + f(u_{}) )\partial_{x_ 3} u_{} -  \int_{\partial C_R}   \nn u_{}\cdot \hat r\,   \partial_{x_3}u_{}    = \int_{C_R}   \partial_{x_3} \, (\,  F(u_{}) - \frac 12 |\nn u_{}|^2 \, ) = 0.
$$

We claim that
$$
\lim_{R\to +\infty } \int_{\partial C_R}   \nn u_{}\cdot \hat r\,  \, \partial_{x_3}u_{}    \, =\, 0 .
$$

Using estimate \equ{for} we have that near the manifold,
$$
 \partial_{x_3}u_{}  \nn u_{} (x) \cdot \hat r =    w'(t)^2 ((\nu -\alpha (-1)^k \frac{\beta_k}{r_{\A}} \hat{r}) \cdot \hat r) \nu_3  + O( \alpha e^{-\sigma |t|} \frac 1{ r^2}).
 $$
Let us consider the $k$-th end, which for large $r$ is expanded as
$$
x_3 =  F_{k,\A}( x_1,x_2) = \alpha^{-1} (a_k \log \alpha r + b_k + O(r^{-1}))
$$
so that
\begin{equation}
(-1)^k \nu =  \frac 1{\sqrt{ 1 + |\nn F_{k,\A}|^2}} ( \nn F_{k,\A} , -1)  =
  \, \frac {a_k}\A \frac{\hat r} r  - e_3 \, +\,  O( r^{-2})\, .
\label{normal} \end{equation}
Then on the portion of $C_R$ near this end we have that
\begin{equation}
(\nu -\alpha (-1)^k \frac{\beta_k}{r_{\A}} \hat{r})
\cdot \hat r\, \nu_3  =
- \alpha^{-1} \frac {a_k +\alpha \beta_k }{ R} + O(R^{-2})  . \ee
In addition, also, for $x_1^2 + x_2^2 = R^2$  we have the expansion
$$
t  = (x_3 - F_{k,\A}(x_1,x_2) -  \beta_k \log\A r + O(1) ) ( 1+ O(R^{-2}))
$$
with the same order valid after differentiation in $x_3$, uniformly in such $(x_1,x_2)$.
Let us choose $\rho = \gamma \log R $ for a large, fixed  $\gamma$. Observe that on $\partial C_R$ the distance between ends is greater than $2\rho$ whenever $\A$ is sufficiently small. 
We get,
$$
\int_{F_{k,\A}(x_1,x_2) + \beta_k \log\A r- \rho}^{F_{k,\A}(x_1,x_2)+  \beta_k \log \A r + \rho}  w'(t)^2 dx_3 =  \int_{-\infty}^{\infty} w'(t)^2 dt + O(R^{-2}) \, .
$$

Because of estimate \equ{end11}
we conclude,
fixing appropriately $\gamma$,
that
 $$
 \int_{ \bigcap_k \{|x_3 - F_{k,\A}|> \rho\} }   \partial_{x_3}u_{}  \nn u_{} (x) \cdot \hat r \, dx_3 = O(R^{-2})\, .
$$
 As a conclusion
$$
\int_{-\infty}^\infty   \partial_{x_3}u_{}  \nn u_{} \cdot \hat r \, dx_3  = - \frac {1}{\A R}\sum_{k=1}^m  (a_k +\A \beta_k)
\int_{-\infty}^\infty w'(t)^2\, dt\,   + O( R^{-2})
$$
and hence
$$
\int_{\partial C_R} \partial_{x_3}u_{}  \nn u_{} (x) \cdot \hat r  = -\frac {2\pi} {\A} \sum_{k=1}^m  (a_k+\A \beta_k)   + O(R^{-1}) \, .
$$
But $\sum_{k=1}^m  a_k = \sum_{k=1}^m \beta_k =0$ and hence  \equ{ee1-1} for $i=3$ follows after letting $R\to \infty$.

\medskip
Let us prove the identity for $i=2$.  We need to carry out now the integration against $\partial_{x_2} u_{}$. In this case we get

$$
\int_{C_R } (\Delta  u_{}  + f(u_{}) )\partial_{x_ 2} u_{}  = \int_{\partial C_R}   \nn u_{}\cdot \hat r\,   \partial_{x_2}u_{}    +  \int_{C_R}   \partial_{x_2} \, (\,  F(u_{}) - \frac 12 |\nn u_{}|^2 \, ) .
$$
We have that
$$
 \int_{C_R}   \partial_{x_2} \, (\,  F(u_{}) - \frac 12 |\nn u_{}|^2 \, ) = \int_{\partial C_R} (\,  F(u_{}) - \frac 12 |\nn u_{}|^2 \, )  n_2
$$
where $n_2 = x_2/r$.
Now, near the ends estimate \equ{for} yields
$$
|\nn u_{}|^2 = |w'(t)|^2 + O( e^{-\sigma |t|} \frac 1{ r^2})
$$
and arguing as before, we get
$$
\int_{-\infty}^\infty  |\nn u_{}|^2 dx_3 = m \int_{-\infty}^\infty  |w'(t)|^2dt  + O( R^{-2}) .
$$
Hence
$$
\int_{\partial C_R} |\nn u_{}|^2 n_2  =  m  \int_{-\infty}^\infty  |w'(t)|^2dt \, \int_{[ r= R ] } n_2  +  O( R^{-1}) \, .
$$
Since
$
\int_{[r= R]} n_2 = 0
$
we conclude that
$$
\lim_{R\to +\infty} \int_{\partial C_R} |\nn u_{}|^2 \, n_2  = 0 .
$$
In a similar way we get
$$
\lim_{R\to +\infty} \int_{\partial C_R} F( u_{}) \, n_2  =0 .
$$
Since  near the ends we have
 $$ \partial_{x_2}u_{} =  w'(t)(\nu_2 -\A (-1)^k \frac{\beta_k}{r_{\A}} \hat{r} e_2)   +  O( \alpha r^{-2} e^{-\sigma |t|}) $$
 and from \equ{normal} $\nu_2 = O( R^{-1})$,  completing the computation as previously done yields
 $$
\int_{\partial C_R}   \nn u_{}\cdot \hat r\,   \partial_{x_2}u_{}  = O( R^{-1}).
$$
As a conclusion of the previous estimates, letting $R\to +\infty$  we finally find the validity of \equ{ee1-1} for $i=2$. Of course the same
argument holds for $i=1$.

\medskip
Finally, for  $i=4$
it is convenient to compute the integral over $C_R$ using cylindrical coordinates. Let us write
$u_{} = u(r,\theta, z)$.
Then

$$
\int_{C_R } (\Delta  u_{}  + f(u_{}) ) \,(x_2 \partial_{x_1} u_{}  - x_1 \partial_{x_1} u_{}) =
$$
$$
\int_0^{2\pi} \int_0^R \int_{-\infty}^\infty  [ u_{zz} +  r^{-1} (ru_r)_r + f(u) ]\, u_\theta  \, r\,d\theta \,dr \,dz\, \  =
$$
$$
- \frac 12 \int_0^{2\pi} \int_0^R \int_{-\infty}^\infty   \, \partial_\theta\,  [ u_{z}^2 +   u_r^2 - 2 F(u)   ]\,  r\,d\theta \,dr \,dz\,   +
 R \int_{-\infty}^{\infty} \int_0^{2\pi}  u_r\, u_\theta (R,\theta, z)\, d\theta\, dz \, = \, $$
 $$
0  + \int_{\partial C_R} u_ru_\theta \ .
$$
 On the other hand, on the portion of $\partial C_R$ near the ends we have
$$
u_r\, u_\theta  = w'(t)^2  R (\nu\cdot \hat r )( \nu\cdot \hat \theta)  +   O(  R^{-2} e^{-\sigma |t|}).
$$
From \equ{normal}  we find
$$ (\nu\cdot \hat r )( \nu\cdot \hat \theta) = O( R^{-3}), $$
hence
$$
u_r\, u_\theta  = w'(t)^2  O(R^{-2}) +   O(  R^{-2} e^{-\sigma |t|})
$$
and finally   $$\int_{\partial C_R} u_r\, u_\theta  = O( R^{-1}). $$
Letting $R\to +\infty$ we obtain relation  \equ{ee1-1} for $i=4$.
The proof is concluded. \qed

\bigskip

\section{  \emph{Negative eigenvalues and their eigenfunctions  for the Jacobi operator}}\label{eigenvaluejacobi}

For the proof of Theorem \ref{teo2} we need to translate the information on the index of the minimal surface $M$ into spectral features of the Jacobi
operator.
Since $M$ has finite total curvature, the index $i(M)$ of the minimal surface $M$ is finite.  We will translate this information
into an eigenvalue problem for the operator $\JJ$. Let
$$
{\tt Q}(k,k) := \int_M |\nn k|^2\, dV \, -\, \int_M |A|^2 k^2\, dV \ .
$$
The number $i(M)$ is, by definition,  the largest dimension for a vector space  $E$ of compactly supported smooth functions in $M$ such that
$$
{\tt Q}(z,z)  < 0 \foral z\in E\setminus \{0\} .
$$

The number $i(M)$ when finite has the following convenient characterization, whose proof is straightforward.
In what follows we fix the function  $q$ as
\be
p(y)\  :=\ \frac 1{1 + r(y)^4}  \ .
\label{py1}\ee
Let
us consider for a large number $R$, the region
$$ M^R = \{y\in M\ / \ r(y) < R\} $$
and the eigenvalue problem
\be
\Delta_M k + |A|^2 k +  \la p(y)\, k = 0 \quad \hbox{ in } M^R,
\label{eigjac}\ee
$$
k=0 \quad \hbox{ on } \partial M^R\ .$$
Let $m_R (p)$ denote the number of negative eigenvalues (counting multiplicities) for this problem. Then we have
\be
i(M) \ =\ \sup_{R>0}  m_R(p) \ .
\label{i(M)}\ee
Let us also consider the eigenvalue problem in entire space
\be
\Delta_M k + |A|^2 k +  \la p(y)\, k = 0 \quad \hbox{ in } M,\quad k\in L^\infty(M)\, .
\label{eigjac1}\ee
We will prove the following result.

\begin{lemma}\label{lemin8}
Problem $\equ{eigjac1}$ has exactly $i(M)$ negative eigenvalues, counting multiplicities.
\end{lemma}

\bigskip
\subsubsection{ A priori estimates in  $M^R$}
For the proof of Lemma \ref{lemin8}, and for later purposes, it is useful to have a priori estimates
uniform in large $R>0$ for the linear problem
\be
  \Delta_M k + |A|^2 k  -\gamma p(y)\, k = f \quad \hbox{ in } M^R,
\label{linjac}\ee
$$
k=0 \quad \hbox{ on } \partial M^R\ .$$

We have the following result.
\begin{lemma} \label{lemin3}
Let $p>1$, $\sigma > 0$. Then
 for  $R_0>0 $  large enough and fixed and $\gamma_0>0$, there exist a $ C>0$ such that for all $R> R_0+1$, $0\le \gamma <\gamma_0$,
any  $f$,  and any  solution  $k$ of problem $\equ{linjac}$, we have that

$(a)$  If
$\|f\|_{{p,4-\frac 4p}} <+\infty$
then
 \be
 \| k\|_{\infty}
\ \le\
C \, [ \, \|f\|_{p, 4-\frac 4p } \, + \, \|k\|_{L^\infty( |y|< 3R_0)}\  ] \, .
\label{e}\ee

\medskip
$(b)$  If  $\|f\|_{{p, 2-\frac 2p -\sigma }} <+\infty$, then
,
\be
\|  D^2k\|_{p, 2-\frac 2p -\sigma  } + \| \,Dk\, \|_{p, 1-\frac 2p -\sigma  }\le
C \, [ \, \|f\|_{p, 2-\frac 2p -\sigma }  \, + \, \|k\|_\infty \  ] \, .
\label{oi}\ee

\quad If $p>2$, we have in addition
\be \| \, (1+ |y|)^{1-\sigma} \,Dk\, \|_{\infty} \le
C \, [ \, \|f\|_{p, 2-\frac 2p -\sigma }  \, + \, \|k\|_\infty \  ] \, .
\label{oii}\ee

\end{lemma}

\proof
Let us consider the equation in $M$
\begin{equation}
\Delta_M \psi   + |A|^2\psi  = -|f|\chi_{|y|<R}, \quad    |y| >  R_0 \, ,
\label{i1}\end{equation}
\be
 \psi (y)=  0 , \quad   |y| =  R_0\, .
\label{i2}\ee
For a large and fixed $R_0$,  solving this problem amounts to doing it on each separate end. As in Lemma \ref{lem8.1}, after a Kelvin's transform
the problem reduces in each end  to solving in a ball in $\R^2$ an equation of the form
 $$
\Delta \ttt \psi + O(|\by|^2 )  D^2\ttt \psi + O(|\by|)D \ttt \psi + O(1)\ttt \psi   = -|\ttt f|\chi_{|\py|>\frac 1R},\quad  |\by| <  \frac 1{R_0} \, ,
$$
$$
 \ttt\psi  (\by)=  0 , \quad     |\by| =  \frac 1 R_0\, .
$$
Enlarging $R_0$ if necessary, this problem has a unique solution, which is also positive. 
 This produces a bounded, positive  solution $\psi$ of
\equ{i1}-\equ{i2} with
$$
\|\psi\|_\infty \le C\|f\|_{p, 4-\frac 4p}\ .
$$
On the other hand, on this end the Jacobi field $z_3 = \nu\cdot e_3$ can be taken positive with $z_3\ge 1$ on $|\py| >R_0$.
Thus the function $\psi + \|k\|_{L^\infty ( |y|= R_0)} z_3 $
is a positive, bounded supersolution for the problem
\equ{linjac} in this end, where $|\py| > R_0$,
and estimate \equ{e} then readily follows.

\medskip
Let us prove now estimate  \equ{oi}.
Fix a large number $R_0>0$ and another number $R>> R_0$.
Consider also a large $\rho>0$ with $3\rho < R$. On a given end we parameterize with Euclidean coordinates $\py \in \R^2$ and get that the
equation satisfied by $k = k(\py )$ reads
$$
\Delta k + O(|\py|^{-2})D^2 k  +  O(|\py|^{-3}) D k + O(|\py|^{-4}) k  =  f, \quad  R_0 < |\py| < R .
$$

Consider the function $k_\rho ( z) = k (\rho z ) $ wherever it is defined. Then
$$
\Delta k_\rho  + O(  \rho^{-2} |  z|^{-2})D^2 k_\rho   +  O( \rho^{-2} |z |^{-3}) D k_\rho   + O(\rho^{-2} |z|^{-4} ) k_\rho =  f_\rho
$$
where  $f_\rho (z) = \rho^{2} f(\rho z) $.
Then interior elliptic estimates (see Theorem 9.11 of \cite{GT})  yield the existence of a constant $C= C(p)$ such that for any sufficiently large $\rho$
\be
\| D k_\rho  \|_{L^p( 1< |z| <2)}+ \| D^2 k_\rho  \|_{L^p( 1< |z| <2)} \, \le  \,  C\, (\, \|k_\rho\|_{L^\infty ( \frac 12 < |z| < 3)}  +  \| f_\rho  \|_{L^p( \frac 12 < |z| < 3)}\, ).
\label{jjj}\ee
Now, 
$$
 \| f_\rho  \|^p_{L^p( \frac 12 < |z| < 3) }\,  =   \rho^{2p} \int_{ ( \frac 12 < |z| < 3)} | f(\rho z)|^p \, dz  \le
$$
$$
C\rho^{p\sigma}\, \int_{ ( \frac 12 < |z| < 3)}  |\rho z|^{2p-2-p\sigma} |f(\rho z)|^p \, \rho^2  dz \, =\,
 C\, \rho^{p\sigma}\,\int_{ ( \frac \rho 2 < |\py| < 3\rho )}|\py|^{2p-2-p\sigma} |f(\py)|^p \,d\py \ .
$$
Similarly
$$
\| D^2 k_\rho  \|^p_{L^p( 1< |z| <2)} \ge   C\, \rho^{p\sigma}\,\int_{ ( \rho  < |\py| < 2\rho )}|\py|^{2p-2-p\sigma} |D^2k(\py)|^p \,d\py \ .
$$
Thus
$$
 \,\int_{ ( \rho  < |\py| < 2\rho )}|\py|^{2p-2-p\sigma} |D^2k(\py)|^p \,d\py\, \le\,  C\,
 \int_{ ( \frac \rho 2 < |\py| < 4\rho )}|\py|^{2p-2-p\sigma} |f(\py)|^p \,d\py \ + \rho^{-p\sigma} \|k\|_\infty^p .
$$
Take  $\rho= \rho_j = 2^j$. Then
$$
\int_{ ( \rho_j  < |\py| < \rho_{j+1} )} |\py|^{2p-2-p\sigma} |D^2k(\py)|^p \,d\py\, \le\,
$$
$$
C\,
 \int_{ (  \rho_{j-1} < |\py| < \rho_{j+2} )}|\py|^{2p-2-p\sigma} |f(\py)|^p \,d\py \ + 2^{-j  p\sigma} \|k\|_\infty^p\, .
$$
Then, adding up these relations wherever they are defined, taking in addition into account
 boundary elliptic estimates which give that for $\rho = \frac R2$,
$$
\| D^2 k_\rho  \|_{L^p( 1< |z| <2)} \, \le  \,  C\, \left
(\, \|k_\rho\|_{L^\infty ( \frac 12 < |z| < 2)}  +  \| f_\rho  \|_{L^p( \frac 12 < |z| < 2)}\,\right )\, ,
$$
plus a local elliptic estimate in a bounded region, we obtain that for some $C>0$ independent of $R$,
$$
\| D^2k\|_{p, 2-\frac 2p -\sigma} \le \,  C\, (\, \|k\|_\infty  + \| f\|_{p, 2-\frac 2p -\sigma} \, ).
$$
The corresponding estimate for the gradient follows
immediately from \equ{jjj}. We have proven \equ{oi}. If $p>2$ we can use Sobolev's embedding to include $\| Dk_\rho\|_{L^\infty ( 1< |z|<2)}$ 
on the left hand side
of \equ{jjj}, and estimate \equ{oii} follows.
The proof is complete. \qed
\medskip

\subsubsection{Proof of Lemma $\ref{lemin8}$}
We will prove first that problem \equ{eigjac1} has at least $i(M)$ linearly independent eigenfunctions
associated to negative eigenvalues
in $L^\infty(M)$.
For all $R>0$ sufficiently large, problem \equ{eigjac} has $ n= i(M)$ linearly independent eigenfunctions
$ k_{1,R},\ldots, k_{n,R}$ associated to negative eigenvalues
$$
\la_{1,R}\le\la_{2,R}\le\cdots \le \la_{n,R} <0 \ .
$$
Through the min-max characterization of these eigenvalues, we see that they can be chosen to define decreasing functions of $R$. On the other hand,
$\la_{1,R}$ must be bounded below. Indeed, for a sufficiently large $\gamma>0$ we have that
$$
|A|^2 -\gamma p \ < \ 0 \quad\hbox{in } M
$$
and by maximum principle we must have $\la_{1,R} > -\gamma$. The eigenfunctions can be chosen orthogonal in the sense that
\be
\int_{M^R}  p\, k_{i,R}\, k_{j,R} \, dV\ = \ 0 \foral i\ne j  \ .
\label{orto}\ee
Let us assume that $\|k_{i,R}\|_\infty = 1$.
Then the a priori estimate in Lemma \ref{lemin3} imply that, passing to a subsequence in $R\to +\infty$, we may assume that
$$ \la_{i,R}\downarrow \la_{i} < 0 ,\qquad k_{i,R}(y)\to k_{i}(y), $$
uniformly on compact subsets of $M$, where $k_{i} \ne 0$ is a bounded eigenfunction of \equ{eigjac1} associated
to the negative eigenvalue $\la_{i}$.
Moreover, relations \equ{orto} pass to the limit and yield
\be
\int_{M}  p\, k_{i}\, k_{j} \, dV\ = \ 0 \foral i\ne j  \ .
\label{orto1-1}\ee
Thus, problem \equ{eigjac1} has at least $n=i(M)$ negative eigenvalues. Let us assume there is a further bounded eigenfunction $k_{n+1}$,
linearly independent of $k_1,\ldots, k_n$, say with
\be
\int_{M}  p\, k_{i}\, k_{n+1} \, dV\ = \ 0 \foral i =1,\ldots, n \ ,
\label{orto1}\ee
associated to a negative eigenvalue $\la_{n+1}$. Then the a priori estimate of Lemma \ref{lem8.1} implies that
$$ \|(1+r^2) \nn k_{n+1}\|\ < +\infty . $$ The same of course holds for the remaining $k_i$'s.
It follows that
$$ {\tt Q}(k, k) < 0 \foral  k\in {\rm span}\, \{ k_1,\ldots , k_{n+1}\}\setminus \{0\} .$$
However, again since $\nabla k_j$ decays fast, the same relation above will hold true for the $k_i$'s replaced by suitable smooth
truncations far away from the origin. This implies, by definition, $i(M) \ge n+1$ and we have have reached a contradiction. The proof  is concluded.
\qed

\bigskip

\setcounter{equation}{0}
\section{ \emph{The proof of Theorem \ref{teo2}}}

In this section we will prove that the Morse index $m( u_\A)$ of the solution we have built in Theorem \ref{teo1} coincides with the index of the surface $M$, as stated in Theorem \ref{teo2}.
We recall that this number is defined as the supremum of all dimensions of vector spaces $E$
of compactly supported smooth functions
for which
$$
\QQ( \psi,\psi ) \ =\ \int_{\R^3} |\nn \psi|^2 - f'(u_\A)\psi^2 <0 \foral \psi \in E\setminus \{0\}   .
$$
We provide next a more convenient characterization of this number, analogous to that for the Jacobi operator of \S \ref{eigenvaluejacobi}.
Let us consider a smooth function $p(x)$ defined in $\R^3$ such that 
$$
p(\A x) = \frac 1{1+ r_\A(y)^4 }\quad \hbox{if }\quad x = y+ (t+ h(\A y))\nu(\A y)\,   \in \NN_\delta,\quad
$$
and such that for positive numbers $a,b$,
$$
\frac a{1+ |\A x'|^4} \ \le \     p( \A x) \ \le \ \frac b{1+ |\A x'|^4} \quad \foral x =(x', x_3) \in \R^3\ .
$$
For each $R>0$, we consider the eigenvalue problem
in the cylinder
$$
{\mathcal C}_R \,=\, \{\,(x',x_3)\ /\ |x'| <  R\A^{-1},\  |x_3| < R\A^{-1})\,\},
$$

\begin{equation}
\Delta \phi+ f'(u_\alpha ) \phi  + \la  p (\A x) \phi = 0 \quad\hbox{in   } {\mathcal C}_R \, ,
\label{R1}\end{equation}
$$
\phi = 0   \quad\hbox{on  }  \partial {\mathcal C}_R\ .
$$
We also consider the problem in entire space
\begin{equation}
\Delta \phi+ f'(u_\alpha ) \phi  + \la  p (\A x) \phi = 0 \quad\hbox{in  }\R^3,\quad \phi\in L^\infty(\R^3).
\label{R2}\end{equation}

Let $m_R ( u_\A)$ be the number of negative eigenvalues $\la$ (counting multiplicities) of this Problem \equ{R1}. Then we readily check that
$$
m(u_\A) = \sup_{R>0} m_R(u_\A) .
$$
On the other hand, we have seen in \S \ref{eigenvaluejacobi} that the index $i(M)$ of the minimal surface can be characterized as
 the number of linearly independent eigenfunctions associated to negative eigenvalues of the problem
\begin{equation}
\Delta  z  + |A|^2z  + \la  p(y) z  = 0 \quad\hbox{in }  M,\quad z\in L^\infty (M)\, ,
\label{M1}\end{equation}
which corresponds to the maximal dimension of the negative subspace in $L^\infty(M)$ for the quadratic form
$$
{\tt Q}(z,z)\, =\, \int_M |\nn_M z|^2 - |A|^2 z^2 \, dV\, .
$$
We shall prove in this section that $m(u_\A) = i(M)$ for any sufficiently small $\A$.

\medskip
The idea of the proof is to put in correspondence  eigenfunctions for negative eigenvalues of problem \equ{R1} for large $R$
with those of problem \equ{M1}. This correspondence comes roughly as follows. If $z$ is such an eigenfunction for problem \equ{M1} then the function
defined near $M_\A$ as
\be
k(y)w'(t),\quad k(y)= z(\A y)
\label{kw}\ee
defines after truncation a negative direction for the quadratic form $\QQ$ on any large ball. Reciprocally, an eigenfunction for negative eigenvalue of
problem \equ{R1} will look for any sufficiently small $\A$ and all large $R$ like a function of the form \equ{kw}.
In the following two lemmas we clarify the action of the operator $L$ on functions of this type, and the corresponding connection at the
level of the quadratic forms $\QQ$ and ${\tt Q}$.

\begin{lemma}\label{lemin10}
Let $k(y)$ be a function of class $C^2$ defined in some open subset ${\mathcal V}$ of $M_\A$.
Let us consider the function $v(x)$ defined for $x \in \NN_\delta$, 
 $y\in {\mathcal V}$ as
$$
v(x)\, =\, v(y,t)\ :=\ k(y)\,w'(t)\, , \quad y\in  {\mathcal V},\quad |t+ h_1(\A y)|\,<\, \rho_\A (y)
$$
where $\rho_\A$ is the function in the definition of $\NN_\A$, $\equ{region}$.
Then $L(v):= \Delta_xv + f'(u_\A)v$ can be expanded as in $\equ{L(v)}$ below. Besides we have
$$
\int_{|t +h_1 |< \rho_\A} \, L(v)\, w'\, dt\ = \
 \ (\,\Delta_{M_\A} k +  \A^2 |A|^2 k\, + \,  \A h a_{ij}^{1,0}\partial_{ij}k \, )\, \int_\R{w'}^2\, dt
$$
\be
+\  O(\A^2 r_\A^{-2})\,\partial_{ij}k \,+ \, O(\A^2 r_\A^{-3})\, \partial_i k \, + O(\A^3 r_\A^{-4})\, k
\, .
\label{r7}\ee
Here  $$a_{ij}^{1,0} = a_{ij}^{1,0}(\A y) = O( r_\A^{-2})\ . $$
The same conclusions hold for the function
$$
v(x)\, =\, v(y,t)\ :=\ k(y)\,w'(t)\, \eta_\delta (y,t)  , \quad y\in  {\mathcal V},\quad |t+ h_1(\A y)|\,<\, \rho_\A (y)
$$
where the cut-off function $\eta_\delta$ is defined in
$\equ{etadelta}.$
\end{lemma}

\proof

Let us recall that
$$
\Delta_x   =  \, \partial_{tt}  + \Delta_{M_\A}  - \A^2[ (t+h)|A|^2 + \Delta_M h]  \partial_t\ -  2\A\, a_{ij}^0\, \partial_jh \partial_{it} \ +
$$
$$
 \A (t+h) \, [ a_{ij}^1\partial_{ij}   -  2\A\, a_{ij}^1\, \partial_jh \partial_{it}  + \A (b_i^1\partial_i   - \A b_i^1 \partial_ih \partial_t)\, ) \, ] \ +
 $$
$$
  \A^3 (t+h)^2b^1_3  \partial_t\, +\, \A^2 [\,a_{ij}^0 + \A(t+h) a_{ij}^1 )\,]\partial_ih\partial_jh \, \partial_{tt} \, .
$$
Hence, using Lemma \ref{formula} in the appendix
we get
$$
\Delta_x v +f'(u_\A) v\  =  \,  k (w'''  + f'(w)w') \, + \, [f'(u_\A)-f'(w)]\,kw'\ +
$$

$$
 w'\, \Delta_{M_\A}k  \,-\,  \A^2[ (t+h_1)|A|^2 + \Delta_M h_1]\,k\, w'' \, - \,2 \A\, a_{ij}^0 \partial_jh \partial_{i}k \, w'' \ +
$$

$$
 \A (t+h) \, [ a_{ij}^1\partial_{ij}k w'  - \A\, a_{ij}^1(\, \partial_jh \partial_{i}k + \partial_ih \partial_{j}k)\, w''\,  +
 \A (b_i^1\partial_ik\,w'   - \A b_i^1 \partial_ih w'')\, \, ] \ +
 $$

\be
  \A^3 (t+h)^2b^1_3 \, k\, w'' \, +\, \A^2 [\,a_{ij}^0 + \A(t+h) a_{ij}^1 )\,]\partial_ih\partial_jh \,k\, w'''\ .
\label{r8}\ee

\medskip
We can expand
$$
a_{ij}^1  = a_{ij}^1 (\A y,0)  + \A (t+h)\,a_{ij}^2 (\A y, \A (t+h)) =: a_{ij}^{1,0} + \A (t+h)a_{ij}^2,
$$
with $a_{ij}^2 = O( r_\A^{-2})$, and similarly
$$
 b_{j}^1 = b_j^1 (\A y,0)  + \A (t+h)\,b_{j}^2 (\A y, \A (t+h))=: b_j^{1,0} + \A (t+h)b_{j}^2,
 $$
 with $b_{j}^2 = O( r_\A^{-3})$.
On the other hand,  let us recall that
$$ u_\A - w = \phi_1 + O(\A^3 r_\A^{-4} e^{-\sigma|t|}) $$
where $\phi_1$ is given by \equ{u1},
\be
\phi_1(y,t) =  \A^2 |A(\A y)|^2 \psi_0 (t)  - \A^2 a_{ij}^0 \pp_i h_0\pp_j h_0 (\A y)\, \psi_1 (t)
\label{phi1}\ee
and $\psi_0$, $\psi_1$ decay exponentially as $|t|\to +\infty$. Hence
$$
[f'(u_\A)-f'(w)]\,w' \ =\ f''(w)\,\phi_1\, w' \ +\ O( \A^3 e^{-\sigma|t|} r_\A^{-4}).
$$
\medskip
Using these considerations and expression \equ{r8}, we can write,
$$
Q\ := \ \Delta_x v +f'(u_\A) v\  =
$$

$$
\, \underbrace{ \Delta_{M_\A}k\, w'\, - \, \A^2 |A|^2\,k\, t w'' + \A^2 \,a_{ij}^0\partial_ih_0\partial_jh_0 \,k\, w''' + \A h a_{ij}^{1,0}\partial_{ij}k  w'}_{Q_1}\, +\,
  \underbrace{ f''(w) \,\phi_1 \,kw'\ }_{Q_2}
$$

$$
\underbrace{ -w''\, \left [   \A a_{ij}^0( \partial_jh \partial_{i}k + \partial_ih \partial_{j}k) + \A^2k \Delta_M h_1
  +
  \A^2 h  a_{ij}^{1,0}( \partial_jh \partial_{i}k + \partial_ih \partial_{j}k)
    \right]}_{Q_3}
    $$

    $$
  \  +\    \underbrace{
    \A t w'\,\left [ a_{ij}^{1,0}\partial_{ij}k  + \A b_i^{1,0}\partial_ik \right ]}_{Q_4}
    $$

    \be \, +\,  \underbrace{ \A^2(t+ h )^2a_{ij}^{2}\partial_{ij}k  w'  +\A^2 (t+h)  a_{ij}^2(\, \partial_jh \partial_{i}k + \partial_ih \partial_{j}k)\, w''}_{Q_5}
 \,+\,
 \underbrace{ O(\A^3e^{-\sigma|t|} r_\A^{-2}) }_{Q_6} .\label{L(v)}\ee
The precise meaning of the remainder $Q_6$ is
$$
Q_6 \ =\   O(\A^3e^{-\sigma|t|} r_\A^{-2})\, \partial_{ij}k \,+ \, O(\A^3e^{-\sigma|t|} r_\A^{-3})\,\partial_{j}k\, .
$$
\bigskip
We will integrate the above relation against $w'(t)$ in the region  $|t+ h_1(\A y)| < \rho_\A(y)$.

\medskip
Let us observe that the terms $Q_i$ for $i=1,\ldots, 4$ are in reality defined for all $t$
and that
\be
\int_{|t +h_1| < \rho_\A} Q_i\,w'\, dt\,  =\,  \int_\R Q_i\,w'\, dt\, + O(\A^3 r_\A^{-4})
\label{r6}\ee
where the remainder means
$$
O(\A^3 r_\A^{-4})\ :=\ O(\A^3 r_\A^{-4})\,\partial_{ij}k + O(\A^3 r_\A^{-4})\,\partial_i k +\,O(\A^3 r_\A^{-4})\,k \ .
$$
Let us observe that
\be
\int_\R\, (Q_3\,+\, Q_4)\, w'\, dt\ = \ 0\ .
\label{r1}\ee
On the other hand, since
$$
\int_\R tw''w'\, dt \ = \ -\, \frac 12 \int_\R {w'}^2\, dt\, ,
$$
we get that
\be
\int_\R\, Q_1 \, w'\, dt\ =  \ (\,\Delta_{M_\A} k + \frac 12 |A|^2 k\, + \,  \A h a_{ij}^{1,0}\partial_{ij}k \, )\, \int_\R{w'}^2\, dt
\, +\, a_{ij}^0\pp_i h_0\pp_jh_0 \, \int_\R w'''\,{w'}\, dt \,.
\label{r2}\ee\
Next we will compute $\int_\R\, Q_2 \, w'\, dt$.
We recall that, setting $L_0(\psi) = \psi'' + f'(w)\psi$, the functions $\psi_0$ and $\psi_1$ in \equ{phi1} satisfy
$$
L_0(\psi_0) = tw'(t),\quad L_0(\psi_1) = w''\ .
$$  Differentiating these equations we get
$$
L_0( \psi_0') + f''(w)w'\psi_0 = (tw')'  ,\quad L_0( \psi_1') + f''(w)w'\psi_1 = w''' \, .
$$
Integrating by parts against $w'$, using $L_0(w')=0$ we obtain
$$
\int_\R f''(w){w'}^2\psi_0    =  - \int_\R tw''w'  = \frac 12 \int {w'}^2  ,\quad \int_\R f''(w){w'}^2\psi_1 = \int_\R w'''w' \, .
$$
Therefore
$$
\int_\R  Q_1\, w'\, dt\, = \,   \int_\R  f''(w) \phi_1 \,k{w'}^2\, dt\, =
$$
$$
 \A^2k |A|^2 \int_\R  f''(w) \psi_0 \,{w'}^2\, dt\, - \A^2 a_{ij}^0 \pp_i h_0\pp_j h_0 \,k\,
\int_\R  f''(w) \psi_1 \,{w'}^2\, dt\, =
 $$
\be
  \A^2k |A|^2  \frac 12 \int_\R {w'}^2   - \A^2 a_{ij}^0 \pp_i h_0\pp_j h_0\, k\, \int_\R w'''w'.
\label{r3}\ee

Thus,  combining relations \equ{r1}-\equ{r3} we get
\be
\int_\R\, (\,Q_1+ \cdots + Q_4) \, w'\, dt\ =  \ (\,\Delta_{M_\A} k +  |A|^2 k\, + \,  \A h a_{ij}^{1,0}\partial_{ij}k \, )\, \int_\R{w'}^2\, dt
\, .
\label{r4}\ee

On the other hand, we observe that
\be
\int_{|t +h_1 |< \rho_\A} \, (Q_5 +Q_6)\, w'\, dt\ = \  O(\A^2 r_\A^{-2})\,\partial_{ij}k \, + \,O(\A^2 r_\A^{-3})\, \partial_i k \
\, .
\label{r5}\ee
Combining relations \equ{r4}, \equ{r5} and \equ{r6}, expansion \equ{r7} follows.
Finally, for $v$ replaced by $\eta_\delta k\, w'$
we have that
$$
\int L(kw\eta_\delta)\eta_\delta \,kw'\, dt \ = \, \int \eta_\delta^2 L(kw)kw\, dt +  \int \eta_\delta (\Delta\eta_\delta  \, kw'+ 2\nabla \eta_\delta \nabla (kw') )\, kw\, dt\ .
$$
The arguments above apply to obtain the desired expansion for  the first integral in the right hand side of the above decomposition. The second
integral produces only smaller order operators in $k$ since $\Delta\eta_\delta$, $\nabla \eta_\delta $ are both of order $O( r_\A^{-4}\A^4)$ inside their supports. The proof is concluded.    \qed

\bigskip
Let us consider now the region
$${\mathcal W} := \{ x\in \NN_\delta \ /\ r_\A(y)< R \},$$
where $R$ is a given large number.

\begin{lemma}\label{lemin11}
Let $k(y)$ be a smooth function in $M_\A$  that vanishes when $r_\A(y) = R$, and
set $v(y,t) := \eta_\delta (y,t)\, k(y)\, w'(t)$. Then
the following estimate holds.
  $$
 \QQ(v,v)   =\ \int_{{\mathcal W}}
  |\nn v\, |^2 \, -\,  f'(u_\A) \,v ^2\, dx \, = \,
  $$
  $$
 \int_{ r_\A(y)< R }\left [ \, |\nn_{M_\A} k|^2 - \A^2 |A(\A y)|^2\, k^2\, \right ]\, dV_\A  \int_\R {w'}^2\, dt
 $$
 \be
 \ +\  O\, \left ( \A \, \int_{ r_\A(y)< R } [\,\,|\nn k|^2 \, +  \, \A^2\,(1+ r_\A^{4})^{-1}\, k^2 \,]\, \,dV_\A \, \right )\, .
 \label{qq}\ee
\end{lemma}

\medskip
\proof
Let us estimate first the quantity
$$
 \int_{{\mathcal W}} L(kw')\, kw'\, dx\ .
$$
Let us express the Euclidean element of volume $dx$ in the coordinates $(y,t)$. Consider one of the charts
 $Y_l(\py),\ \py \in {\mathcal U}_l$ of $M$, $l=1,\ldots,N$ introduced in \equ{chart}, which induce corresponding charts in $M_\A$
 as in \equ{chart2}. Then, dropping the index $l$,
we  compute the element of volume through the change of coordinates
\be
x = X(\py, t)=   \A^{-1} Y(\A \py)  + (t+h) \nu(\A \py ) \ .
\label{cv}\ee
Below, as in subsequent computations, the computation of the integral in the entire region is performed by localization through smooth
partition of unity $\xi_1,\ldots, \xi_m$ subordinated to the covering $Y_l({\mathcal U}_k)$, $l=1,\ldots, N$ of $M$,
namely with the support of $\xi_l$ is contained in $Y_l({\mathcal U}_k)$ and $\sum_l \xi_l \equiv 1$.
We perform typically a computation of an integral of a function $g(x)$ defined in $\NN_\delta$ as
 \be
 \int_{\NN_\delta} g(x)\, dx =  \sum_{l=1}^N
 \int_{\NN_\delta} \xi_l (\A y) \, g(y,t)\, dx(y,t) \ .
\label{dc}\ee
Let us keep this in mind particularly for estimates obtained by integration by parts relative to local variables for $y$.

\medskip
Let us consider the coordinates $X$ in \equ{cv}. Then
we have that
$$
DX (\py, t)\ =\
\left [ \pp_1 Y+  \A (t+h) \, \pp_1 \nu  + \A \pp_1 h \, \nu  \ |\   \pp_2 Y +\A(t+h)  \, \pp_2 \nu  + \A \pp_2 h\,   \nu \ |\  \nu\, \right ]\, (\A y),
$$
hence
$$
\det DX (\py, t)\, =\,
\det \left [ \pp_1 Y+  \A (t+h) \, \pp_1 \nu   \ |\   \pp_2 Y +\A(t+h)  h\, \pp_2 \nu \ |\  \nu\, \right ]\, (\A y)\ =\ $$
$$
\det \left [ \pp_1 Y   \ |\   \pp_2 Y   \ |\  \nu\, \right ]\, +\,
\A^2 (t+h)^2 \det \left [  \, \pp_1 \nu   \ |\  \pp_2 \nu  \ |\  \nu\, \right ]\,  +\,
$$
$$
\A (t+h)\, \left \{ \det \left [ \pp_1 Y  \ |\  \pp_2 \nu   \ |\  \nu\, \right ]\, +\,
 \det \left [ \pp_1 \nu  \ |\  \pp_2 Y   \ |\  \nu\, \right ]\, \,\right \}\, .$$
Since mean curvature of $M$ vanishes and the Gauss curvature equals  $|A|^2$, we obtain
$$
 dx\, =\, \left |\det DX (\py, t)\,\right |\,d\py\, dt \ =  \
 $$
 $$
 (\, 1+ \A^2 (t+h)^2\,|A|^2 ) \,  \left | \det \left [ \pp_1 Y   \ |\   \pp_2 Y   \ |\  \nu\, \right ]\,\right |\,(\A \py)\,   d\py\, dt\ =
 $$
$$
(\, 1+ \A^2 (t+h)^2\, | A (\A y) |^2\, ) \,\, dV_\A (y)\, dt \ .
$$
Using this, we estimate
$$
 I\ =\ \int_{\mathcal W}
 L(kw')\, kw'\, dx\ = \
 $$
 $$
  \int_{r_\A (y) < R}\int_{|t+h_1|< \rho_\A }  \left [\, L(kw')\, kw'\,\right ]\,
 (\, 1+ \A^2 (t+h)^2\,
 | A (\A y) |^2\, ) \,\, dV_\A (y)\, dt\ .
$$
According to Lemma    \ref{lemin10},
$$
I \ =\  \underbrace{ \int (\Delta_{M_\A } k + \A^2|A(\A y)|^2k )\,k \,\, dV_\A (y)  }_{I_1} \ +
 $$
 $$
\underbrace{  \int_{r_\A (y) < R}   \, [\,  O(\,\A\,r_\A^{-2} \log r_\A ) k \partial_{ij} k
+  O(\,\A^2\,r_\A^{-3}) k  \partial_{j} k   + O(\,\A^3\,r_\A^{-4})k^2 \,]\, \,dV_\A (y)}_{I_2}
 $$
 $$
\underbrace{ \A^2 \int_{\mathcal W} L(kw')\, kw'\,(t+h)^2\,
 | A (\A y) |^2\,\, dV_\A (y)\, dt }_{I_3}\ .
$$

 Integrating by parts, we see that
 $$
 -I_1 =  \int_{r_\A (y) < R}  \, [\, |\nn_{M_\A} k|^2 -  \A^2|A(\A y)|^2 \,]\, \,dV_\A (y)\, .
 $$
 The quantity $I_2$ involves some abuse of notation since it is expressed in local coordinates for $y$ associated to each chart, and the total should be understood in the sense \equ{dc}.
 Integrating by parts in those coordinates, we get
 $$
 I_2 =  \int_{r_\A (y) < R}  [ O(\,\A\,r_\A^{-2} \log r_\A ) \partial_{i}k \pp_jk \ + \  O(\,\A^2\,r_\A^{-3}\log{r_\A} )\pp_i k \, k \,+\, O(\,\A^3\,r_\A^{-4})k^2 \,]\, \,dV_\A (y)\, .
 $$
 Now,
  $$
  |\,\A^2\,O(r_\A^{-3}\log{r_\A})\,\pp_i k \, k\, | \,\le \,
  \, C\, [ \, \A |\nn k|^2 +  \A^3(1+ r_\A^{4})^{-1}\, k^2\, ]\, ,
 $$
 and hence  we have
  \be
 |I_2|\  \le \   C \,  \A\, \int_{r_\A (y) < R}  [ \,|\nn k|^2 \, +  \,\A^2\,(1+ r_\A^{4})^{-1}\, k^2 \,]\, \,dV_\A (y)\,
 \label{Ii2}\ee
where $C$ is independent of $\nu$ and all small $\A$. Finally, to deal with the term $I_3$, we consider the expression \equ{L(v)}
for $L(kw')$ and integrate by parts once the terms
involving second derivatives of $k$. Using that $|A|^2 = O(r_\A^{-4})$ we then get that
  \be
 |I_3|\  \le \   C\A \, \int_{r_\A (y) < R}  [\,\,|\nn k|^2 \, +  \, \A^2\,(1+ r_\A^{4})^{-1}\, k^2 \,]\, \,dV_\A (y)\, .
 \label{I3}\ee
 The same considerations above hold for $kw'$ replaced by $\eta_\delta kw'$, at this point we observe that since $\eta_\delta kw'$ satisfies Dirichlet boundary  conditions, we have
 $$
 -\int_{{\mathcal W}} L(kw'\eta_\delta)\, \eta_\delta kw'\, dx\, =\, \int_{{\mathcal W}}
  |\nn (kw'\eta_\delta)\, |^2 \, -\, \A^2 |\eta_\delta kw' \, |^2\, dx \, = \,
  $$
  $$
 \int_{ r_\A(y)< R }\left \{ \, |\nn_{M_\A} k|^2 - \A^2 |A(\A y)|^2\, k^2\, \right\}\, dV_\A  \int_\R {w'}^2\, dt
 $$
 $$
 \ +\  O\, \left ( \A \, \int_{ r_\A(y)< R } [\,\,|\nn k|^2 \, +  \, \A^2\,(1+ r_\A^{4})^{-1}\, k^2 \,]\, \,dV_\A \, \right )
 $$
 and estimate \equ{qq} has been established. This concludes the proof. \qed

\bigskip
After Lemma \ref{lemin11},  the inequality
\be
m(u_\A) \ \ge\  i(M) \label{i1-1}\ee
for small $\A$ follows at once.
Indeed,  we showed in \S \ref{eigenvaluejacobi} that the Jacobi operator has exactly $i(M)$ linearly independent
bounded eigenfunctions  $\hat{z}_i$ associated to negative eigenvalues $\la_i$ of the weighted problem in entire space $M$.
According to the theory developed in \S \ref{jacobi}, we also find that $\nn \hat{z}_i = O(r^{-2})$, hence   we may assume
\be
{\tt Q}(\hat{z}_i, \hat{z}_j)\, = \, \la_i  \int_M  q \, \hat{z}_i\,\hat{z}_j \, dV\ .
\label{qzi}\ee
Let us set $k_i(y):= \hat{z}_i(\A y)$.
According to Lemma \ref{lemin11}, setting $v_i(x) = k_i(y) \, w'(t) \eta_\delta$ and changing variables  we get
\be
\QQ(v_i,v_j) =   \A^2 {\tt Q}(\hat{z}_i, \hat{z}_j) \int_\R {w'}^2 + O(\A^3) \sum_{l=i,j} \int_M |\nn \hat{z}_l|^2 + (1+r^4)^{-1} \hat{z}_l^2 \, dV \, .
\label{qvi}\ee
From here and relations \equ{qzi}, we find that the quadratic form $\QQ$ is negative on the space spanned by the functions
$v_1,\ldots, v_{i(M)}$.  The same remains true for the functions $v_i$ smoothly truncated around $r_\A (y) =R$, for very large $R$.
We have proven then inequality \equ{i1-1}.

\bigskip
In what remains this section we will carry out the proof of the inequality
\be
m(u_\A) \ \le\  i(M). \label{i2-4}\ee

Relation \equ{qvi} suggests  that associated to a negative eigenvalue $\la_i$ of  problem \equ{M1}, there is an eigenvalue of \equ{R1}
 approximated by $\sim \la_i \A^2$. We will show next that negative eigenvalues of problem \equ{R1} cannot exceed a size $O(\A^2)$.

\begin{lemma} \label{lemin12}
There exists a $\mu>0$ independent of $R>0$ and all small $\A$ such that if $\la$ is an eigenvalue of problem $\equ{R1}$ then
$$\la \ge -\mu\,\A^2\ .$$
\end{lemma}

\proof
Let us denote
$$\QQ_\Omega (\psi,\psi)\, := \, \int_\Omega |\nn \psi|^2 - f'(u_\A)\psi^2\ . $$
Then if $\psi(x)$ is any function that vanishes for $|x'| > R\A^{-1}$ then we have
$$
\QQ(\psi,\psi) \, \ge\,  \QQ_{\NN_\delta\cap \{r_\A(y) < R\}}\,(\psi,\psi)\, +\, \gamma\int_ {\R^3 \setminus\NN_\delta }\psi^2
$$
where $\gamma>0$ is independent of $\A$ and $R$.
We want to prove that for some $\mu>0$ we have in $\Omega = \NN_\delta\cap \{r_\A(y) < R\}$ that
\be
\QQ_{\Omega }\,(\psi,\psi)\,\ge \, -\mu\A^2 \int_\Omega \frac{\psi^2 } {1 + r_\A^4}\, dx\ .
\label{need}\ee
Equivalently, let us consider the eigenvalue problem

\be
L(\psi) + \la p(\A x) \psi = 0 \quad \hbox{in } \Omega,
\label{eg}\ee
$$
\psi =0 \quad \hbox{on } r_\A = R, \quad \pp_n\psi =0 \quad \hbox{on } |t+h_1| = \rho_\A .
$$
Then we need to show that for any eigenfunction $\psi$ associated to a negative eigenfunction, inequality \equ{need} holds.
Here $\pp_n$ denotes normal derivative. Let us express this boundary operator in terms of the coordinates $(t,y)$.
Let us consider the portion of $\partial \NN_\delta$ where
\be
t+h_1(\A y) = \rho_\A(y) .
\label{bdy}\ee
We recall that for some $\gamma>0$, $\rho_\A(y) = \rho(\A y) =  \gamma\log (1 + r_\A(y))$.
Relation $\equ{bdy}$ is equivalent to
\be
z - h_0(\A y) - \rho_\A(y) = 0
\label{bdy1}\ee
where $z$ denotes the normal coordinate to $M_\A$. Then, for $\nn = \nn_x$,  we have that a normal vector to the boundary at a point satisfying \equ{bdy1} is
$$
n \ =\  \nabla z  - \nn_{M_\A}( h_0 + \rho ) = \nu (\A y)   - \A \nn_M ( h_0 + \rho) (\A y) .
$$
Now, we have that
$
\partial_t \psi   = \nabla_x \psi\cdot \nu (\A y) .
$
Hence, on points \equ{bdy}, condition $\partial_n \psi=0$ is equivalent to
\be
\pp_t \psi - \A \nn_M(h_0 + \rho) \cdot \nabla_{M_\A} \psi   = 0\, ,
\label{bc}\ee
and similarly, for
\be
t+h_1(\A y) = \rho_\A(y) .
\label{bdy2}\ee
it corresponds to
\be
\pp_t \psi - \A \nn_M(h_0 - \rho)\cdot \nabla_{M_\A} \psi   = 0\ .
\label{bc2}\ee
Let us consider a solution $\psi$ of problem \equ{eg}. We decompose
$$
\psi \ = \  \ k(y) w'(t)\eta_\delta \ + \ \ps
$$
where $\eta_\delta$ is the cut-off function \equ{etadelta} and
$$
\int_{|\tau+ h_1(\A y)| < \rho_\A(y) }  \ps(y,\tau)\, w'(\tau)\, d\tau  = 0 \foral y\in M_\A\cap\{ r_\A(y) < R\} ,
$$
 namely
 \be
 k(y)  = \frac{ \int_{|\tau+ h_1(\A y)| < \rho_\A(y) }  \psi(y,\tau)\, w'(\tau)\, d\tau}{\int_\R  w'(t)^2\eta_\delta \, dt }\ .
 \label{k(y)}\ee
Then we have
$$
\QQ_\Omega ( \psi,\psi) = \QQ_\Omega(\ps,\ps) +  \QQ_\Omega(kw'\eta_\delta ,kw'\eta_\delta)  + 2 \QQ_\Omega(kw'\eta_\delta, \ps).
$$
Since $\ps$ satisfies the same boundary conditions as $\psi$ we have that
$$
\QQ_\Omega(\ps,\ps)  = -\int_\Omega (\,\ps \Delta_x \ps + f'(u_\A) \ps^2)\, dx\  .
$$
Thus,
$$
\QQ_\Omega(\ps,\ps) \, =\,   -\int _{ r_\A < R} \int_{|t+ h_1| < \rho_\A }
[\ps\, \Delta_x \ps\, + f'(u_\A)\ps^2 ]\, (1+ \A^2 (t+h)^2 |A|^2 )\,dV_\A \, dt \,.
$$
Let us fix a smooth function $H(t)$ with $H(t) =+1$ if $t>1$, $H(t) =-1$ for $t<-1$. Let us write
$$
-\Delta_x\ps -f'(u_\A) \ps\, =\, -\partial_{tt}\ps - f'(w)\ps  + \A \pp_t\, \left [ \nn_M (h_0 + H(t) \rho)  \cdot \nn_{M_\A}\ps \right ]
$$
$$
-  \Delta_{M_\A}\ps + B(\ps).
$$
Then, integrating by parts in $t$, using the Neumann boundary condition, we get that
the integral

$$
I \ := \
$$
$$
-\int_{|t+ h_1| < \rho_\A }\left [\partial_{tt}\ps  +f'(w)\ps -\A \pp_t\, \left (\nn_M (h_0 + H(t) \rho)  \cdot \nn_{M_\A}\ps \right ) \right ]\,
 \ps\, (1+ \A^2 (t+h)^2 |A|^2 ) \, dt
$$

$$
\ =\ \int_{|t+ h_1| < \rho_\A }\left [ \partial_{t}\ps  \, - \, \A \, \left (\nn_M (h_0 + H(t) \rho)  \cdot \nn_{M_\A}\ps \right )\, \right ]\,
 \pp_t \ps\, (1+ \A^2 (t+h)^2 |A|^2 ) \, dt\
 \
 $$
 $$
-\int_{|t+ h_1| < \rho_\A }f'(w)\,
  \ps^2\, (1+ \A^2 (t+h)^2 |A|^2 ) \, dt\
 \
 $$
 $$
 +\int_{|t+ h_1| < \rho_\A }\left [ \partial_{t}\ps  \, - \, \A \, \left (\nn_M (h_0 + H(t) \rho)  \cdot \nn_{M_\A}\ps \right )\, \right ]\,
  \ps\, 2\A^2 (t+h)^2 \, |A|^2 \, dt\
$$
$$
\ =\ \int_{|t+ h_1| < \rho_\A } [\, |\pp_t\ps|^2  - f'(w)|\ps|^2 \, ]\,(1+o(1)\, ) \, +\, \A  O( r_\A^{-1} ) \nn_{M_\A}\ps\, \pp_t\ps  +  o(1)\pp_t \ps\,\ps \, dt\ .
$$
Now we need to make use of the following fact: there is a $\gamma>0$ such that if $a>0$ is a sufficiently large number, then
for any smooth function $\xi(t)$ with $\int_{-a}^a \xi\, w' \, dt = 0$ we have that
\begin{equation}
\label{ineq}
\int_{-a}^a {\xi'}^2 - f'(w) \xi^2  \, dt\,\ \ge \ \gamma \int_{-a}^a {\xi'}^2+ \xi^2  \, dt\, .
\end{equation}

Inequality (\ref{ineq}) is just a perturbation of the inequality
(\ref{q}).  We leave the details to the reader.

Hence
\be
I\ \ge \  \frac \gamma 2 \int_{|t+ h_1| < \rho_\A } [\, |\pp_t\ps|^2  + |\ps|^2 \, ]\, dt\,  + \, \int_{|t+ h_1| < \rho_\A } \A  O( r_\A^{-1} ) \nn_M\ps\, \pp_t\ps\, dt \, .
\label{I1-10}\ee
On the other hand, for the remaining part, integrating by parts in the $y$ variable the terms that involve two derivatives of $\ps$ we
get that
$$
II\, :=\, -\int_{|t+ h_1| < \rho_\A } \, dt\, \int_{ r_\A (y) < R}(\Delta_{M_\A}\ps  + B\ps )\, \ps \, (1+ \A^2 (t+h)^2 |A|^2 ) \,dV_\A(y)\,\ge
$$
\be
\int_{|t+ h_1| < \rho_\A } \, dt\, \int_{ r_\A (y) < R} |\nn_{M_\A}\ps|^2 \,+\, o(1)\, (\, \ps^2 + |\pp_t \ps| ^2  |\nn_{M_\A} \ps|^2\,)\ .
\label{I2-10}\ee
Using estimates \equ{I1-10}, \equ{I2-10}, we finally get
\be
\QQ_\Omega (\ps, \ps)\ \ge\ 3\mu \int_\Omega ( |\pp_t\ps |^2 + |\nn_{M_\A}\ps|^2 + \ps^2)\, dx\ ,
\label{ie1}\ee
for some $\mu>0$.

Now, we estimate the crossed term. We have
$$-\QQ_\Omega (\ps, kw'\eta_\delta) \,=\, \int_\Omega L(kw'\eta_\delta )\,\ps\, (1 + \A^2(t+h)^2|A|^2)\, dV_\A\, dt\,. $$
Let us consider expression \equ{L(v)} for $L(kw')$, and let us also consider the fact that
$$ L(\eta_\delta kw') = \eta_\delta L(kw') +  2\nn \eta_\delta \nn (kw') + \Delta\eta_\delta \, kw', $$ with the last two terms
producing a first order operator in $k$ with exponentially small size, at the same time with decay $O(r_\A ^{-4})$.
Thus all main contributions come from the integral
$$
I = \int_\Omega  \eta_\delta L(kw')\, \ps\,(1 + \A^2(t+h)^2|A|^2)\, dV_\A\, dt.
$$
Examining the expression \equ{L(v)},  integrating by parts once in $y$ variable those terms involving two derivatives in $k$, we see
that most of the terms obtained produce straightforwardly quantities of the type
$$
\theta := o(1) \int_{M_\A}  (|\nn k|^2 + \A^2|A|^2k^2\, ) \, dV_\A  + o(1) \int_\Omega  (|\ps|^2 + |\nn\ps|^2) .
$$
In fact
we have
$$
I \ = \ \underbrace{\int_\Omega   \Delta_{M_\A} k\, w'\,\eta_\delta \, \ps \,  dV_\A\,dt }_{I_1} +  \underbrace{
\int_\Omega  \A^2 \,a_{ij}^0\partial_ih_0\partial_jh_0 \,k\, w'''
\, \ps dV_\A\, \, dt  }_{I_2}\ +\
$$
$$
\underbrace{\int_\Omega  f''(w) \,\phi_1 \,kw'\,  \ps dV_\A\, \, dt  }_{I_3}\, + \, \theta .
$$
On the other hand, the orthogonality definition of $\ps$ essentially eliminates $I_1$. Indeed,
$$
I_1 = - \int_\Omega   \Delta_{M_\A} k\, w'\,(1-\eta_\delta) \, \ps \,  dV_\A\,dt =
\int_\Omega   \nn_{M_\A} k\, w'\,[(1-\eta_\delta) \,\nn_{M_\A} \ps  -\nn\eta_\delta \ps) \,  dV_\A\,dt = \theta \, .
$$
On the other hand, for a small, fixed number $\nu>0$ we have
$$
|I_2 |\le   C\A^2 \int_\Omega   \frac 1 {1+ r_\A^{2}}\,|k|\, |w'''|
\, |\ps| \,dV_\A\, \, dt  \,\le \, C  \nu^{-1}\A^2  \int_{M_\A}   \frac 1 {1+ r_\A^{4}}\,k^2\,dV_\A  + \nu\int_{\Omega}  |\psi|^2\, dx\ .
$$
A similar control is valid for $I_3$ since $\phi_1 = O(\A^2 r_\A^{-2})$.
We then get
\be
I\ \ge \ -  C  \nu^{-1}\A^2  \int_{M_\A}   \frac 1 {1+ r_\A^{4}}\,k^2\,dV_\A - \nu \int_\Omega |\ps |^2 .
\label{ie2}\ee
Finally, we recall that from
Lemma \ref{lemin11},
\be
\QQ_\Omega (  kw'\eta_\delta , kw'\eta_\delta )  =
 \int_{ r_\A(y)< R }\left [ \, |\nn_{M_\A} k|^2 - \A^2 |A(\A y)|^2\, k^2\, \right ]\, dV_\A  \int_\R {w'}^2\, dt \, +\, \theta .
 \label{ie3}\ee
From estimates \equ{ie1}, \equ{ie2}, \equ{ie3}, we obtain that if   $\nu$ is chosen sufficiently small, then
$$
\QQ_\Omega( \psi,\psi) \ge  - C\,\A^2  \int_{M_\A} \frac 1 {1+ r_\A^{4}}\,k^2\,dV_\A\  \ge -  \mu\,\A^2  \int_{\Omega }   \frac 1 {1+ r_\A^{4}}\,|\psi|^2\, dx\ ,
$$
for some $\mu>0$ and inequality \equ{need} follows. \qed

\bigskip

In the next result, we  show that an eigenfunction with negative eigenvalue  of problem \equ{R1} or \equ{R2}
decays exponentially, away from the interface of $u_\A$.

\begin{lemma}\label{lemin100}
let $\phi$ be a solution of either $\equ{R1}$ or $\equ{R2}$ with $\la \le 0$. Then $\phi$ satisfies in the subregion of $\NN_\A$ where it is defined that
\be
|\phi (y,t)| \, \le \, C\,\|\phi\|_\infty \, e^{-\sigma|t|}\
\label{end0}\ee
where $\sigma>0$ can be taken arbitrarily close to $\min\{\sigma_+, \sigma_-\}.$ The number $C$ depends on $\sigma$ but it is independent
of small $\A$ and large $R$. We have, moreover, that for $|\A x'| > R_0$,
\be
|\phi (x)| \le  C\sum_{j=1}^m e^{-\sigma| x_3 - \A^{-1} ( F_k(\A x') + \beta_j\A  \log|\A x'| \, )\, | } \ .
\label{end1}\ee
where $R_0$ is independent of $\A$.   Finally, we have that
\be
|\phi (x)| \le  C\, e^{- \sigma \frac \delta \A}  \quad\hbox{for } \dist(x,M_\A) > \frac \delta\A.
\label{end2}\ee
\end{lemma}

\proof
Let $\phi$ solve problem \equ{R1} for a large $R$.
Let us consider  the region between two consecutive ends $M_{j,\A}$ and $M_{j+1,\A}$.
For definiteness, we assume that this region lies inside $S_+$ so that
$
f'(u_\A)$ approaches $\sigma_+^2$ inside it.
So, let us consider the region $S$ of points $x= (x',x_3)$ such that $r_\A(x) > R_0$ for a sufficiently large but fixed $R_0>0$ and
$$
( a_j +  \A \beta_j) \log \A|x'| + b_j  + \A \gamma   \,  < \, \A x_3 \, < \,( a_{j+1} +  \A \beta_{j+1}\,) \,\log \A|x'| \,+\, b_{j+1} \, -\, \A \gamma .
$$
In terms of the coordinate $t$ near $M_{j,\A}$, saying that
 $$ \A x_3 \sim  ( a_j +  \A \beta_j) \log \A|x'| + b_j  + \A \gamma  $$ is up to lower order terms, the same as saying $t\sim \gamma$, similarly near
 $M_{j+1,\A}$. Thus given any small number $\tau>0$ we can choose $\gamma $ sufficiently large but fixed, independently of all $R_0$ sufficiently large
 and any small $\A$, such that
 $$ f'(u_\A)\, <\, - ( \sigma_+ - \tau )^2 \quad \hbox{ in } S. $$

Let us consider, for $x\in S$  and $\sigma =\sigma_+ - 2\tau$ the function
$$ v_1(x):=  e^{-\sigma [\,x_3 -    \A^{-1}( a_j +  \A \beta_j) \log \A|x'| + b_j ]}  \, +\,
 e^{-\sigma ( \A^{-1}[ a_{j+1} +  \A \beta_{j+1}) \log \A|x'| + b_{j+1} )  -  x_3 ]} \ .$$
Then $v$ has the form
$$
v_1 =     A_1 e^{-\sigma x_3 } r^{A_2}     +  B_1 e^{\sigma x_3} r^{-B_2}, \quad r=|x'|,
$$
so that
$$
\Delta v_1 =   A_2^2 r^{-2}  r^{A_2}\, A_1e^{-\sigma x_3}    +  B_2^2 r^{-2}\, B_1 r^{-B_2}e^{\sigma x_3}  +  \sigma^2 v_1 <
$$
$$
[\,  \A^2 A_2^2 R_0^{-2} +  \A^2 B_2^2 R_0^{-2}\,  + \,  \sigma^2\, ]\, v_1 \, .
$$
Here
$$ A_2  = \sigma \A^{-1}( a_j +  \A \beta_j) ,\quad  B_2 =   \sigma \A^{-1}( a_{j+1} +  \A \beta_{j+1}). $$
Hence, enlarging $R_0$ if necessary, we achieve
$$
\Delta v_1 + f'(u_\A) v_1 < 0  \quad \hbox{ in } S.
$$
Therefore $v$ so chosen is a positive supersolution  of
\be \Delta v  +f'(u_\A) v  + \la  p( \A x)v  \le  0 \quad \hbox{in } S . \label{er}\ee
Observe that the definition of $v$ also achieves that
$$
\inf_{ \partial S\setminus \{ r_\A = R_0\}} \ \ge \gamma\, >\, 0\
$$
where $\gamma$ is independent of $\A$. Now, let us observe that the function $ v_2 = e^{-\sigma (|x'| - \frac {R_0}\A) }$
also satisfies, for small $\A$, inequality \equ{er}.
As a conclusion, for $\phi$, solution of \equ{R1}, we have that
\be
|\phi(x)|\ \le \  C\,\|\phi\|_{\infty} \, [ v_1(x) + v_2(x) ]\foral x\in S, r_\A(x) < R  .
\label{poto}\ee
Using the form of this barrier, we then obtain the validity of estimate \equ{end1}, in particular that of \equ{end0},
in the subregion of $\NN_\delta$ in the positive $t$ direction of  $M_{j,\A}$  and $M_{j+1,\A}$ when $r_\A(y) > R_0$. The remaining subregions
of  $\NN_\delta \cap \{r_\A(y) > R_0\} $ are dealt with in a similar manner.
Finally, to prove the desired estimate for $r_\A (y) < r_0$ we consider  the region
where $|t|< \frac {2\delta}\A $ assuming that the local coordinates are well defined there.
In this case we use, for instance in the region
$$
\nu < t < \frac  {2\delta}\A
$$
for $\nu>0$ large and fixed, a barrier of the form
$$v(y,t) =  e^{-\sigma t} + e^{-\sigma( \frac {2\delta}\A - t )} \, .$$
It is easily seen that for small $\A$ this function indeed satisfies
$$ \Delta_x v - f'(u_\A)v < 0 $$ where $\sigma$ can be taken arbitrarily close to $\sigma_+$.
We conclude that
$$
|\phi(y,t)|\ \le \ C\|\phi\|_{\infty}  e^{-\sigma t}\quad \hbox{for } \nu < t < \frac  {\delta}\A \ .
$$
Thus estimate \equ{end0} holds true. Inequality \equ{end2} follows from maximum principle.

\medskip
Finally, for a solution of problem \equ{R2} the same procedure works, with only minor difference introduced.
Estimate \equ{poto} can be obtained after adding a growing barrier. Indeed, we obtain
$$
|\phi(x)| \, \le\,  C\,\|\phi\|_{\infty} \, [ v_1(x) + v_2(x) +\ve v_3(x) ]\foral x\in S
$$
with $v_3(x) = \ve e^{\sigma|x'|}$, and then we let $\ve\to 0$. We should also use $\ve e^{\sigma x_3}$ to deal with the region above
the last end $M_m$ and similarly below $M_1.$ We then use the controls far away to deal with the comparisons at the
second step. The proof is concluded. \qed

\medskip
\subsection{ The proof of inequality \equ{i2}}
Let us assume by contradiction  that there is a sequence $\A =\A_n \to 0$ along which
$$
m(u_\A) \ > \  i(M)=: N \, . $$
This implies that for some sequence $R_n\to +\infty$ we have that,  for all $R>R_n$,  Problem \equ{R1} has at least $N+1$ linearly independent eigenfunctions
$$
\phi_{1,\A, R},\ldots , \phi_{N+1, \A, R}
$$
associated to negative eigenvalues
$$
\la_{1,\A,R} \le \la_{2,\A,R}\le \cdots \le \la_{N+1,\A,R}\, <0\, .
$$
We may assume that $\|\phi_{i,\A, R}\|_\infty = 1$ and that
$$
\int_{\R^3}  p (\A x)\, \phi_{i,\A, R}\, \phi_{j,\A, R}\, dx \ = \ 0 \foral i,j=1,\ldots, N+1,\quad i\ne j.
$$
Let us observe that then the estimates in Lemma \ref{lemin100} imply that the contribution to the above integrals
of the region outside $\NN_\delta$ is small. We have at most
\be
\int_{\NN_\delta }  p (\A y )\, \phi_{i,\A, R}\, \phi_{j,\A, R}\, dx \ = \ O(\A^3)  \foral i,j=1,\ldots, N+1,\quad i\ne j.
\label{ortophi}\ee

From the variational characterization of the eigenvalues, we may also assume that $\la_{i,\A,R}$ defines a decreasing function of $R$.
On the other hand,
from Lemma  \ref{lemin12} we know that $\la_{i,\A, R} = O(\A^2)$, uniformly in $R$, so that we write for convenience
$$
\la_{i,\A, R}\ =\ \mu_{i,\A, R} \,\A^2 ,\quad \mu_{i,\A, R} < 0.
$$

\bigskip
 We may assume
 $\mu_{i,\A, R} \to \mu_{i,\A} <0$ as $R\to +\infty$.  We will prove that $\phi_{i,\A, R}$ converges, up to subsequences,
 uniformly over compacts to a nonzero bounded limit $\phi_{i,\A}$ which is an eigenfunction with eigenvalue $\mu_{i,\A}\A^2$ of Problem
 \equ{R2}. We will then take limits when $\A\to 0$ and find a contradiction with the fact that $\JJ$ has at $i(M)$ negative eigenvalues.

\medskip
We fix an index $i$  and consider the corresponding pair $\phi_{i,\A, R}$, $\mu_{i,\A, R}$,
 to which temporarily we  drop the subscripts ${i,\A, R}$.

\medskip
Note that by maximum principle, $|\phi|$ can have values that stay away from zero only inside
$\NN_\delta$ 
Besides, from  Lemma \ref{lemin100}, $\phi = O( e^{-\sigma|t|})$ in $\NN_\delta$.
We observe then that since $\la$ remains bounded, local elliptic estimates imply the stronger assertion
\be
|D^2\phi| + |D\phi| + |\phi| \ \le \ C\, e^{-\sigma|t|} \quad\hbox{in } \NN_\delta\ .
\label{expdec}\ee
In particular, considering its dependence in $R$, $\phi$ approaches up to subsequences, locally uniformly in $\R^3$ a limit.
We will prove by suitable estimates that that limit is nonzero. Moreover, we will show that $\phi \approx z(\A y) w'(t)$ in $\NN_\delta$ where $z$ is an eigenfunction with negative eigenvalue $\approx \mu$ 
of the Jacobi operator $\JJ$.

\medskip
First, let us localize $\phi$ inside $\NN_\delta$.
Let us  consider the cut-off function $\eta_\delta$ in \equ{etadelta},
and  the function
$$
\ttt\phi  = \eta_\delta \phi .
$$
Then $\ttt\phi$ satisfies

\be
L(\ttt \phi)\, + \, \mu \A^2 q(\A x)
 \ttt \phi \, = \,  E_\A: = -2\nabla\eta_\delta \nabla \phi - \Delta\eta_\delta \phi
\label{eee}\ee
with $L(\ttt\phi)=\Delta\ttt\phi\, +\, f'(u_\A)\ttt\phi$.
Then from \equ{expdec} we have that for some $\sigma >0$,
$$
|E_\A| \le C \A^3 e^{-\sigma|t|}(1+ r_\A^4)^{-1}.
$$
Inside $\NN_\delta$ we write in $(y,t)$ coordinates equation \equ{eee} as
\be
L_* (\ttt\phi) +  B (\ttt \phi ) + \la p (\A y)
 \ttt \phi =   E_\A
\label{ee1}\ee
where
$$
L_* (\ttt\phi) = \pp_{tt}\ttt \phi+ \Delta_{M_\A}\ttt\phi + f'(w(t))\ttt\phi\ .
$$
Extending $\ttt \phi$ and $E_\A$  as zero, we can regard equation \equ{ee1}
as  the solution of a problem in entire $M_\A \times \R$ for an operator $\LL$ that interpolates $L$ inside $\NN_\delta$
with $L_*$ outside.  More precisely $\ttt \phi$ satisfies
\be
\LL (\ttt \phi) \, :=\, L_* (\ttt\phi) +  \BB (\ttt \phi ) + \la p(\A y)
 \ttt \phi =   E_\A \quad\hbox{in } M_\A^R \times \R,
\label{ee2}\ee
where
for a function $\psi(y,t)$ we denote
\be
\BB (\psi)  \, := \,
\left\{
\begin{matrix}
\chi\, B(\psi)  &  \quad\hbox{if }  |t+h_1(\A y)| < \rho_\A(y) + 3  \\
 0 &\quad\hbox{ otherwise }
\end{matrix}
\right .  \
\label{BB-1}\ee
 and $$ \chi(y,t) = \zeta_1(y+ (t+h)\nu_\A(y))$$
  with $\zeta_1$ the cut off function defined by \equ{zetan} for $n=1$.
 In particular, $L=\LL$ in $\NN_\delta$.

\medskip
Now, we  decompose
\be
\ttt \phi(y,t)\, =\,  \vp(y,t) \,+ \, k(y) \, \eta_\delta \, w'(t)
\label{dec1}\ee
where
$$
k(y)\ = \ -\, w'(t)\, \frac{\int_\R \ttt\phi(y,\cdot)\, w'\, d\tau }{\int_\R \eta_\delta {w'}^2\, d\tau }
$$
so that
$$
\int_\R  \vp(y,t)\, w'(t)\, dt \ = \ 0 \foral y\in M_\A^R.
$$
From \equ{expdec}, $k$ is  a bounded function, of class $C^2$ defined on $M_\alpha^R$ with first and second derivatives
uniformly bounded independently of large $R$. A posteriori we expect
that $k$ has also  bounded  smoothness  as a function of  $\alpha y$, which means in particular that $Dk = O(\A)$. We will see that this is indeed the case.

\medskip
The function $\vp$ satisfies the equation
\be
\LL(  \vp ) + \mu\A^2  p (\A y )\,\vp   \ = \ - \LL(kw') + \,  E_\A \, - \, \mu\A^2  p \,k\, w'\, \quad \hbox{ in} M_\A^R \times \R\ .
\label{eqtt}\ee
We observe that the expansion \equ{L(v)}  holds true globally in $M_\A^R\times \R$  for $\LL(kw')$ replacing $L(kw')$.  We also have
the validity of expansion \equ{r7}  for the corresponding projection, namely

$$
\int_\R \, \LL(kw')\, w'\, dt\ = \
 \ (\,\Delta_{M_\A} k +  \A^2 |A|^2 k \,  )\, \int_\R{w'}^2\, dt
$$
\be
+\  O(\A r_\A^{-2})\,\partial_{ij}k \,+ \, O(\A^2 r_\A^{-3})\, \partial_i k \, + O(\A^3 r_\A^{-4})\, k
\, .
\label{rr}\ee
Thus, integrating equation \equ{eqtt} against $w'$ we find that $k$ satisfies
$$
\Delta_{M_\A} k +  \A^2 |A|^2 k\, + \,  \mu\A^2\, p(\A y)\, k
\ +
$$
$$
  O(\A r_\A^{-2})\,\partial_{ij}k \,+ \, O(\A^2 r_\A^{-3})\, \partial_i k \, + O(\A^3 r_\A^{-4})\, k \ =
$$
\be
O( \A^3 r_\A^{-4}) -     \frac 1{ \int_\R {w'}^2}  \int_\R \BB(\vp)\, w'\, dt , \quad y\in M_\A^R\ .
\label{rrr}\ee
 Let us consider the function $z(y)$ defined in $M$ by the relation $k(y)= z(\A y)$.
Then  \equ{rrr} translates in terms of $z$ as
$$
\Delta_{M} z \,+  \, |A(y)|^2 z\, + \,  \mu\, q (y)\, z    \ = \
$$
\be
\A \left [O( r^{-2})\,\partial_{ij}z \,+ \, O( r^{-3})\, \partial_i z \, + O( r^{-4})\, z  +  O( r^{-4})\right ]
  +  {\mathcal B}    \quad y\in M^R\ .
\label{rrrr}\ee
where
\be
{\mathcal B}(y) \, := \, \frac 1{ \int_\R {w'}^2} \,\A^{-2} \int_\R \BB(\vp)( \A^{-1}y,t)\, w'\, dt ,\quad y\in M^R\ .
\label{bbb}\ee
In other words we have that $k(y) = z(\A y)$,
 where $z$ solves ``a perturbation'' of the eigenvalue equation for the Jacobi operator that we treated
in \S \ref{eigenvaluejacobi}. We need to make this assertion precise, the basic element being to
prove  that the operator ${\mathcal B}[z ]$ is ``small''. For this we will derive estimates for $\vp$ from equation \equ{eqtt}.

\medskip
We shall refer to the decomposition $Q_1+\cdots +Q_6$ in \equ{L(v)} to identify different terms in $\LL(kw')$.
Let us consider the decomposition
$$
\vp = \vp_1 +\vp_2, $$
where $\vp_1$  solves the linear problem for the operator $L_*$ and the part of $\LL(kw')$ that ``does not contribute to projections'',
namely
$$
Q_3 +Q_4 =
 -w''\, \left [   \A a_{ij}^0( \partial_jh \partial_{i}k + \partial_ih \partial_{j}k) + \A^2k \Delta_M h_1
  +
  \A^2 h  a_{ij}^{1,0}( \partial_jh \partial_{i}k + \partial_ih \partial_{j}k)
    \right]
    $$
 \be
  \  +\
    \A t w'\,\left [ a_{ij}^{1,0}\partial_{ij}k  + \A b_i^{1,0}\partial_ik \right ] .
    \label{q34}\ee
 More precisely
$\vp_1$ solves  the equation

\be
L_*(\vp_1)\, +\, \A^2\mu\, p \, \vp_1    = Q_3 + Q_4 \quad \hbox{in } M_\A^R \times \R.
\label{star}\ee
This problem can indeed be solved:
according to the linear theory developed, there exists a unique solution to the problem
$$
L_*(\vp_1) + \mu\A^2  p \,\vp_1 = Q_3 + Q_4 + c(y)w'(t)\quad\hbox{in }  M_\A^R \times \R,
$$
such that
$$
\int_\R \vp_1\, w'\, dt = 0 \foral y\in M_\A^R
$$
and 
\be
\|D^2 \vp_1\|_{p,1,\sigma}
\, +\,  \|\, D \vp_1 \|_{ \infty ,1,\sigma} + \,  \|\vp_1 \|_{ \infty ,1,\sigma} \, \le \, \|Q_3+ Q_4\|_{p,1,\sigma} \, \le \,C\A \ .
\label{cotavp1}\ee
But since
$$
\int_\R (Q_3 + Q_4)\, w'\, dt \ =\ 0 \foral y\in M_\A^R
$$
it follows that actually $c(y)\equiv 0$, namely $\vp_1$ solves
equation \equ{star}.

\medskip
We claim that  $\vp_2$ has actually a smaller size than $\vp_1$. Indeed
$\vp_2$ solves the equation
\be
 L_*(\vp_2)+ \BB (\vp_2)+  \mu\A^2  p \vp_2 =  E_\A  -\BB(\vp_1)  -(Q_1 + Q_2 + Q_5 + Q_6) - \mu\A^2 q\,kw' \quad\hbox{in } M^R_\A \times \R.
\label{ecua}\ee

Now, we have that
$$
Q_1 + Q_2 + Q_5 + Q_6 \ =
$$

$$
\,  \left [ \Delta_{M_\A}k\, + \, \A h a_{ij}^{1,0}\partial_{ij}k\right ] \, w' \, + \,  \A^2\, \left[  -  |A|^2\,k\, t w'' +  \,a_{ij}^0\partial_ih_0\partial_jh_0 \,k\, w'''
  \right . \ +
$$

$$
 \left .   \A^{-2} f''(w) \,\phi_1 \,kw'\,+\,  (t+ h )^2a_{ij}^{2}\partial_{ij}k  w'  +  2(t+h)  a_{ij}^2\,  \partial_ih \partial_{j}k)\, w''
 \, \right ]\,\ +\,
 $$

 \be
 \A^3\, \left [ O( e^{-\sigma|t|} r_\A^{-2})\, \partial_{ij}k \,+ \, O( e^{-\sigma|t|} r_\A^{-3})\,\partial_{j}k\,\right ]\  =
\label{L(v)1}
\ee

$$
O(\A^2  r_\A^{-2}\, \log ^2 r_\A\, e^{-\sigma|t|}  )\, + \,  \rho(y) w'(t)\ ,
$$
for a certain function $\rho(y)$. 
On the other hand, let us recall that
$$
B\,= \,  (f'(u_\A)-f'(w)) - \A^2[ (t+h_1)|A|^2 + \Delta_M h_1]  \partial_t\ -  \A\, a_{ij}^0(\, \partial_jh \partial_{it} +
\partial_ih \partial_{jt})\ +
$$

$$
 \A (t+h) \, [ a_{ij}^1\partial_{ij}   - 2\A\, a_{ij}^1\partial_ih \partial_{jt} + \A (b_i^1\partial_i   - \A b_i^1 \partial_ih \partial_t)\, ) \, ] \ +
 $$

\be
  \A^3 (t+h)^2b^1_3  \partial_t\, +\, \A^2 [\,a_{ij}^0 + \A(t+h) a_{ij}^1 )\,]\partial_ih\partial_jh \, \partial_{tt}
\label{B}\ee
 Thus  the order of
$\BB(\vp_1) $ carries  both an extra $\A$ and an extra $r_\A^{-1}$ over those of $\vp_1$,
in the sense that
\be \|\BB(\vp_1)\|_{p,2,\sigma} \le C\A^2 . \label{pp1}\ee
From relations \equ{L(v)1} and \equ{pp1}   we find that $\vp_2$ satisfies an equation of the form
\be
 L_*(\vp_2)+ \BB (\vp_2)+  \mu\A^2  q\vp_2 =   g  + c(y)\, w' \quad\hbox{in } M^R_\A \times \R
\label{ecua1}\ee
where  for arbitrarily small $\sigma'>0$ we have
$$\| g \|_{p,2-\sigma',\sigma} \le C\A^2 . $$
Since $\vp_2$ satisfies $\int_\R \vp_2\, w'\, dt\, \equiv \, 0$, the linear theory for the operator $L_*$ yields then that
\be
\|D^2\vp_2 \|_{p,2-\sigma',\sigma}\, +\,  \|\, D \vp_2 \|_{ \infty ,2-\sigma',\sigma} + \,  \|\, \vp_2 \|_{ \infty ,2-\sigma',\sigma} \, \le \,  C\A^2\, ,
\label{bdvp2}\ee
which compared with $\equ{cotavp1}$ gives us the claimed extra smallness:
\be \|\BB(\vp_2)\|_{p,3-\sigma',\sigma} \le C\A^3 . \label{pp1-1}\ee

\medskip
Let us decompose in \equ{bbb}
$$
{\mathcal B} = {\mathcal B}_1 + {\mathcal B}_2
$$
where
\be
{\mathcal B}_l := \, \frac 1{ \int_\R {w'}^2} \,\A^{-2} \int_\R \BB(\vp_l)( \A^{-1}y,t)\, w'\, dt ,\quad l=1,2.
\label{bbb1}\ee
From Lemma \ref{lemin4} 
we get that
\be
\|{\mathcal B}_1\|_{p,  2 -\frac 2p - \sigma'} \le C\A^{-2 } \|  \BB(\vp_1) \|_{p, 2 ,\sigma} \le C
\label{bg-1}\ee
and
\be
\|{\mathcal B}_2\|_{p,  3 -\frac 2p -2\sigma'} \le C\A^{-2 } \|  \BB(\vp_2) \|_{p, 3 -\sigma',\sigma} \le C\A
\label{bg2}\ee

Now, we apply the estimate in part (b) of Lemma \ref{lemin3} to equation \equ{rrrr} and then get for $z(y) = k(\frac y\A)$ the estimate
\be
\|  D^2z\|_{p, 2-\frac 2p - 2\sigma'  } + \| \, (1+ |x|)^{1-2\sigma'} \,Dz\, \|_{\infty} \le
C \, [ \, \|f\|_{p, 2-\frac 2p - 2\sigma' }  \, + \, \|z\|_\infty \  ] \,
\label{eeq}\ee
where
$$
f= \A \left [O( r^{-2})\,\partial_{ij}z \,+ \, O( r^{-3})\, \partial_i z \, + O( r^{-4})\, z  +  O( r^{-4})\right ]
  +  {\mathcal B}   .
$$
Then from estimate \equ{eeq} it follows that for small $\A$,
\be
\|  D^2z\|_{p, 2-\frac 2p - 2\sigma'  } + \| \, (1+ |x|)^{1-2\sigma'} \,Dz\, \|_{\infty} \le C\A.
\label{35}\ee
Using this new information, let us go back to equation \equ{star} and to the expression \equ{q34}  for $Q_3 + Q_4$.
The terms contributing the largest sizes in this function can be bounded by

$$
 C\,\A \, e^{-\sigma|t|}\, \left[ \, \frac {|D k|} {1 + r_\A}  + \frac {|D^2 k|} {1 + r_\A^2}
\, \right ]\ .
$$

Now, we compute
$$
 (1+r_\A(y)^2 )^p \int_{B(y, 1)}    \frac {|D^2 k|^p } {(1 + r_\A^2)^p} \, dV_\A\  \le
 $$
 $$
 C \A^{2p -2}
 \int_{B(y, \A)}|D^2 z|^p \, dV \ \le \ C \A^{2p -2} \|  D^2z\|_{p, 2-\frac 2p - 2\sigma'  }\ \le\  C\A^{2p -2}  \, ,
$$
and
$$
 (1+r_\A(y)^{2-2\sigma'} )^p \int_{B(y, 1)}    \frac {|D k|^p } {(1 + r_\A)^p } \, dV_\A\  \le
 $$
 $$
 C\|\, |D k|\,(1 + r_\A)^{1-2\sigma'} \,\|^p_\infty \, =\,  C\, \A^p  \, \|\, |D z|\,(1 + r)^{1-2\sigma'} \,\|^p_\infty \, \le\,  C\A^{p}\, .
 $$
As a conclusion, from expression \equ{q34} we obtain  that
$$
\| Q_3+ Q_4\|_{p, 2- 2\sigma', \sigma} \le C\A^{2} ,
$$
and therefore a substantial reduction of the size of $\vp_1$, compared with \equ{cotavp1}, we  have
\be
\|D^2 \vp_1\|_{p, 2-2\sigma' ,\sigma}
\, +\,  \|\, D \vp_1 \|_{ \infty ,2-2\sigma' ,\sigma}\, \le \, C\A^{2} \ ,
\label{cotavp2}\ee
hence, using again Lemma \ref{lemin4} we get
\be
\|{\mathcal B}_1\|_{p,  3 -\frac 2p - 3\sigma'} \le C\A^{-2} \|  \BB(\vp_1) \|_{p, 3- 2\sigma' ,\sigma} \le C\A^{}
\label{bg}\ee
which matches the size we initially found for ${\mathcal B}_2$ in \equ{bg2}.

\bigskip
We recall that $\phi =\phi_{i,\A,R}$ has a uniform $C^1$ bound  \equ{expdec}.
Thus, passing to a subsequence if necessary, we may assume that
 $$ \phi_{i,\A,R}\to \phi_{i,\A}\quad\hbox{as } R\to +\infty, $$ locally uniformly, where  $\phi_{i,\A}$ is bounded and  solves
 \be
 \Delta \phi_{i,\A} + f'(u_\A) \phi_{i,\A} + \mu_{i,\A}\,\A^2 \,p(\A x)\, \phi_{i,\A} \, =\, 0\quad\hbox{in } \R^3.
 \label{P11-2}\ee
Let us return to equation \equ{rrrr} including the omitted subscripts.
Thus  $k= k_{i,\A,R}$ satisfies the local uniform convergence in $M^\A$, 
$$
k_{i,\A,R}(y) =  c\int_{ |t+h_1| < \rho_\A } \phi_{i,\A, R} \,w'\, dt \to   c\int_{ |t+h_1| < \rho_\A } \phi_{i,\A} \,w'\, dt\, =:\,
k_{i,\A}(y)\ .
$$
We have that $z = z_{i,\A,R}$ satisfies
$$
\Delta_{M} z_{i,\A,R} \,+  \, |A(y)|^2 z_{i,\A,R}\, + \,  \mu_{i,\A,R}\, q(y)\, z_{i,\A,R}    \ = \
$$
$$
\A \left [O( r^{-2})\,\partial_{ij}z_{i,\A,R} \,+ \, O( r^{-3})\, \partial_i z_{i,\A,R} \, + O( r^{-4})\, z_{i,\A,R}  +  O( r^{-4})\right ]
\  +
$$
$$
{\mathcal B}_{i,\A,R},    \quad y\in M^R\ ,
$$
where
\be
\|{\mathcal B}_{i,\A, R} \|_{p,  3 -\frac 2p - 3\sigma'} \ \le\ C\A^{ }
\label{bg3-1}\ee
with arbitrarily small $\sigma'>0$ and $C$ independent of $R$.
We apply now the estimates in Lemma \ref{lemin3} for some $1<p < 2$ and find that for $C$ independent of $R$ we have
$$
\|z_{i,\A,R}\|_{L^\infty (M^R)} \le  C \, [ \, \|z_{i,\A,R}\|_{L^\infty ( r < R_0)} + O(\A)\, ] \
$$
or equivalently
\be
\|k_{i,\A,R}\|_{L^\infty (M_\A^R)} \le  C \, [ \, \|k_{i,\A,R}\|_{L^\infty ( r_\A < R_0)} + O(\A)\, ] \ .
\label{ppp}\ee
Since from \equ{dec1} we have that
\be
\phi_{i,R,\A} (y,t)\, =\,  \vp_{i,R,\A}(y,t) \,+ \, k_{i,R,\A}(y) \,  \, w'(t)\quad \hbox{in } \NN_\delta,\quad r_\A(y) \le R,
\ee
where we have uniformly in $R$
$$
|\vp_{i,R,\A}(y,t)|  = O (\A\, e^{-\sigma|t|} r_\A^{-2}),
$$
while
$\phi_{i,R,\A} = O( e^{-\frac a\A})$ outside $\NN_\delta$, and
$$
\|\phi_{i,R,\A}\|_\infty  = 1,
$$
then
$$
\|k_{i,\A,R}\|_{L^\infty (M_\A^R)}\ \ge \ \gamma>0
$$
uniformly in $R$. Thus from \equ{ppp}, the  limit  $k_{i,\A}$ as $R\to +\infty$ cannot be zero.
We have thus found that  $\phi_{i,\A}$ is non-zero.
Moreover, we observe the following: Since the functions
$$Z_{i} := \partial_iu_\A,\quad i=1,2,3,\quad Z_{4} := -x_2\pp_1 u_\A + x_1\pp_2 u_\A $$
are bounded solutions of \equ{R2} for $\la =0$,
we necessarily have that
\be
\int_{\R^3} p(\A x)\, Z_{j} \phi_{i,\A}\, dx \ =\ 0 ,\quad j=1,2,3,4.
\label{ortho1}\ee
Let
\[ \hat{Z}_i= \sum_{l=1}^4 d_{il} Z_{l}, i=1, \ldots, J \, .\]
Then we also have
\be
\int_{\R^3} p(\A x)\, \hat{Z}_{i} \phi_{i,\A}\, dx \ =\ 0 ,\quad i=1,\dots, J \, .
\label{ortho}
\ee
Now we want to let $\A \to 0$. $z_{i,\A}$ satisfies
$$
\Delta_{M} z_{i,\A} \,+  \, |A(y)|^2 z_{i,\A}\, + \,  \mu_{i,\A}\, p(y)\, z_{i,\A}    \ = \
$$
$$
\A \left [O( r^{-2})\,\partial_{ij}z_{i,\A} \,+ \, O( r^{-3})\, \partial_i z_{i,\A} \, + O( r^{-4})\, z_{i,\A}  +  O( r^{-4})\right ]
\,  +
{\mathcal B}_{i,\A},    \quad y\in M\ ,
$$
with
\be
\|{\mathcal B}_{i,\A} \|_{p,  3 -\frac 2p - 3\sigma'} \ \le\ C\A^{ }.
\label{bg3}\ee
Moreover,
$$
\|z_{i,\A}\|_{L^\infty (M)} \le  C \, [ \, \|z_{i,\A,R}\|_{L^\infty ( r < R_0)} + O(\A)\, ] \ .
$$
Since we also have that
$$
\|D^2 z_{i,\A}\|_{L^p (M)} \le  C \, \ ,
$$
Sobolev's embedding implies that passing to a subsequence in $\A$, $z_{i,\A}$ converges as $\A\to 0$, uniformly over compact subsets of $M$ to
a non-zero bounded solution $\bar{z}_i$ of the equation
$$
\Delta_{M} \bar{z}_{i} \,+  \, |A(y)|^2 \bar{z}_{i}\, + \,  \mu_{i}\, q(y)\, \bar{z}_{i}    \ = \ 0\quad\hbox{in } M,
$$
with $\mu_i\le 0$.

\medskip
Now, we have that
$$
\phi_{i,\A} = z_{i,\A} (\A y ) \, w'(t)\, + \vp_i(y,t) \quad \hbox{in } \NN_\delta
$$
where
$$
|\vp_i (y,t)| \le C\A e^{-\sigma|t|} \, .
$$
We recall that
$$
\int_{\NN_\A}  q(\A y)\phi_{i,\A}\, \phi_{j,\A} \, dx \ =  O( \A) \foral i\ne j\, .
$$
 Since on $\NN_\delta$,
 $$dx =   ( 1+ \A^2 |A|^2(t+h) )dV_\A \, dt\, ,$$
we get then that
$$
\int_{M_\A}  q(\A y) z_{i,\A} (\A y ) \, z_{j,\A} (\A y ) \, dV_\A  =  O( \A)
$$
or
$$
\int_{M}  q( y) z_{i,\A} (y) \, z_{j,\A} ( y ) \, dV   =  O( \A^3)\,  \foral i\ne j.
$$
We conclude, passing to the limit,  that the $z_i$'s $i=1,\ldots, N+1$ satisfy
$$
\int_M q\, \bar{z}_{i} \bar{z}_{j} \, dV \, =\, 0 \foral i\ne j .
$$
Since, as we have seen in \S \ref{eigenvaluejacobi}, this problem has exactly $N=i(M)$ negative eigenvalues, it follows that
$\mu_{N+1}=0$, so that that $z_{N+1}$ is a  bounded Jacobi field.

\medskip
But we recall that, also
$$
Z_{i} = z_{i} (\A y) w'(t) + O(\A e^{-\sigma |t|} )\foral i=1,\ldots, J,
$$
hence the orthogonality relations $\equ{ortho}$ pass to the limit to yield
$$
\int_M q\, \hat{z}_{i}\cdot \bar{z}_{N+1} \, dV \, =\, 0, \quad i=1,\ldots, J\, .
$$
where $\hat{z}_{i}$'s are the $J$ linearly independent Jacobi fields. We have thus reached a contradiction with the non-degeneracy assumption
for $M$ and the proof of $m(u_\A) = i(M)$ is concluded.

\medskip
Finally, the proof of the non-degeneracy of $u_\A$ for all small $\A$ goes along the same lines. Indeed, the above arguments are also valid for
a bounded eigenfunction in entire space, in particular for $\mu=0$.
If we assume that a bounded solution $Z_{5}$ of equation \equ{R2}
is present, linearly independent from  $Z_{1},\ldots, Z_{4}$, then we assume
that \be \int_{\R^3} p(\A x) \, Z_{5}\, \hat{Z}_i\, dx \ =\ 0  \quad i=1,\dots, J. \label{yy}\ee
Thus,
if  in the same way as before, we have that in $\NN_\delta$,
$$
Z_5 = z_5 (\A y) w'(t) + \vp
$$
with $\vp$ orthogonal to $w'(t)$ for all $y$ and  $\vp$  small with size $\A$ and uniform exponential decay in $t$. The function $z_{\A}$
solves an equation of the form \equ{bg3}, now for $\mu=0$. In the same way as we did before, it converges uniformly on compacts to a non-zero
limit which is a bounded Jacobi field. But the orthogonality  \equ{yy} passes to the limit, thus implying the existence of
at least $J+1$ linearly independent Jacobi field. We have  reached a contradiction that finishes the proof of Theorem \ref{teo2}.
\qed

\bigskip
\setcounter{equation}{0}
\section{\emph{Further comments and open questions}}\label{22}

\subsection{Symmetries} As it is natural, the invariances of the surface are at the same time inherited from the construction. If $M$ is a catenoid, revolved around the
$x_3$ axis, the solution in Theorem \ref{teo1} is radial in the first two variables,
$$u_\A (x) = u_\A \left ( \, |x'| ,\, x_3\right ).$$ This is a consequence of the construction. The invariance of the Laplacian under
rotations and the autonomous character of the nonlinearity imply that the entire proof can be carried out in spaces of functions with this radial symmetry. More generally, if $M$ is invariant a group of linear isometries, so will be the solution found, at least
in the case that $f(u)$ is odd. This assumption allows for odd reflections. The Costa-Hoffmann-Meeks surface is invariant under a discrete group
constituted of combination of dihedral symmetries and reflections to which this remark apply.

\bigskip

\subsection{Towards a classification of finite Morse index solutions}

 \ \\ Understanding
 bounded, entire solutions of nonlinear elliptic equations in $\R^N$ is a problem that has always been at the center of  PDE research.
This is the context of  various  classical results in PDE literature like the Gidas-Ni-Nirenberg theorems  on radial symmetry of one-signed solutions, Liouville type theorems, or the achievements around De Giorgi conjecture. In  those results, the geometry of level sets of the solutions turns out to be a posteriori very simple (planes or spheres). More challenging seems the problem of classifying solutions with finite Morse index, in a model as simple as the Allen-Cahn equation. While the solutions predicted by Theorem \ref{teo1} are generated in an asymptotic setting, it seems plausible  that they contain germs of generality, in view of parallel facts in the theory of minimal surfaces. In particular we believe that the following two statements hold true for a  a bounded solution  $u$ to equation $\equ{ac}$ in $\R^3$.

\medskip
(1)  {\em  If  $u$ has finite Morse index and
$\nabla u(x)\neq  0$
outside a bounded set, then
  each level set of $u$ must have outside a large ball a finite number of components, each of them  asymptotic  to either
 a plane or to a catenoid. After a rotation of the coordinate system, all these components  are graphs of  functions of the same two variables.  }

\medskip
 (2) {\em If $u$  has Morse index equal to one. Then $u$  must be axially symmetric, namely after a rotation and a translation, $u$ is radially symmetric in two of its variables. Its level sets have two ends, both of them  catenoidal.
}

 \medskip
 It is worth mentioning that a {\em balancing formula} for the ``ends'' of level sets to the Allen-Cahn equation is available in $\R^2$, see \cite{gui-cross}. An extension of such a formula to $\R^3$ should involve the  configuration (1) as its basis.
The condition
 of finite Morse index can probably be replaced by the energy growth $\equ{dirichlet}$.

\medskip
 On the other hand, (1) should not hold if the condition $\nn u\ne 0$ outside a large ball is violated.
 For instance, let us consider the octant $\{x_1,x_2,x_3\ge 0\}$ and the odd nonlinearity $f(u) =(1-u^2)u$.  Problem \equ{ac} in the octant
  with zero boundary data can be solved by a super-subsolution scheme (similar to that in \cite{fife}) yielding a positive solution. Extending by successive odd reflections to the remaining octants, one generates an entire solution (likely to have finite Morse index), whose zero level set does not have the characteristics above: the condition $\nn u\ne 0$  far away corresponds to {\em embeddedness of the ends.}

\medskip
Various rather general conditions on a minimal surface imply that it
is a catenoid. For example, R. Schoen \cite{schoen 1}  proved
that a complete embedded minimal surface in $\R^3$ with two ends
must be catenoid (and hence it has  index one).  One may wonder if a
bounded solution to \equ{ac} whose zero level set has only  two ends
is radially symmetric in two variables. On the other hand a one-end
minimal surface is forced to be a plane \cite{hm22}. We may wonder
whether or not  the zero level set lies on a half space implies
that the solution depends on only one variable. 


\medskip
These questions seem rather natural generalizations of  that by De Giorgi, now
on the classification finite Morse index entire
solutions of \equ{ac}. The case in which the minimal surfaces have
finite topology but infinite total curvature, like the helicoid, are
natural objects to be considered. While  results parallel to that in
Theorem \ref{teo1} may be expected possible, they may have rather
different nature. The condition of diverging ends in $\beta$ is not
just technical. If it fails a solution may still be associated to
the manifold but interactions between neighboring interfaces, which
are inherent to the Allen-Cahn equation but not to the minimal
surface problem, will come into play.   The case of infinite
topology may also give rise to very complicated patterns, we refer
to Pacard and Hauswirth  \cite{Hau-Pacard} and references therein
for recent result on construction of minimal surfaces in this
scenario.

\bigskip
\setcounter{equation}{0}
\section{ \emph{Appendix}}

In this appendix we carry out the computations that lead to Lemma \ref{lapfinal}.

\subsection{Coordinates near $M$ and the Euclidean Laplacian }
Let us consider the smooth map

\be
(y,z)\in M\times \R  \ \longmapsto \  x= \ttt X(y,z) = y + z \nu(y) \in \R^3.
\label{x}\ee
Let $\OO$ be a set as in the statement of Lemma \ref{lapfinal}, and consider the subset of $M\times \R$ defined as
$$
\ttt \OO = \{ \, (\A y,\A (t+h(y))\, ) \in M\times \R \ / \  (t,y)\in \OO  \, \}.
$$
Then $\ttt X|_{\ttt \OO}$ is one to one, and
$$
\ttt \OO \subset \{ (y,z)\in M\times \R \ /\ |z| < \delta\, \log (1+ r(y))\}.
$$
Since along ends $\partial_i \nu = O(r^{-2})$ so that $ z\partial_i \nu$ is uniformly small in $\ttt O$, it follows that $\ttt  X$
is actually a diffeomorphism onto is image, $\ttt \NN = \ttt X( \ttt \OO) =\A \NN$.

\medskip
The Euclidean Laplacian $\Delta_x$ can be computed in such a region by the well-known formula in terms of the coordinates $(y,z)\in \ttt \OO$
as

\be
\Delta_x =   \partial_{zz}+ \Delta_{M_z}  - H_{M_z}\partial_z, \quad x = \ttt X(y,z),\quad (y,z)\in \OO
\label{laplacian}\ee
where $M_z$ is the manifold
$$
M_z = \{ y + z\nu (y) \ /\ y\in M \}.
$$
Local coordinates $y= Y_k(\py)$, $\py\in \R^2$ as  in (\ref{Yk})
induce natural local coordinates in $M_z$.
The metric $g_{ij}(z)$ in $M_z$ can then be computed as
\be
g_{ij} (z) =  \left < \pp_i Y, \pp_j Y \right>   +  z(\left < \pp_i Y , \pp_j \nu \right> + \left < \pp_j Y, \pp_i \nu \right> ) +  z^2 \left < \pp_i \nu , \pp_j \nu \right>
\label{z}\ee
or
$$
g_{ij}(z) =   g_{ij} +   z\,O(r^{-2})\,   +  z^2O( r^{-4})\, .
$$
where these relations can be differentiated.
Thus we find from the expression of $\Delta_{M_z}$ in local coordinates that
\be
\Delta_{M_z} = \Delta_M +  z a_{ij}^1(y, z)\partial_{ij}  + zb_i^1(y,z)\partial_i , \quad y= Y(\py)
\label{DeltaMz}\ee
where $a_{ij}^1, b_i^1$ are smooth functions of their arguments.
Let us examine this expansion closer  around the ends of $M_k$  where $y= Y_k(\py)$ is chosen as in \equ{Yk}.
In this case, from \equ{z} and \equ{gij} we find
 $$
g^{ij}(z) =   g^{ij} +  z\, O( r^{-2} ) + z^2O(r^4) + \ldots
$$
Then we find that for large
$r$,
\be
\Delta_{M_z}  =  \Delta_M + z\, O(r^{-2}) \partial_{ij}   +  z O(r^{-3})\partial_{i}.
\label{expMz}\ee

\medskip
Let us consider the remaining term in the expression for the Laplacian, the mean curvature $H_{M_z}$.
We have the validity of the formula

$$
H_{M_z} =  \sum_{i=1}^2 \frac {k_i}{1- k_i z }   = \sum_{i=1}^2 k_i +  k_i^2 z +  k_i^3 z^2 + \cdots
$$
where $k_i$, $i=1,2$ are the principal curvatures.
Since $M$ is a minimal surface, we have that $k_1 + k_2 = 0$.  Thus
$$
|A|^2 =  k_1^2 + k_2^2 = -2k_1k_2 = -2K
$$
where $|A|$ is the Euclidean norm of the second fundamental form, and $K$ the Gauss curvature.
As $r\to +\infty$ we have seen that
$k_i= O(r^{-2})$ and hence
$
|A|^2=   O ( r^{-4} ).
$
More precisely, we find for large $r$,
$$
H_{M_z} =  |A|^2 z +   {z^2}O({r^{-6}}).
$$


\medskip

Thus we have found the following expansion for the Euclidean Laplacian,
\be
\Delta_x  =  \partial_{zz} + \Delta_M   - z|A|^2\partial_z  + B
\label{euclidean}\ee
where expressed in local coordinates in $M$ the operator $B$ has the form
\be
B = z\, a_{ij}^1(y,z)\partial_{ij} + z\, b_i^1(y,z)\partial_i   + z^2 b_3^1(y,z)\partial_z
\ee
with $a_{ij}^1,\ b_i^1,\  b_3^1$ smooth functions. Besides,  we find  that
\be
a_{ij}^1(y,z) = O(r^{-2}), \quad  b_i^1(y,z)= O(r^{-3}),\quad   b_i^1(y,z) =O(r^{-6}) ,
\label{growth}\ee
uniformly in $z$ for $(y,z) \in \ttt \OO.$
Moreover, the way these coefficients are produced from the metric yields for instance that
$$
a_{ij}^1(y,z) = a_{i,j}^1(y,0) + za_{i,j}^{(2)}(y,z), \quad a_{i,j}^{2}(y,z) = O(r^{-3}),
$$
$$
b_i^1(y,z) = b_i^1(y,0) + z b_i^{(2)}(y,z), \quad b_i^{(2)}(y,z) = O(r^{-4})\, .
$$



We summarize the discussion above. Let us consider the parameterization  in \equ{x} of the region $\ttt{\mathcal N}$.

\begin{lemma} \label{lem2.1}
The Euclidean Laplacian  can be expanded in $\ttt {\mathcal N}$ as 
$$
\Delta_x =  \partial_{zz}+ \Delta_{M_z}  - H_{M_z}\partial_z \ =
 $$
 $$
 \partial_{zz} + \Delta_M -z\, |A|^2\partial_z\, +\, z\, [ a_{ij}^1(\py , z)\partial_{ij} + b_i^1(\py , z) \partial_i ]\,
 + z^2b_3^1(\py , z)\partial_z ,
$$
$$
\Delta_M  =  a_{ij}^0\partial_{ij} + b_i^0 \partial_i,\quad   x=\ttt X(y,z),\quad  (y,z)\in \ttt \OO,
$$
where $a_{ij}^l$, $b_j^l$ are smooth, bounded functions, with the index $k$ omitted.
In addition, for $k=1,\ldots, m$,
$$
a_{ij}^l = \delta_{ij}\delta_{0l}  + O(r^{-2}),\quad  b_i^l = O(r^{-3}), \quad  b_3^1 = O(r^{-6}) \ ,
$$
as $r= |\py|\to \infty$, uniformly in $z$ variable.
\end{lemma}

\bigskip
\subsection{Laplacian in expanded variables}
Now we consider the expanded minimal surface
$
M_\A = \A^{-1} M
$
for a small number $\A$. We have that
${\mathcal N} = \A^{-1}\ttt {\mathcal N}$.
We describe  ${\mathcal N}$  via  the coordinates
\be
x= X(y,z) : =   y + z\nu_\A (y),\quad (y,z) \in \A^{-1}\ttt \OO.
\label{x1}\ee
Let us observe that
$$ X(y, z) = \A^{-1} \ttt X (\A y, \A z) $$
where $ \ttt x = \ttt X(\ttt y, \ttt z) = \ttt y + \ttt z \nu(\ttt y) , $
where the coordinates in ${\mathcal N}_\delta$ previously dealt with.
We want to compute the Euclidean Laplacian in these coordinates associated to $M_\A$.
Observe that
$$
\Delta_x [u(x)]\,|_{x= X(y,z)}  = \A^2 \Delta_{\ttt x} [ u(\A^{-1}\ttt x)\, ]\,|_{\ttt x= \ttt X( \A y,\A z) }
$$
and that the term in the right hand side is the one we have already computed.
In fact  setting
$ v(y,z) := u( y+ z\nu_\A( y)  )$,
we get

\begin{equation}
\Delta_x u\,|_{x= X(y,z)}\, = \, \A^2 (\Delta_{\ttt y, M_{\ttt z}} +
\partial_{\ttt z \ttt z}  - H_{M_{\ttt z}} \partial_{\ttt z} ) \, [
v(\A^{-1}\ttt y , \A^{-1} \ttt z) ] \,\left . \right|_{ (\ttt y
,\ttt z)= (\A y, \A z)} \, . \label{form}
\end{equation}

\bigskip
We can then use the discussion summarized in Lemma \ref{lem2.1} to obtain a representation of $\Delta_x$ in ${\mathcal N}$ via the coordinates
$ X(y,t)$ in \equ{x1}. Let us consider the local coordinates  $Y_{k\A}$ of $M_\A$  in \equ{chart2}.
.

\begin{lemma} \label{lem2.2}
In ${\mathcal N}$  we have  
$$
\Delta_x =  \partial_{zz}+ \Delta_{M_{\A,z}}  - H_{M_{\A,z}}\partial_z \ =
 $$
 $$
 \partial_{zz} + \Delta_{M_\A}  - \A^2z\, |A(\A y)|^2\partial_z\, +\, \A z\, [ a_{ij}^1(\A\py ,\A z)\partial_{ij} + \A b_i^1(\A\py ,\A z) \partial_i ]\,
 + \A^3 z^2b_3^1(\A\py ,\A z)\partial_z ,
$$

$$
\Delta_{M_\A}  =  a_{ij}^0(\A y)\partial_{ij} + b_i^1(\A y) \partial_i,\quad  (y,z)\in  \A^{-1}\ttt \OO,\quad y=Y_{k\A}(\py)
$$
where $a_{ij}^l$, $b_j^l$ are smooth, bounded functions.
In addition, for $k=1,\ldots, m$,
$$
a_{ij}^l = \delta_{ij}\delta_{0l}  + O(r_\A^{-2}),\quad  b_i^l = O(r_\A^{-3}), \quad  b_3^1 = O(r_\A^{-6}) \ ,
$$
as $ r_\A (y) = |\A \py|\to \infty$, uniformly in $z$ variable.
\end{lemma}


\subsection{The proof of Lemma \ref{lapfinal}}\label{fin}

Let us consider
a function $u$ defined in   ${\mathcal N}$, expressed in coordinates $x= X(y,z)$, and consider the expression of $u$ in the coordinates
$ x= X_h(y,t)$, namely the function $v(\py,t)$ defined by the relation in local coordinates $y=Y_k(\py)$,
$$
v(\py, z- h(\A \py) ) = u(\py ,z),
$$
(by slight abuse of notation we are denoting just by $h $ the function $h\circ Y_k$).
Then we compute
$$
\partial_i u =  \partial_i v - \A \partial_t v \partial_ih, \quad \partial_z u = \partial_t v,
$$
$$
\partial_{ij} u   =
\partial_{ij} v -
 \A \partial_{it}v\partial_jh
- \A \partial_{jt}v\partial_ih +
\A^2 \pp_{tt}v\partial_ih\partial_jh
- \A^2 \partial_t v \pp_{ij}h \, .
$$

Observe that, in the notation for coefficients in Lemma \ref{lem2.2},
$$ a_{ij}^0 \partial_{ij} h  +  b_i^0 \partial_ih = \Delta_M h ,\quad   a_{ij}^0 \partial_{ij} v  +  \A b_i^0 \partial_iv = \Delta_{M_\A} v\ .$$
We find then
$$
\Delta_x   =  \, \partial_{tt}  + \Delta_{M_\A}  - \A^2[ (t+h)|A|^2 + \Delta_M h]  \partial_t\ -
2 \A \, a_{ij}^0\, \partial_jh \partial_{it}  \ +
$$
$$
 \A (t+h) \, [ a_{ij}^1\partial_{ij}   -2 \A\, a_{ij}^1\,\partial_ih \partial_{jt} \, +\, \A (b_i^1\partial_i   - \A b_i^1 \partial_ih \partial_t)\,  \, ] \ +
 $$

\be
  \A^3 (t+h)^2b^1_3  \partial_t\, +\, \A^2 [\,a_{ij}^0 + \A(t+h) a_{ij}^1 \,]\partial_ih\partial_jh \, \partial_{tt}
\label{laplacianfinal1}\ee

\noindent
where all the coefficients are understood to be evaluated at $\A\py$ or $(\A \py , \A(t+h(\A \py) )$.
 The desired properties of the coefficients have already been established. The proof of Lemma \ref{lapfinal} is concluded.
 \qed

\bigskip
{\bf Acknowledgments:}  We are indebted to N. Dancer for a useful discussion, and to F. Pacard for pointing us out reference \cite{traizet}.
The first author has been partly supported
by research grants Fondecyt 1070389 and FONDAP, Chile.  The second
author has been supported by  Fondecyt grant 1050311, and FONDAP, Chile. He thanks the Department of Mathematics of Politecnico di Torino for its hospitality while part of the work was completed.  The research of the third author is partially supported by a General Research Fund  from RGC of Hong Kong and  Focused Research Scheme
 from CUHK.


\begin{thebibliography}{aaaa}




\bibitem{cabre} L. Ambrosio and X. Cabr\'e, {\em Entire solutions of semilinear elliptic equations in $\R^3$ and a conjecture of De Giorgi,} Journal Amer. Math. Soc. 13 (2000), 725--739.








\bibitem{caffarelli1} L. Caffarelli, A. C\'ordoba, {\em Uniform
convergence of a singular perturbation problem,} Comm. Pure Appl.
Math. XLVII (1995), 1--12.

\bibitem{caffarelli2} L. Caffarelli, A. C\'ordoba, {\em
 Phase transitions: uniform regularity of the intermediate layers.} J. Reine Angew. Math. 593 (2006), 209--235.

\bibitem{Costa1} C.J. Costa, {\em Imersoes minimas en $\R^3$ de genero un e curvatura total finita. PhD thesis}, IMPA, Rio de Janeiro, Brasil (1982).

\bibitem{Costa2} C.J. Costa, {\em Example of a complete minimal immersions in $\R^3$ of genus one and three embedded ends}, {Bol. Soc. Bras. Mat. } 15(1-2)(1984), 47--54.



\bibitem{cabreterra}
X. Cabr\'e, J. Terra {\em
Saddle-shaped solutions of bistable diffusion equations in all of $\R^{2m}$. } Preprint 2008.


\bibitem{dancer}
 E. N. Dancer, {\em Stable and finite Morse index solutions on $\bold R\sp n$} or on bounded domains with small diffusion. Trans. Amer. Math. Soc. 357 (2005), no. 3, 1225--1243



\bibitem{fife}
H. Dang, P.C. Fife, L.A.  Peletier,  {\em Saddle solutions of the bistable diffusion equation.} Z. Angew. Math. Phys. 43 (1992), no. 6, 984--998.



\bibitem{dg}
E. De Giorgi, {\em
 Convergence problems for functionals and operators,} Proc. Int. Meeting
on Recent Methods in Nonlinear Analysis (Rome, 1978), 131--188, Pitagora, Bologna(1979).







\bibitem{dkwdg}
M. del Pino, M. Kowalczyk,  J. Wei, {\em On De Giorgi's Conjecture
in Dimensions $N \geq 9$}, preprint 2008, arXiv.org/0806.3141.

\bibitem{dkwdg-n}
M. del Pino, M. Kowalczyk, J. Wei, {\em A counterexample to a
conjecture by De Giorgi in large dimensions}, Comp. Rend.
Mathematique 346 (2008), 23-24, 1261-1266.



\bibitem{FV} A. Farina and E. Valdinoci, {\em The state of art for a conjecture of De Giorgi and related questions}, to appear in ``Reaction-Diffusion Systems and Viscosity Solutions'', World Scientific, 2008.






\bibitem{gg}
N. Ghoussoub and C. Gui, {\em On a conjecture of De Giorgi and some related problems,}
Math. Ann. 311 (1998), 481--491.

\bibitem{GT}  D. Gilbarg and N. S. Trudinger, Elliptic Partial Differential Equations of Second order, 2nd edition, Springer-Verlag, 1983.





\bibitem{gui-cross} C. Gui, {\em Hamiltonian identities for elliptic partial differential equations}, J. Funct. Anal. 254 (2008), no. 4, 904--933


\bibitem{gulliver} R. Gulliver. {\em Index and total curvature of complete minimal surfaces.} Proc. Symp. Pure
Math., 44, 207"212, (1986).

\bibitem{Hau-Pacard}  L. Hauswirth and F. Pacard, {\em
Higher genus Riemann minimal surfaces}, Invent. Math. 169(2007), 569-620.


\bibitem{hm1} D. Hoffman and W.H. Meeks III, {\em A complete embedded minimal surface in $\R^3$ with genus one and three ends}, {J. Diff. Geom.} 21 (1985), 109--127.


\bibitem{hm2} D. Hoffman and W.H. Meeks III, {\em The asymptotic behavior of properly embedded minimal surfaces of finite topology}, {J. Am. Math. Soc. } 4(2) (1989), 667--681.

\bibitem{hm22}    D. Hoffman and W. H. Meeks III. {\em The strong halfspace theorem for minimal surfaces.}
Inventiones Math., 101, 373"-377, 1990.

\bibitem{hm3} D. Hoffman and W.H. Meeks III, {\em Embedded minimal surfaces of finite topology}, {Ann. Math.} 131 (1990), 1--34.

\bibitem{hk1} D. Hoffman and H. Karcher, {\em Complete embedded minimal surfaces of finite total curvature.} In Geometry V, Encyclopaedia Math. Sci., vol. 90, pp. 5-93, 262-272. Springer Berlin (1997).

 \bibitem{tonegawa} J.E. Hutchinson, Y. Tonegawa, {\em Convergence of phase interfaces in the van der Waals-Cahn-Hilliard theory}. Calc. Var. Partial Differential Equations 10 (2000), no. 1, 49--84.

\bibitem{jerison}
D. Jerison and R. Monneau, {\em Towards a counter-example to a conjecture of De Giorgi in high dimensions,} Ann. Mat. Pura Appl. 183 (2004), 439--467.



\bibitem{kap1} N. Kapuleas, {\em Complete embedded minimal surfaces of finite total curvature}, J. Differential Geom. 45 (1997), 95--169.



\bibitem{lopez-ros} F. J. L\'opez, A. Ros, {\em On embedded complete minimal surfaces of genus zero}, J. Differential Geom. 33 (1991), no. 1, 293--300.

%


%

\bibitem{modica1}
L. Modica,  {\em Convergence to minimal surfaces problem and global solutions of $\Delta u  = 2(u^3 -u)$.} Proceedings of the International Meeting on Recent Methods in Nonlinear Analysis (Rome, 1978), pp. 223--244, Pitagora, Bologna, (1979).



\bibitem{Mo} F. Morabito, {\em Index and nullity of the Gauss map of the Costa-Hoffman-Meeks surfaces}, preprint 2008.


\bibitem{Nay1} S. Nayatani, {\em Morse index and Gauss maps of complete minimal surfaces in Euclidean $3-$space}, {\em Comm. Math. Helv.}  68(4)(1993), 511--537.

\bibitem{Nay2} S. Nayatani, {\em Morse index  of complete minimal surfaces. In: Rassis, T.M. (ed.) The Problem of Plateau}, pp. 181--189(1992)

\bibitem{osserman} R. Osserman, {\em A survey of minimal surfaces}, Math. Studies 25, Van Nostrand, New York 1969.

\bibitem{pacard}
F. Pacard and M. Ritor\'e, {\em From the constant mean curvature hypersurfaces to
the gradient theory of phase transitions,} J. Differential Geom. 64 (2003), no. 3, 359--423.

\bibitem{perez-ros}
J. P\'erez, A. Ros, {\em The space of properly embedded minimal surfaces with finite total curvature.}
 Indiana Univ. Math. J. 45 (1996), no. 1, 177--204.


\bibitem{savin}
O. Savin, {\em Regularity of flat level sets in phase transitions}, Annals of Mathematics, 169 (2009), 41--78.

\bibitem{schoen 1} R. Schoen, {\em Uniqueness, symmetry, and embeddedness of minimal surfaces}, J. Differential Geom. 18 (1983), 791--809.

\bibitem{traizet} M.Traizet, {\em An embedded minimal surface with no symmetries}, J. Differential Geom. 60 (2002) 103--153.

















\end{thebibliography}
\end{document}